# The Lanczos algorithm for matrix functions

*a handbook for scientists*

Tyler Chen

## Abstract

Lanczos-based methods have become standard tools for tasks involving matrix functions. Progress on these algorithms has been driven by several largely disjoint communities, resulting many innovative and important advancements which would not have been possible otherwise. However, this also has resulted in a somewhat fragmented state of knowledge and the propagation of a number of incorrect beliefs about the behavior of Lanczos-based methods in finite precision arithmetic.

This monograph aims to provide an *accessible introduction* to Lanczos-based methods for matrix functions. The intended audience is scientists outside of numerical analysis, graduate students, and researchers wishing to begin work in this area. Our emphasis is on conceptual understanding, with the goal of providing a starting point to learn more about the remarkable behavior of the Lanczos algorithm. Hopefully readers will come away from this text with a better understanding of how to think about Lanczos for modern problems involving matrix functions, particularly in the context of finite precision arithmetic.

**Acknowledgements** I am grateful to Anne Greenbaum, Cameron Musco, and Christoper Musco for our many conversations on the Lanczos algorithm that greatly influenced my own understanding. Thanks also to Raphael Meyer for many suggestions that were incorporated in the second version.

# Contents

# 1 Introduction

The Lanczos algorithm [Lan50] is one of the most remarkable algorithms in numerical analysis. As with all Krylov subspace methods (KSMs), the Lanczos algorithm iteratively extracts information about a matrix through a sequence of matrix-vector products. However, in contrast to KSMs for non-symmetric problems, the Lanczos algorithm, which is designed specifically for symmetric matrices, utilizes a beautiful connection between the problem symmetry and orthogonal polynomials to avoid much of the expensive memory and arithmetic overheard incurred by general-purpose KSMs. Unfortunately, the theoretical elegance of the Lanczos algorithm is not without its practical costs—the behavior of Lanczos in finite precision arithmetic is very different than what one might expect from the exact arithmetic theory.

A great deal of effort has gone into obtaining and disseminating theory and wisdom about KSMs for linear systems and eigenvalue problems, the classical use cases for KSMs. As such, there are presently many excellent books and monographs on these topics [Gre97; Saa03; LS13] including includes several books which focus entirely on Lanczos-based methods, with a heavy emphasis on their behavior in finite precision arithmetic [CW12; Meu06; MS06]. Thus, while there are many interesting questions about KSMs for these classical tasks which remain open [CLS24], readers have plenty of options, at varying levels of technicality, to learn about these methods.

On the other hand, over the past several decades, KSMs have become increasingly standard tools for other tasks such as low-rank approximation, applying matrix functions, and approximating spectral densities. While these more modern tasks are in many ways closely related to eigenvalue problems and linear systems, the literature on the use of KSMs for these tasks remains somewhat fragmented. This has, in turn, lead to the unfortunate propagation of a number of incorrect beliefs about the behavior of KSMs, particularly the Lanczos algorithm in finite precision arithmetic.

Through a combination of theory and numerical examples, we aim to dispel with some of the most common misconceptions about Lanczos-based methods for matrix functions. The intended audience is scientists outside of numerical analysis, graduate students, and researchers wishing to work in this area.

Our emphasis is on conceptual understanding, with the goal of providing researchers from fields adjacent to numerical analysis a *starting point* to learn more about the behavior of Lanczos-based methods. Hopefully readers will come away from this text with a better understanding of how to *think* about Lanczos-based methods for modern problems involving matrix functions.

This monograph does not aim to be an authoritative reference on Lanczos-based methods. Indeed, while we do state and prove a number of useful results, more technical theorems, particularly those about algorithms run in finite precision arithmetic, are often summarized in an informal fashion so as to avoid obfuscating the *essence* of the theorem with the theorem itself.

## 1.1 Key quantities

Throughout, $\mathbf{A}$ will be a $d \times d$ real-symmetric matrix with eigendecomposition

$$\mathbf{A} = \sum_{i=1}^{d} \lambda_i \mathbf{u}_i \mathbf{u}_i^\mathsf{T}. \tag{1.1}$$

Here $\{\lambda_i\}$ are the set of (real) eigenvalues of $\mathbf{A}$ and $\{\mathbf{u}_i\}$ the corresponding orthonormal eigenvectors.

The star of this monograph is the matrix function.[1]

**Definition 1.1.** The matrix function $f(\mathbf{A})$, induced by $f : \mathbb{R} \to \mathbb{R}$ and $\mathbf{A}$, is defined as

$$f(\mathbf{A}) := \sum_{i=1}^{d} f(\lambda_i) \mathbf{u}_i \mathbf{u}_i^\mathsf{T}.$$

Common matrix functions including the matrix sign, logarithm, exponential, square root, and inverse square root, each of which has many applications throughout the mathematical sciences.

The Lanczos algorithm falls into a broader class of algorithms called Krylov subspace methods (KSMs), which make use of the information from a special subspace called the Krylov subspace.

**Definition 1.2.** The $k$-th Krylov subspace generated by the matrix $\mathbf{A}$ and vector $\mathbf{b}$ is defined as

$$\mathcal{K}_k(\mathbf{A}, \mathbf{b}) := \operatorname{span}\{\mathbf{b}, \mathbf{A}\mathbf{b}, \ldots, \mathbf{A}^{k-1}\mathbf{b}\}.$$

[1] Definition 1.1 is compatible with the standard definition of matrix polynomials. Indeed, one verifies $\mathbf{A}^k = \underbrace{\mathbf{A}\mathbf{A}\cdots\mathbf{A}}_{k \text{ times}}$, and linearity extends this observation to arbitrarily polynomials.

The $k$-th Krylov subspace consists of polynomials of $\mathbf{A}$, of degree less than $k$, applied to $\mathbf{b}$. Indeed, by definition, an arbitrary vector in the Krylov subspace looks like

$$c_0\mathbf{b} + c_1\mathbf{A}\mathbf{b} + \cdots + c_{k-1}\mathbf{A}^{k-1}\mathbf{b}, \qquad c_0, \ldots, c_{k-1} \in \mathbb{R}, \tag{1.2}$$

which is $p(\mathbf{A})\mathbf{b}$, where $p(x) = c_0 + c_1 x + \cdots + c_{k-1}x^{k-1}$. The Krylov subspace is also spanned by polynomials of this form, so we have an equivalent characterization

$$\mathcal{K}_k(\mathbf{A}, \mathbf{b}) = \{p(\mathbf{A})\mathbf{b} : \deg(p) < k\}. \tag{1.3}$$

This makes it clear that KSMs and polynomials are intimately related, and a major theme of this monograph is that this inherent connection should be leveraged.

Throughout, we will also make use of two important density functions. In both definitions, $\delta(x)$ is the Dirac delta unit point mass centered at zero.

**Definition 1.3.** The spectral density $\varphi(x; \mathbf{A})$ corresponding to $\mathbf{A}$ is defined as

$$\varphi(x) = \varphi(x; \mathbf{A}) := \sum_{i=1}^{d} d^{-1}\delta(x - \lambda_i).$$

**Definition 1.4.** The eigenvector spectral density $\psi(x; \mathbf{A}, \mathbf{b})$ corresponding to $\mathbf{A}$ and a vector $\mathbf{b}$ is defined as

$$\psi(x) = \psi(x; \mathbf{A}, \mathbf{b}) := \sum_{i=1}^{d} |\mathbf{u}_i^\mathsf{T}\mathbf{b}|^2\delta(x - \lambda_i).$$

The spectral density of $\mathbf{A}$ encodes the location and multiplicity of the eigenvalues of $\mathbf{A}$, while the eigenvector spectral density encodes the relationship between $\mathbf{A}$ and $\mathbf{b}$ observable to KSMs. In particular, for any orthogonal matrix $\mathbf{Q}$ it holds that

$$\psi(x; \mathbf{A}, \mathbf{b}) = \psi(x; \mathbf{Q}\mathbf{A}\mathbf{Q}^\mathsf{T}, \mathbf{Q}\mathbf{b}). \tag{1.4}$$

Thus, the eigenvector density function contains the information about $(\mathbf{A}, \mathbf{b})$ that is invariant to orthogonal transforms. In many settings this allows us to understand more carefully the underlying behavior of KSMs. When $\mathbf{b}$ is a unit-vector, then the eigenvector spectral density is unit-mass.

## 1.2 Problems of interest

This monograph is centered on tasks relating to matrix functions for which Lanczos-based methods are widely used. We now briefly introduce these problems.

**Linear systems** In Chapter 5, we discuss Lanczos-based methods for solving linear systems $\mathbf{A}\mathbf{x} = \mathbf{b}$. These include the well-known conjugate gradient and MINRES algorithms. The aim of the aim of this section is to highlight some concepts which are also relevant to understanding the other topics of interest.

**Action of matrix functions** In Chapter 6, we discuss methods for approximating $f(\mathbf{A})\mathbf{b}$. This can be viewed as a generalization of solving a linear system which corresponds to $f(x) = x^{-1}$.

**Quadratic forms and trace approximation** In Chapter 7, we discuss methods for approximating $\mathbf{b}^{\mathsf{T}} f(\mathbf{A})\mathbf{b}$. This is in many ways conceptually distinct from the task of computing the action of a matrix function. Through the use of stochastic trace estimation, methods for quadratic forms of matrix functions yield estimates for approximating $\operatorname{tr}(f(\mathbf{A}))$.

**Spectrum approximation** In Chapter 8, we discuss methods for obtaining coarse-grain approximations to the spectral density $\varphi(x; \mathbf{A})$. This is closely related to trace approximation.

## 1.3 Notation

Unless otherwise stated $\|\cdot\|$ will denote the Euclidean norm for vectors and the spectral norm for matrices. The spectrum of $\mathbf{A}$ is denoted by $\Lambda = \{\lambda_1, \ldots, \lambda_d\}$, the convex closure of the spectrum by $\mathcal{I} = [\lambda_{\min}, \lambda_{\max}]$, and the condition number of $\mathbf{A}$ by $\kappa(\mathbf{A})$. We will write $\Lambda(\mathbf{B})$ and $\mathcal{I}(\mathbf{B})$ to denote the analogous quantities for a matrix $\mathbf{B}$. Given a scalar function $f(x)$ and a set $S$, we define $\|f\|_S := \max_{x \in S} |f(x)|$; i.e. the infinity norm of a function $f$ over a set $S$. Finally, the $k$-th canonical basis vector is denoted by $\mathbf{e}_k$.

# 2 The Arnoldi and Lanczos algorithms

In this section, we describe the Arnoldi and Lanczos algorithms for generating a orthonormal basis $\mathbf{q}_0, \mathbf{q}_1, \ldots, \mathbf{q}_{k-1}$ for the Krylov subspace $\mathcal{K}_k(\mathbf{A}, \mathbf{b})$. First, however, we make a few simple observations about Krylov subspaces. Clearly Krylov subspaces are nested in the sense that

$$\mathcal{K}_0(\mathbf{A}, \mathbf{b}) \subseteq \mathcal{K}_1(\mathbf{A}, \mathbf{b}) \subseteq \cdots \subseteq \mathcal{K}_k(\mathbf{A}, \mathbf{b}). \tag{2.1}$$

Eventually (and certainly for some $k \leq d$) we must have that $\mathcal{K}_{k+1}(\mathbf{A}, \mathbf{b}) = \mathcal{K}_k(\mathbf{A}, \mathbf{b})$. At this point, the Krylov subspace stabilizes.

**Lemma 2.1.** Suppose $\mathcal{K}_{k+1}(\mathbf{A}, \mathbf{b}) = \mathcal{K}_k(\mathbf{A}, \mathbf{b})$. Then $\mathcal{K}_{k+j}(\mathbf{A}, \mathbf{b}) = \mathcal{K}_k(\mathbf{A}, \mathbf{b})$ for all $j > 0$.

*Proof.* If $\mathcal{K}_{k+1}(\mathbf{A}, \mathbf{b}) = \mathcal{K}_k(\mathbf{A}, \mathbf{b})$, then $\mathbf{A}^k\mathbf{b}$ can be written as a linear combination of $\mathbf{b}, \mathbf{A}\mathbf{b}, \ldots, \mathbf{A}^{k-1}\mathbf{b}$. But then $\mathbf{A}^{k+1}\mathbf{b} = \mathbf{A}(\mathbf{A}^k\mathbf{b})$ can be written as a linear combination of $\mathbf{A}\mathbf{b}, \ldots, \mathbf{A}^k\mathbf{b}$, and hence of $\mathbf{b}, \mathbf{A}\mathbf{b}, \ldots, \mathbf{A}^{k-1}\mathbf{b}$. The result then follows by repeating this argument. ■

The index at which the Krylov subspace stabilizes, sometimes called the grade, can be exactly characterized in terms of the eigenvector spectral density.

**Lemma 2.2.** Let $d'$ be the number of distinct points of support in $\psi(x; \mathbf{A}, \mathbf{b})$. Then $\dim(\mathcal{K}_k(\mathbf{A}, \mathbf{b})) = \min\{k, d'\}$.

*Proof.* Suppose $k \leq d'$. Then for any nonzero polynomial $p(x)$ of degree at most $k-1$, there exists an index $i$ such that $p(\lambda_i) \neq 0$ and $|\mathbf{u}_i^\mathsf{T}\mathbf{b}| > 0$. Hence $p(\mathbf{A})\mathbf{b} \neq \mathbf{0}$, so by (1.2), $\mathbf{b}, \mathbf{A}\mathbf{b}, \ldots, \mathbf{A}^{k-1}\mathbf{b}$ are linearly independent so $\dim(\mathcal{K}_k(\mathbf{A}, \mathbf{b})) = k$. Now, suppose $k = d' + 1$. Let $p(x)$ be a nonzero polynomial with roots at each of the $d'$ points of support of $\psi(x; \mathbf{A}, \mathbf{b})$. Then $p(\mathbf{A})\mathbf{b} = \mathbf{0}$, so $\mathbf{b}, \mathbf{A}\mathbf{b}, \ldots, \mathbf{A}^{d'}\mathbf{b}$ are linearly dependent. The result for $k > d' + 1$ follows by Lemma 2.1. ■

## 2.1 The Arnoldi algorithm

Perhaps the simplest approach to obtaining an orthonormal basis for the Krylov subspace is to construct a (non-orthogonal) basis and then orthogonalize it (e.g. using Gram–Schmidt). While this is fine in theory, in practice it can lead to numerical issues if the resulting basis vectors are not well-conditioned.

An alternative approach is to alternate between orthogonalizing and expanding the Krylov subspace. In particular, suppose we have already obtain an orthonormal basis $\mathbf{q}_0, \ldots, \mathbf{q}_{n-1}$ for $\mathcal{K}_n(\mathbf{A}, \mathbf{b})$. It is not too hard to verify that

$$\mathcal{K}_{n+1}(\mathbf{A}, \mathbf{b}) = \operatorname{span}\{\mathbf{q}_0, \ldots, \mathbf{q}_{n-1}, \mathbf{A}\mathbf{q}_{n-1}\}. \tag{2.2}$$

Thus, we can obtain an orthonormal basis for $\mathcal{K}_{n+1}(\mathbf{A}, \mathbf{b})$ by orthogonalizing $\mathbf{A}\mathbf{q}_{n-1}$ against the previous basis vectors (or identifying that $\mathbf{A}\mathbf{q}_{n-1}$ is in the span of $\mathbf{q}_0, \ldots, \mathbf{q}_{n-1}$, in which case the Krylov subspace has stabilized, and terminating). This results in the Arnoldi algorithm (Algorithm 2.3) [Arn51].

**Algorithm 2.3** (Arnoldi).

ARNOLDI($\mathbf{A}, \mathbf{b}, k$)
  $\mathbf{q}_0 = \mathbf{b}/\|\mathbf{b}\|$
  **for** $n = 0, 1, \ldots, k-1$ **do**
    $\mathbf{y}_{n+1} = \mathbf{A}\mathbf{q}_n$
    **for** $i = 0, \ldots, n$ **do**
      $h_{i,n} = \mathbf{q}_i^\mathsf{T}\mathbf{y}_{n+1}$
    $\mathbf{z}_{n+1} = \mathbf{y}_{n+1} - (h_{0,n}\mathbf{q}_0 + \cdots + h_{n,n}\mathbf{q}_n)$
    $h_{n+1,n} = \|\mathbf{z}_{n+1}\|$  ▷ terminate if $h_{n+1,n} = 0$
    $\mathbf{q}_{n+1} = \mathbf{z}_{n+1}/h_{n+1,n}$
  **return** $\{\mathbf{q}_n\}, \{h_{i,j}\}$

Observe that

$$h_{n+1,n}\mathbf{q}_{n+1} = \mathbf{z}_{n+1} = \mathbf{y}_{n+1} - (h_{0,n}\mathbf{q}_0 + \cdots + h_{n,n}\mathbf{q}_n), \qquad n \geq 0. \tag{2.3}$$

Thus, using that $\mathbf{y}_{n+1} = \mathbf{A}\mathbf{q}_n$ we obtain a recurrence

$$\mathbf{A}\mathbf{q}_n = h_{0,n}\mathbf{q}_0 + \cdots + h_{n,n}\mathbf{q}_n + h_{n+1,n}\mathbf{q}_{n+1}, \qquad n \geq 0. \tag{2.4}$$

Define

$$\mathbf{Q}_k := \begin{bmatrix} | & | & & | \\ \mathbf{q}_0 & \mathbf{q}_1 & \cdots & \mathbf{q}_{k-1} \\ | & | & & | \end{bmatrix}, \qquad \mathbf{H}_k := \begin{bmatrix} h_{0,0} & h_{0,1} & \cdots & h_{0,k-1} \\ h_{1,0} & h_{1,1} & \cdots & h_{1,k-1} \\ & \ddots & \ddots & \vdots \\ & & h_{k-1,k-2} & h_{k-1,k-1} \end{bmatrix}. \tag{2.5}$$

Then (2.4) can be written in matrix form as

$$\mathbf{A}\mathbf{Q}_k = \mathbf{Q}_k\mathbf{H}_k + h_{k,k-1}\mathbf{q}_k\mathbf{e}_k^\mathsf{T}, \tag{2.6}$$

where $\mathbf{e}_k$ is the $k$-th canonical basis vector.

Let $\mathbf{H}_{k+1,k}$ be the $(k+1)\times k$ matrix obtained by appending the $1\times k$ row-vector $h_{k,k-1}\mathbf{e}_k^\mathsf{T}$ to the bottom of $\mathbf{H}_k$. Then we can more compactly write (2.6) as

$$\mathbf{A}\mathbf{Q}_k = \mathbf{Q}_{k+1}\mathbf{H}_{k+1,k}. \tag{2.7}$$

While equivalent to (2.6), (2.7) is sometimes more convenient due to its compact form. We also note that, since $\mathbf{q}_0, \ldots, \mathbf{q}_k$ are orthonormal,

$$\mathbf{Q}_k^\mathsf{T}\mathbf{A}\mathbf{Q}_k = \mathbf{Q}_k^\mathsf{T}\mathbf{Q}_k\mathbf{H}_k + h_{k,k-1}\mathbf{Q}_k^\mathsf{T}\mathbf{q}_k\mathbf{e}_k^\mathsf{T} = \mathbf{H}_k. \tag{2.8}$$

## 2.2 From symmetry, Lanczos

Since $\mathbf{q}_0, \ldots, \mathbf{q}_k$ are orthonormal,

$$\mathbf{Q}_k^\mathsf{T}\mathbf{A}\mathbf{Q}_k = \mathbf{Q}_k^\mathsf{T}\mathbf{Q}_k\mathbf{H}_k + h_{k,k-1}\mathbf{Q}_k^\mathsf{T}\mathbf{q}_k\mathbf{e}_k^\mathsf{T} = \mathbf{H}_k. \tag{2.9}$$

Suppose $\mathbf{A}$ is symmetric; that is $\mathbf{A} = \mathbf{A}^\mathsf{T}$. Then, using (2.9), we see that

$$\mathbf{H}_k = \mathbf{Q}_k^\mathsf{T}\mathbf{A}\mathbf{Q}_k = \mathbf{Q}_k^\mathsf{T}\mathbf{A}^\mathsf{T}\mathbf{Q}_k = (\mathbf{Q}_k^\mathsf{T}\mathbf{A}\mathbf{Q}_k)^\mathsf{T} = (\mathbf{H}_k)^\mathsf{T}; \tag{2.10}$$

that is, $\mathbf{H}_k$ is also symmetric. By construction $\mathbf{H}_k$ is zero below the first sub-diagonal,[1] so $\mathbf{H}_k^\mathsf{T}$ is zero above the first super-diagonal. Therefore $\mathbf{H}_k$ must be tridiagonal!

This means that the majority of the coefficients produced by the Arnoldi algorithm are actually zero *a priori*. Skipping computing these coefficients (and taking advantage of the symmetry to save an additional inner product) results in the Lanczos algorithm (Algorithm 2.4) [Lan50]. While the Arnoldi and Lanczos algorithms are mathematically equivalent (i.e. produce the same quantities in exact arithmetic) when $\mathbf{A}$ is symmetric, the Lanczos algorithm requires significantly fewer arithmetic operations.

**Algorithm 2.4** (Lanczos).

LANCZOS($\mathbf{A}, \mathbf{b}, k$)
    $\mathbf{q}_0 = \mathbf{b}/\|\mathbf{b}\|$, $\beta_{-1} = 0$, $\mathbf{q}_{-1} = \mathbf{0}$
    **for** $n = 0, 1, \ldots, k-1$ **do**
        $\mathbf{y}_{n+1} = \mathbf{A}\mathbf{q}_n - \beta_{n-1}\mathbf{q}_{n-1}$
        $\alpha_n = \mathbf{q}_n^\mathsf{T}\mathbf{y}_{n+1}$

[1] This is called upper-Hessenberg.

6: $\mathbf{z}_{n+1} = \mathbf{y}_{n+1} - \alpha_n \mathbf{q}_n$
7: orthogonalize against $\mathbf{q}_0, \ldots, \mathbf{q}_n$ ▷ optional
8: $\beta_n = \|\mathbf{z}_{n+1}\|_2$ ▷ terminate if $\beta_n = 0$
9: $\mathbf{q}_{n+1} = \mathbf{z}_{n+1}/\beta_n$
10: **return** $\{\mathbf{q}_n\}, \{\alpha_n\}, \{\beta_n\}$.

The basis vectors produced by Lanczos are orthonormal and satisfy a symmetric three term recurrence

$$\mathbf{A}\mathbf{q}_n = \beta_{n-1}\mathbf{q}_{n-1} + \alpha_n \mathbf{q}_n + \beta_n \mathbf{q}_{n+1}, \qquad n \geq 0 \tag{2.11}$$

with initial conditions $\mathbf{q}_{-1} = \mathbf{0}$ and $\beta_{-1} = 0$. The coefficients $\{\alpha_n\}$ and $\{\beta_n\}$ defining the three term recurrence are also generated by the algorithm. This recurrence can be written in matrix form as

$$\mathbf{A}\mathbf{Q}_k = \mathbf{Q}_k \mathbf{T}_k + \beta_{k-1}\mathbf{q}_k \mathbf{e}_k^\mathsf{T} \tag{2.12}$$

where

$$\mathbf{Q}_k := \begin{bmatrix} | & | & & | \\ \mathbf{q}_0 & \mathbf{q}_1 & \cdots & \mathbf{q}_{k-1} \\ | & | & & | \end{bmatrix}, \quad \mathbf{T}_k := \begin{bmatrix} \alpha_0 & \beta_0 & & \\ \beta_0 & \alpha_1 & \ddots & \\ & \ddots & \ddots & \beta_{k-2} \\ & & \beta_{k-2} & \alpha_{k-1} \end{bmatrix}. \tag{2.13}$$

We might also write a compact matrix form

$$\mathbf{A}\mathbf{Q}_k = \mathbf{Q}_{k+1}\mathbf{T}_{k+1,k}, \tag{2.14}$$

where $\mathbf{T}_{k+1,k}$ is the $(k+1) \times k$ matrix defined by appending $\beta_{k-1}\mathbf{e}_k^\mathsf{T}$ below $\mathbf{T}_k$.

**Remark 2.5.** If the orthogonalization step in Line 7 of Algorithm 2.4 is omitted, the Lanczos algorithm can behave very differently in finite precision arithmetic than exact arithmetic.

For many applications, the observation in Remark 2.5 is not an issue; in fact, one of the main goals of this monograph is to provide insight into the behavior of Lanczos in finite precision arithmetic. We discuss the behavior of Lanczos in finite precision arithmetic in detail in Chapter 4, and we discuss the impacts on algorithms for relevant tasks as they arise.

## TL;DR

The Arnoldi and Lanczos algorithms iteratively produce an orthonormal basis for the Krylov subspace. The Lanczos algorithm avoids many of the inner products in the Arnoldi algorithm by takeing advantage of symmetry.

# 3 Orthogonal Polynomials

There is a fundamental relationship between orthogonal polynomials and the Lanczos algorithm which will be critical in our analysis of KSMs, particularly in finite precision arithmetic. We now recall some basic facts that will be relevant to our presentation, and for a more detailed treatment of the "beautiful mathematical relationships between matrices, moments, orthogonal polynomials, quadrature rules and the Lanczos [...] algorithm" we suggest [GM09].

Suppose $\mu(x)$ is a non-negative unit mass density function[1] supported on a subset of the real line with finite moments (i.e. $\int |x|^n \mu(x)\mathrm{d}x < \infty$ for all $n \geq 0$). The density function $\mu(x)$ induces an inner product $\langle \cdot, \cdot \rangle_\mu$ and a norm $\|\cdot\|_\mu$ defined by

$$\langle f, g \rangle_\mu := \int f(x)g(x)\mu(x)\mathrm{d}x, \qquad \|f\|_\mu := \langle f, f \rangle_\mu^{1/2}. \tag{3.1}$$

Associated with this inner product is an orthonormal sequence $p_0(x), p_1(x), \ldots$ of polynomials with $\deg(p_n) = n$ and positive leading coefficient. It is a remarkable fact that these polynomials satisfy the symmetric three-term recurrence,[2]

$$xp_n(x) = \beta_{n-1}p_{n-1}(x) + \alpha_n p_n(x) + \beta_n p_{n+1}(x), \qquad n \geq 0 \tag{3.2}$$

with initial conditions $p_0(x) = 1$, $p_{-1}(x) = 0$ and recurrence coefficients $\{\alpha_n\}$ and $\{\beta_n\}$ determined by $\mu(x)$.

## 3.1 Orthogonal polynomials and Lanczos

It is by no means a coincidence that the three-term recurrence (3.2) looks reminiscent of the Lanczos recurrence (2.11). Indeed, it is well-known that the Lanczos algorithm is closely related to the orthogonal polynomials of the eigenvector density function $\psi(x; \mathbf{A}, \mathbf{b})$ (Definition 1.4).

> **Theorem 3.1.** The recurrence coefficients of the orthogonal polynomials of $\psi(x; \mathbf{A}, \mathbf{b})$ are identical to the coefficients generated by Lanczos run on

[1] Here we are a bit loose with the term "function", and allow densities which include point masses; for example, the eigenvector density function $\psi(x; \mathbf{A}, \mathbf{b})$.

[2] This fact can be proved similarly to how we derive the Lanczos algorithm in Chapter 2.

$(\mathbf{A}, \mathbf{b})$, provided $\|\mathbf{b}\| = 1$. Moreover, the orthogonal polynomials $\{p_n(x)\}$ of $\psi(x; \mathbf{A}, \mathbf{b})$ are related to the Lanczos basis vectors $\{\mathbf{q}_n\}$ by $\mathbf{q}_n = p_n(\mathbf{A})\mathbf{b}$ for each $n \geq 0$.

*Proof.* This fact can (and should!) be verified by constructively computing the orthogonal polynomials of $\psi(x; \mathbf{A}, \mathbf{b})$ using the so-called Stieltjes algorithm (which is essentially just a continuous version of Lanczos) [Gau06]. Since $\langle f, g \rangle_\psi = \mathbf{b}^\mathsf{T} f(\mathbf{A})^\mathsf{T} g(\mathbf{A})\mathbf{b}$, this yields the Lanczos algorithm. ■

The fundamental relationship between the Lanczos algorithm on $(\mathbf{A}, \mathbf{b})$ and the orthogonal polynomials of $\psi(x; \mathbf{A}, \mathbf{b})$ is critical for understanding many aspects of the Lanczos algorithm, particularly its behavior in finite precision arithmetic. This is the main motivation for our study of orthogonal polynomials. In particular, in Section 3.2.1 we discuss a fundamental equivalence between matrices, moments, and quadrature (each of which will subsequently be defined).

## 3.2 Matrices, moments, and quadrature

In this section, we define three important mathematical quantities relating to $\mu(x)$: the Jacobi matrix, the polynomial moments, and the Gaussian quadrature rules. The title of this section pays homage to several papers and a book by Gene Golub and Gérard Meurant [GM93; GM97; GM09] which explore the relationships between these quantities in great detail.

The coefficients of the orthogonal polynomials of $\mu(x)$ will be important.

**Definition 3.2.** Define the symmetric tridiagonal matrix giving the first $2k - 1$ recurrence coefficients for the orthogonal polynomials of $\mu(x)$ by

$$\mathbf{M}_k = \mathbf{M}_k(\mu) := \begin{bmatrix} \alpha_0 & \beta_0 & & \\ \beta_0 & \ddots & \ddots & \\ & \ddots & \ddots & \beta_{k-2} \\ & & \beta_{k-2} & \alpha_{k-1} \end{bmatrix}. \tag{3.3}$$

Let $s$ be the smallest (possibly infinite) value for which $\beta_{s-1} = 0$. Then $\mathbf{M}_s(\mu)$ is called the Jacobi matrix corresponding to $\mu(x)$.

Another important property of a density function $\mu(x)$ is its polynomial moments.

**Definition 3.3.** Given a family of polynomials $\{q_n(x)\}$ with $\deg(q_n) = n$, the modified moments of $\mu(x)$ with respect to $\{q_n(x)\}$ are

$$m_n = m_n(\mu; \{q_n(x)\}) := \int q_n(x)\mu(x)\mathrm{d}x, \qquad n \geq 0. \tag{3.4}$$

Mathematically, the choice of polynomials $\{q_n(x)\}$ is not important (we can simply perform a change of basis), so in many theoretical settings it is common to use the monomials $q_n(x) = x^n$. However, in numerical settings the monomials are very poorly conditioned and it is almost always more appropriate to work with other families of polynomials, for instance Chebyshev polynomials [BB98; Gau04]. This will be particularly relevant in our discussion of the behavior of the Lanczos algorithm in finite precision arithmetic in Chapter 4.

Finally, we introduce the concept of Gaussian quadrature which allows us to integrate polynomials, and hence functions well-approximated by polynomials. This is the foundation for understanding methods involving quadratic forms Chapter 7 and will also be critical for understanding the behavior of Lanczos in finite precision arithmetic.

**Definition 3.4.** The $k$-point Gaussian quadrature rule $\mu_k(x)$ for $\mu(x)$ is defined by $\mu_k(x) := \psi(x; \mathbf{M}_k, \mathbf{e}_1)$, where $\mathbf{M}_k = \mathbf{M}_k(\mu)$ as in Definition 3.2.

From the definition of the eigenvector density Definition 1.4, we see the $k$-point Gaussian quadrature rule is supported on $k$ points, the eigenvalues of $\mathbf{M}_k$, and the corresponding weights are the squares of the first components of the eigenvectors of $\mathbf{M}_k$.

It is well-known that the Gaussian quadrature rule satisfies a moment-matching property, a proof of which we will provide in Section 3.3.2.

**Theorem 3.5.** For all polynomials $p(x)$ with $\deg(p) < 2k$, the $k$-point Gaussian quadrature rule $\mu_k(x)$ for $\mu(x)$ satisfies

$$\int p(x)\mu_k(x)\mathrm{d}x = \int p(x)\mu(x)\mathrm{d}x.$$

The following fact helps us understand how the cumulative distribution of the Gaussian quadrature relates to the cumulative distribution of the original density function.

**Theorem 3.6.** Let $\theta_1, \ldots, \theta_k$ be the support of $\mu_k(x)$ and denote by $M(x)$ and $M_k(x)$ the cumulative distribution functions of $\mu(x)$ and $\mu_k(x)$ respectively. Then, for each $i = 1, 2, \ldots, k$,

$$\lim_{\theta \to \theta_i^-} M_k(\theta) \leq M(\theta_i) \leq \lim_{\theta \to \theta_i^+} M_k(\theta).$$

In particular, since the Gaussian quadrature is supported on exactly $k$ points, the distribution function $M_k(x)$ is piecewise constant with jumps at exactly $k$ points. Theorem 3.6 asserts that at each of these jumps $M_k(x)$ goes from below $M(x)$ to above $M(x)$. This also implies that between consecutive jumps of $M_k(x)$, $M(x)$ goes from below $M_k(x)$ to above $M_k(x)$. This fact is visualized in Figure 3.1.

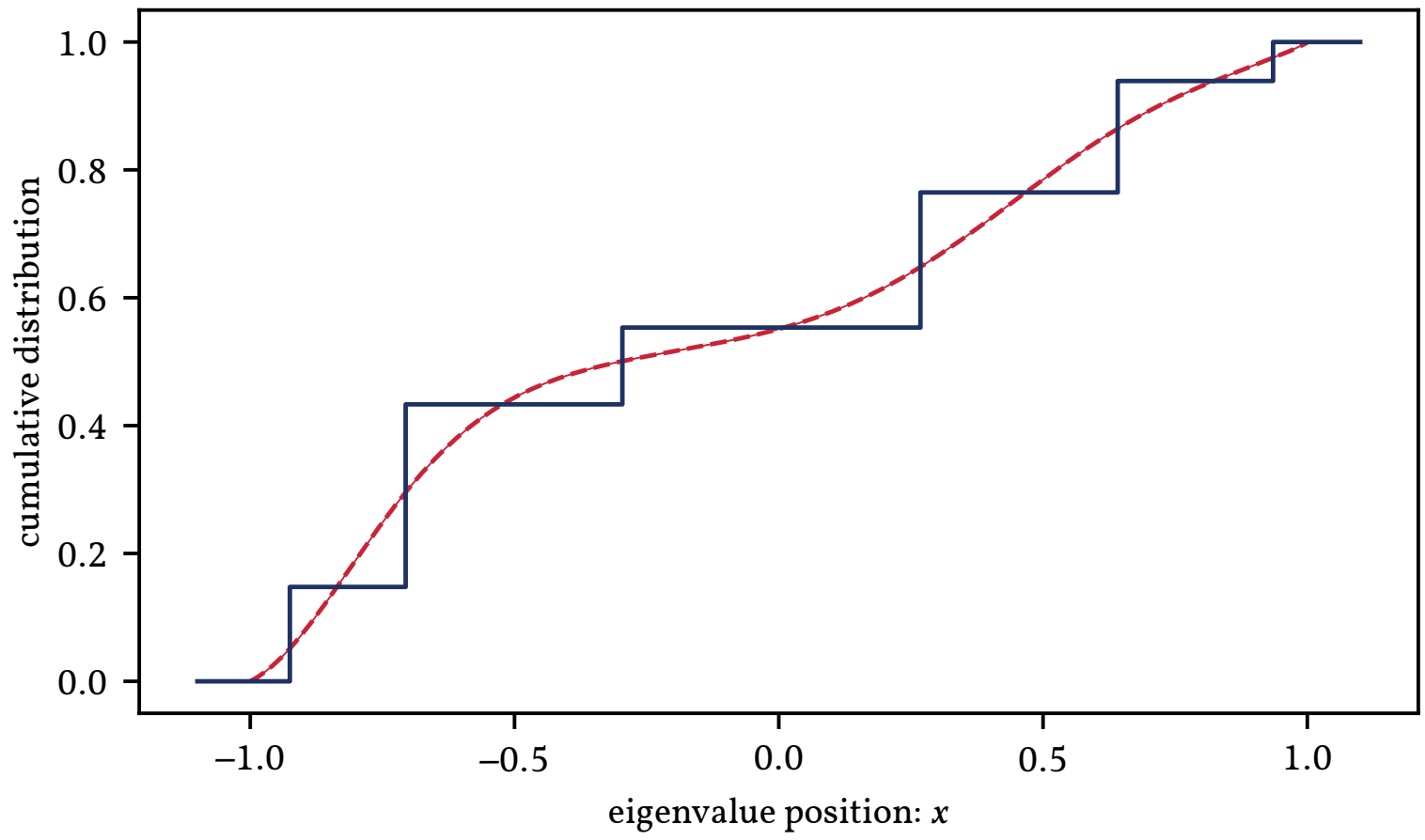

**Figure 3.1:** Cumulative distribution $M(x)$ ( ) and corresponding $k$-point Gaussian quadrature cumulative distribution $M_k(x)$ ( ) for $k = 6$. *Takeaway*: Observe that $M_k(x)$ jumps at $k$ points, and at each of these jumps goes from below $M(x)$ to above $M(x)$. Similarly, between consecutive jumps, $M(x)$ goes from below $M_k(x)$ to above $M_k(x)$.

Theorem 3.6 is an immediate consequence of the following property about the cumulative distributions of two density functions with matching moments; see for instance [KS72, Theorem 22.1].

**Theorem 3.7.** Suppose $\mu(x)$ and $\nu(x)$ are density functions with matching moments through degree $2k-1$. Let $M(x)$ and $N(x)$ be their respective cumulative distribution functions. Then $M(x) - N(x)$ is identically zero or changes signs at at least $2k-1$ points.

### 3.2.1 A fundamental equivalence

The three quantities we have just defined are fundamentally related; see e.g. [Gau04; GM09].

**Theorem 3.8.** There is a one-to-one correspondence between the following quantities:

- the upper-leftmost $k \times k$ submatrix of the Jacobi matrix for $\mu(x)$
- the (modified) moments of $\mu(x)$ through degree $2k - 1$
- the $k$-point Gaussian quadrature rule $\mu_k(x)$ for $\mu(x)$

This equivalence is critical to understanding many aspects of the Lanczos algorithm, particularly its behavior in finite precision arithmetic, and we will return to this equivalence throughout this monograph.

**Remark 3.9.** Exact knowledge of any one of the three quantities in Theorem 3.8 implies exact knowledge of the other two. However, the maps between these quantities may be very poorly conditioned [Gau04; OST07]. Therefore, small changes to one of the quantities (e.g. the moments) can lead to large changes in the others (e.g. the tridiagonal matrix).

## 3.3 Some properties of orthogonal polynomials

Orthogonal polynomials satisfy many remarkable properties, a few of which we recall here. For more detailed treatments we turn readers to [Gau04; GM09].

**Lemma 3.10.** Let $\pi_k(x) \propto p_k(x)$ be monic. Then $\|\pi_k(x)\|_\mu$ is minimal among all degree-$k$ monic polynomials.

*Proof.* An arbitrary degree-$k$ polynomial can be written as $p(x) = c_0 p_0(x) + \cdots + c_k p_k(x)$. By orthonormality of the $\{p_n\}$ we have $\|p\|_\mu^2 = |c_0|^2 + \cdots + |c_k|^2$. If $p(x)$ is monic the choice of $c_k$ is determined. All other coefficients should be set to zero. ■

**Lemma 3.11.** The roots of $p_k(x)$ are real, simple, and located in the convex closure of the support of $\mu(x)$.

*Proof.* Let $\{\theta_j\}_{j=0}^{k'-1}$ be the points in the convex closure of the support of $\mu(x)$ at which $p_k$ changes signs. Then,

$$\int p_k(x) \prod_{n=0}^{k'-1} (x - \theta_n)\mu(x)\mathrm{d}x \neq 0 \tag{3.5}$$

since the integrand does not change signs. This implies $k' = k$ since $p_k$ is orthogonal to all polynomials of lower degree. ■

**Lemma 3.12.** Let $\tau_1, \dots, \tau_{k-1}$ be the roots of $p_{k-1}(x)$ and $\theta_1, \dots, \theta_k$ the roots of $p_k(x)$. Then $\theta_1 < \tau_1 < \theta_2 < \cdots < \tau_{k-1} < \theta_k$.

**Lemma 3.13.** Let $(c, d)$ be an open interval for which $\mu(x) = 0$ for all $x \in (c, d)$. Then $p_k(x)$ has at most one zero in $(c, d)$.

*Proof.* Let $\pi_k(x) \propto p_k(x)$ be monic. Assume, for the sake of contradiction, there are two distinct zeros $\theta_i$ and $\theta_j$ in $(c, d)$. Write the remaining zeros as $\theta_n$. Now note

$$\int \pi_k(x) \prod_{n \neq i,j} (x - \theta_n)\mu(x)\mathrm{d}x = \int (x - \theta_i)(x - \theta_j) \prod_{n \neq i,j} (x - \theta_n)^2 \mu(x)\mathrm{d}x. \tag{3.6}$$

The left hand expression is zero because $\pi_k(x)$ is orthogonal to all polynomials of lower degree. However, the right hand side expression is nonzero because the integrand does not change signs on the support of $\mu(x)$. This is a contradiction. ■

**Lemma 3.14.** Suppose $\theta$ is a root of $p_k(x)$. Then

$$\left[p_0(\theta), p_1(\theta), \dots, p_{k-1}(\theta)\right]^{\mathsf{T}} \tag{3.7}$$

is an eigenvector of $\mathbf{M}_k$ with eigenvalue $\theta$.

*Proof.* In matrix form, (3.2) becomes

$$x\,\mathbf{P}_k(x) = \mathbf{P}_k(x)\mathbf{M}_k + \beta_{k-1}p_k(x)\mathbf{e}_k^{\mathsf{T}}. \tag{3.8}$$

By assumption, $p_k(\theta) = 0$. Therefore, evaluating each side of the above equality at $\theta$ and taking the transpose gives the result. ■

**Lemma 3.15.** Let $\theta_1, \dots, \theta_k$ be the roots of $p_k(x)$. For each $i = 1, \dots, k$ define $s_i := (\sum_{n=0}^{k-1} p_n(\theta_i)^2)^{1/2}$. Now, define the matrix

$$\mathbf{S}_k = \begin{bmatrix} p_0(\theta_1) & p_0(\theta_2) & \cdots & p_0(\theta_k) \\ p_1(\theta_1) & p_1(\theta_2) & \cdots & p_1(\theta_k) \\ \vdots & \vdots & & \vdots \\ p_{k-1}(\theta_1) & p_{k-1}(\theta_2) & \cdots & p_{k-1}(\theta_k) \end{bmatrix} \begin{bmatrix} 1/s_1 & & & \\ & 1/s_2 & & \\ & & \ddots & \\ & & & 1/s_k \end{bmatrix}.$$

Then $\mathbf{S}_k$ is orthogonal and $\mathbf{M}_k = \mathbf{S}_k \mathbf{\Theta}_k \mathbf{S}_k^\mathsf{T}$, where $\mathbf{\Theta}_k = \operatorname{diag}(\theta_1, \dots, \theta_k)$.

*Proof.* By Lemma 3.11, there are $k$ distinct roots of $p_k(x)$. Thus, Lemma 3.14 gives all eigenvectors of $\mathbf{M}_k$, and since $\mathbf{M}_k$ is symmetric, the eigenvectors are orthogonal. The result follows by normalizing these eigenvectors. ■

**Lemma 3.16.** It holds that $p_k(x) = \det(x\mathbf{I} - \mathbf{M}_k) / \prod_{n=0}^{k-1} \beta_n$.

*Proof.* By Lemma 3.14 we have that $p_k(x) \propto \det(x\mathbf{I} - \mathbf{M}_k)$. That the leading coefficient is $1/\prod_{n=0}^{k-1} \beta_n$ follows from the recurrence formula (3.2). ■

**3.3.1 Polynomial approximation and interpolation** We begin by noting that orthogonal projection onto the orthogonal polynomial basis produces the best approximation with respect to the $\mu$-norm.

**Lemma 3.17.** Let $p(x)$ be the polynomial minimizing $\|p - f\|_\mu$ among all polynomials of degree less than $k$. Then

$$p(x) = \begin{bmatrix} | & | & & | \\ p_0(x) & p_1(x) & \cdots & p_{k-1}(x) \\ | & | & & | \end{bmatrix} \begin{bmatrix} \langle f, p_0 \rangle_\mu \\ \langle f, p_1 \rangle_\mu \\ \vdots \\ \langle f, p_{k-1} \rangle_\mu \end{bmatrix}.$$

*Proof.* For any polynomial $q(x)$ of degree less than $k$,

$$\|f - p\|_\mu^2 \le \|f - p\|_\mu^2 + \|p - q\|_\mu^2 = \|(f - p) + (p - q)\|_\mu^2 = \|f - q\|_\mu^2. \tag{3.9}$$

Here we have used that $f(x) - p(x)$ is orthogonal to all polynomials of degree less than $k$. ■

Lemma 3.17 allows us to express low-degree polynomials in terms of $\mathbf{M}_k$.

**Lemma 3.18.** Suppose $p(x)$ is a polynomial of degree less than $k$. Then

$$p(x) = \begin{bmatrix} | & | & & | \\ p_0(x) & p_1(x) & \cdots & p_{k-1}(x) \\ | & | & & | \end{bmatrix} p(\mathbf{M}_k)\mathbf{e}_1.$$

*Proof.* By Lemma 3.15, we have that

$$p(\mathbf{M}_k)\mathbf{e}_1 = \mathbf{S}_k p(\mathbf{\Theta}_k)\mathbf{S}_k^\mathsf{T}\mathbf{e}_1.$$

Since $p_0(x) = 1$, a direct computation shows that

$$\mathbf{e}_m^\mathsf{T}\mathbf{S}_k p(\mathbf{\Theta}_k)\mathbf{S}_k^\mathsf{T}\mathbf{e}_1 = \sum_{\ell=1}^{k} p(\theta_\ell)p_{m-1}(\theta_\ell)/s_\ell^2 = \int p(x)p_{m-1}(x)\mu_k(x)\mathrm{d}x.$$

But now, since $p(x)$ and $p_{m-1}(x)$ are both polynomials of degree less than $k$, then $p(x)p_{m-1}(x)$ is a polynomial of degree at most $2k-1$. Hence, by Theorem 3.5, we have that

$$\int p(x)p_{m-1}(x)\mu_k(x)\mathrm{d}x = \int p(x)p_{m-1}(x)\mu(x)\mathrm{d}x = \langle p, p_{m-1}\rangle_\mu.$$

The result then follows from Lemma 3.17. ■

Finally, we state a characterization of the polynomial which interpolates $f(x)$ at the zeros of the degree-$k$ orthogonal polynomial $p_k(x)$ of $\mu(x)$.

**Lemma 3.19.** Let $p(x)$ be the polynomial interpolant to $f(x)$ of degree less than $k$ at the zeros of $p_k(x)$. Then

$$p(x) = \begin{bmatrix} | & | & & | \\ p_0(x) & p_1(x) & \cdots & p_{k-1}(x) \\ | & | & & | \end{bmatrix} f(\mathbf{M}_k)\mathbf{e}_1.$$

*Proof.* By Lemma 3.16, the zeros of $p_k(x)$ are the eigenvalues of $\mathbf{M}_k$. If $p(x) = f(x)$ at the eigenvalues of $\mathbf{M}_k$ then $p(\mathbf{M}_k) = f(\mathbf{M}_k)$. The result then follows from Lemma 3.18. ■

**3.3.2 Gaussian quadrature** We now have the tools required to prove Theorem 3.5.

*Proof of Theorem 3.5.* Let $\theta_1, \dots, \theta_k$ be the zeros of $p_k(x)$, and obtain weights $\omega_1, \dots, \omega_k$ by solving the linear system of equations

$$\begin{bmatrix} p_0(\theta_1) & p_0(\theta_2) & \cdots & p_0(\theta_k) \\ p_1(\theta_1) & p_1(\theta_2) & \cdots & p_1(\theta_k) \\ \vdots & \vdots & & \vdots \\ p_{k-1}(\theta_1) & p_{k-1}(\theta_2) & \cdots & p_{k-1}(\theta_k) \end{bmatrix} \begin{bmatrix} \omega_1 \\ \omega_2 \\ \vdots \\ \omega_k \end{bmatrix} = \begin{bmatrix} \int p_0(x)\mu(x)\,\mathrm{d}x \\ \int p_1(x)\mu(x)\,\mathrm{d}x \\ \vdots \\ \int p_{k-1}(x)\mu(x)\,\mathrm{d}x \end{bmatrix}. \tag{3.10}$$

Note that the right hand side is $\mathbf{e}_1$, since $p_n(x)$ is orthogonal to $p_0(x) = 1$ with respect to $\mu(x)$ for all $n \geq 1$. Denote by $\mathbf{P}_k^\mathsf{T}$ the coefficient matrix and recall Lemma 3.15 asserts the orthogonal eigenvector matrix $\mathbf{S}_k$ of $\mathbf{M}_k$ satisfies $\mathbf{S}_k = \mathbf{P}_k^\mathsf{T} \operatorname{diag}(1/s_1, \dots, 1/s_k)$ for some coefficients $\{s_i\}$. Then, since $p_0(x) = 1$, Lemma 3.15 implies

$$\boldsymbol{\omega} = (\mathbf{P}_k\mathbf{P}_k^\mathsf{T})^{-1}\mathbf{P}_k\mathbf{e}_1 = \operatorname{diag}(1/s_1^2, \dots, 1/s_k^2)\mathbf{1}.$$

Therefore, the $k$-point Gaussian quadrature rule satisfies

$$\mu_k(x) = \sum_{i=1}^{k} |\mathbf{e}_1^\mathsf{T}\mathbf{s}_i|^2\, \delta(x - \theta_i) = \sum_{i=1}^{k} \omega_i\, \delta(x - \theta_i). \tag{3.11}$$

Thus, we could have equivalently defined the Gaussian quadrature rule by its support and enforcing that it integrates polynomials of degree less than $k$ through (3.10).

Let $p(x)$ be an arbitrary polynomial of degree less than $2k$. By the Euclidean algorithm, we can decompose $p(x) = q(x)p_k(x) + r(x)$, where $\deg(q) < k$ and $\deg(r) < k$. Since $p_k(x)$ is an orthogonal polynomial of $\mu$, it is orthogonal to all polynomials of lower degree. Thus,

$$\int p(x)\mu(x)\mathrm{d}x = \int q(x)p_k(x)\mu(x)\mathrm{d}x + \int r(x)\mu(x)\mathrm{d}x = \int r(x)\mu(x)\mathrm{d}x. \tag{3.12}$$

On the other hand, since $p_k(\theta_i) = 0$ for each $i = 1, 2, \dots, k$,

$$\int p(x)\mu_k(x)\mathrm{d}x = \sum_{i=1}^{k} \omega_i p(\theta_i) = \sum_{i=1}^{k} \omega_i\big(q(\theta_i)p_k(\theta_i) + r(\theta_i)\big) = \sum_{i=1}^{k} \omega_i r(\theta_i). \tag{3.13}$$

Since $r(x)$ is degree at most $k-1$, it can be expressed as a linear combination of $p_0(x), \dots, p_{k-1}(x)$. Therefore, by the way we have obtained the weights $\omega_1, \dots, \omega_k$ in (3.10),

$$\int r(x)\mu(x)\mathrm{d}x = \sum_{i=1}^{k} \omega_i r(\theta_i). \tag{3.14}$$

Combining (3.12)–(3.14) gives the result. ■

## 3.4 Chebyshev polynomials

Two important orthogonal polynomial families are the Chebyshev polynomials of the first and second kind. These families of polynomials are respectively defined by the recurrences,

$$T_0(x) := 1, \quad T_1(x) := x, \quad T_n(x) := 2xT_{n-1}(x) - T_{n-2}(x), \quad n \geq 2 \tag{3.15}$$

$$U_0(x) := 1, \quad U_1(x) := 2x, \quad U_n(x) := 2xU_{n-1}(x) - U_{n-2}(x), \quad n \geq 2. \tag{3.16}$$

Both families of Chebyshev polynomials have explicit formulas. These formulas are often easier to work with than the recurrences (3.15) and (3.16).

**Lemma 3.20.** For all $k \geq 0$ and all $\theta \in \mathbb{R}$

$$T_k(\cos(\theta)) = \cos(k\theta), \qquad U_k(\cos(\theta))\sin(\theta) = \sin((k+1)\theta).$$

**Lemma 3.21.** For all $k \geq 0$ and $x \in \mathbb{R}$

$$T_k(x) = \left(\left(x + \sqrt{x^2-1}\right)^k + \left(x - \sqrt{x^2-1}\right)^k\right)\Big/2$$

$$U_k(x) = \left(\left(x + \sqrt{x^2-1}\right)^{k+1} - \left(x - \sqrt{x^2-1}\right)^{k+1}\right)\Big/\left(2\sqrt{x^2-1}\right).$$

Define the density functions $\mu_T(x)$ and $\mu_U(x)$, each supported on $[-1, 1]$, by

$$\mu_T(x) := \frac{1}{\pi}\frac{1}{\sqrt{1-x^2}}, \qquad \mu_U(x) := \frac{2}{\pi}\sqrt{1-x^2}. \tag{3.17}$$

Up to scaling, the Chebyshev polynomials are the orthogonal polynomials of the respective densities.

**Lemma 3.22.** For all $k \geq 0$,

$$\langle T_j, T_k \rangle_{\mu_T} = \begin{cases} 1 & j = k = 0 \\ 1/2 & j = k > 0 \\ 0 & j \neq k \end{cases}, \qquad \langle U_j, U_k \rangle_{\mu_U} = \begin{cases} 1 & j = k \\ 0 & j \neq k \end{cases}.$$

Our analysis repeatedly makes use of Chebyshev polynomial approximations to a function.

**Definition 3.23.** The degree-$k$ Chebyshev approximant to $f(x)$ is

$$p(x) = c_0 T_0(x) + 2\sum_{n=1}^{k} c_n T_n(x), \qquad c_n = \int f(x) T_n(x) \mu_T(x) \mathrm{d}x, \quad n \geq 0.$$

**Lemma 3.24.** The degree-$k$ Chebyshev approximant minimizes $\|f-p\|_{\mu_T}$ among all degree-$k$ polynomials.

*Proof.* This is an immediate consequence of Lemmas 3.17 and 3.22. ■

The Chebyshev polynomials satisfy many other properties. In fact, there are entire books just about Chebyshev polynomials [MH03; Riv20]! We now list a few of the properties that are most relevant to this monograph.

**Lemma 3.25.** Fix $k \geq 0$ and suppose $p(x)$ is a polynomial such that $\deg(p) \leq k$ and $\|p\|_{[-1,1]} \leq 1$. Then

$$\forall x \in \mathbb{R} \setminus (-1,1): \quad |T_k(x)| \geq |p(x)|.$$

In deriving theoretical guarantees for Lanczos in finite precision, we will repeatedly make use of the fact that Chebyshev polynomials are relatively small in the vicinity of $[-1,1]$.

**Lemma 3.26.** For all $k \geq 0$,

$$\|T_k\|_{[-1,1]} \leq 1, \qquad \|U_k\|_{[-1,1]} \leq k+1.$$

**Lemma 3.27.** For all $k \geq 1$, with $\eta_k = 1/(2k^2)$,

$$\|T_k\|_{[-1-\eta_k,1+\eta_k]} \leq 2, \qquad \|U_k\|_{[-1-\eta_k,1+\eta_k]} \leq 2(k+1).$$

*Proof.* A full proof is contained in [CT24, Lemma 3.7]. However, we can get the intuition for why this should be true relatively easily. Suppose $k$ is large and approximate

$$\left((1+1/(2k^2))^2-1\right)^{1/2} \approx \left(1+k^{-2}-1\right)^{1/2} \approx k^{-1}$$

so that

$$\left(1+1/(2k^2)+\left((1+1/(2k^2))^2-1\right)^{1/2}\right)^k \approx \left(1+k^{-1}\right)^k \approx \mathrm{e} \tag{3.18}$$

$$\left(1+1/(2k^2)-\left((1+1/(2k^2))^2-1\right)^{1/2}\right)^k \approx \left(1-k^{-1}\right)^k \approx \mathrm{e}^{-1}. \tag{3.19}$$

Then, by Lemma 3.21,

$$T_k(1+1/(2k^2)) \approx (\mathrm{e}+\mathrm{e}^{-1})/2 \approx 1.543 \leq 2. \tag{3.20}$$

A similar argument gives a similar result for $U_k(1+1/(2k^2))$. ■

## TL;DR

The Lanczos algorithm is, in effect, computing the orthogonal polynomials of the eigenvector density function. This, and the equivalence between moments, Jacobi matrices, and Gaussian quadrature will help us understand the Lanczos algorithm. Orthogonal polynomials, and in particular the Chebyshev polynomials, have many nice properties.

# 4 Lanczos in finite precision arithmetic

The computational savings of the Lanczos algorithm come from the fact that $\mathbf{A}\mathbf{q}_n$ is automatically orthogonal to $\mathbf{q}_0, \ldots, \mathbf{q}_{n-1}$, and hence does not need to be orthogonalized against these vectors. In finite precision arithmetic tiny rounding errors mean this is no longer the case, and if this orthogonalization is omitted, these small rounding errors can rapidly propagate.

We use an overline to denote quantities produced by the Lanczos algorithm in finite precision arithmetic. In particular, let $\overline{\mathbf{Q}}_k$ and $\overline{\mathbf{T}}_k$ be the outputs of the Lanczos algorithm run in finite precision arithmetic. The following are commonly observed effects of finite precision arithmetic. These effects are illustrated in Figure 4.1.

- **Loss of orthogonality:** The matrix $\overline{\mathbf{Q}}_k$ can be far from orthogonal, although its exact arithmetic counterpart $\mathbf{Q}_k$ is orthogonal. In fact, even linear independence may be lost.
- **Divergence of tridiagonal matrix:** The entries of the tridiagonal matrix $\overline{\mathbf{T}}_k$ produced in finite precision arithmetic may look nothing like their exact arithmetic counterparts from $\mathbf{T}_k$.
- **Ghost Ritz values:**[1] There may be multiple eigenvalues of $\overline{\mathbf{T}}_k$ near a single isolated eigenvalue of $\mathbf{A}$. This is not possible for $\mathbf{T}_k$, whose eigenvalues must interlace those of $\mathbf{A}$ (see Lemma 3.13).

Because the Lanzcos algorithm behaves drastically differently in finite precision arithmtic than exact arithmetic, there has been a widespread hesitance towards Lanczos-based approaches for many problems, at least without complicated/costly reorthogonalization schemes [JP94; Sil+96; Aic+03; Wei+06; UCS17; GWG19, etc.]. As we will discuss more explicitly in the relevant sections, this hesitance is often unfounded, and many Lanczos-based methods work fine, even in the presence of the phenomena described above.

The purpose of this section is to provide an overview of the rich theory known about the Lanczos algorithm in finite precision arithmetic. We begin, in

[1] The eigenvalues of $\mathbf{T}_k$ are commonly referred to as Ritz values.

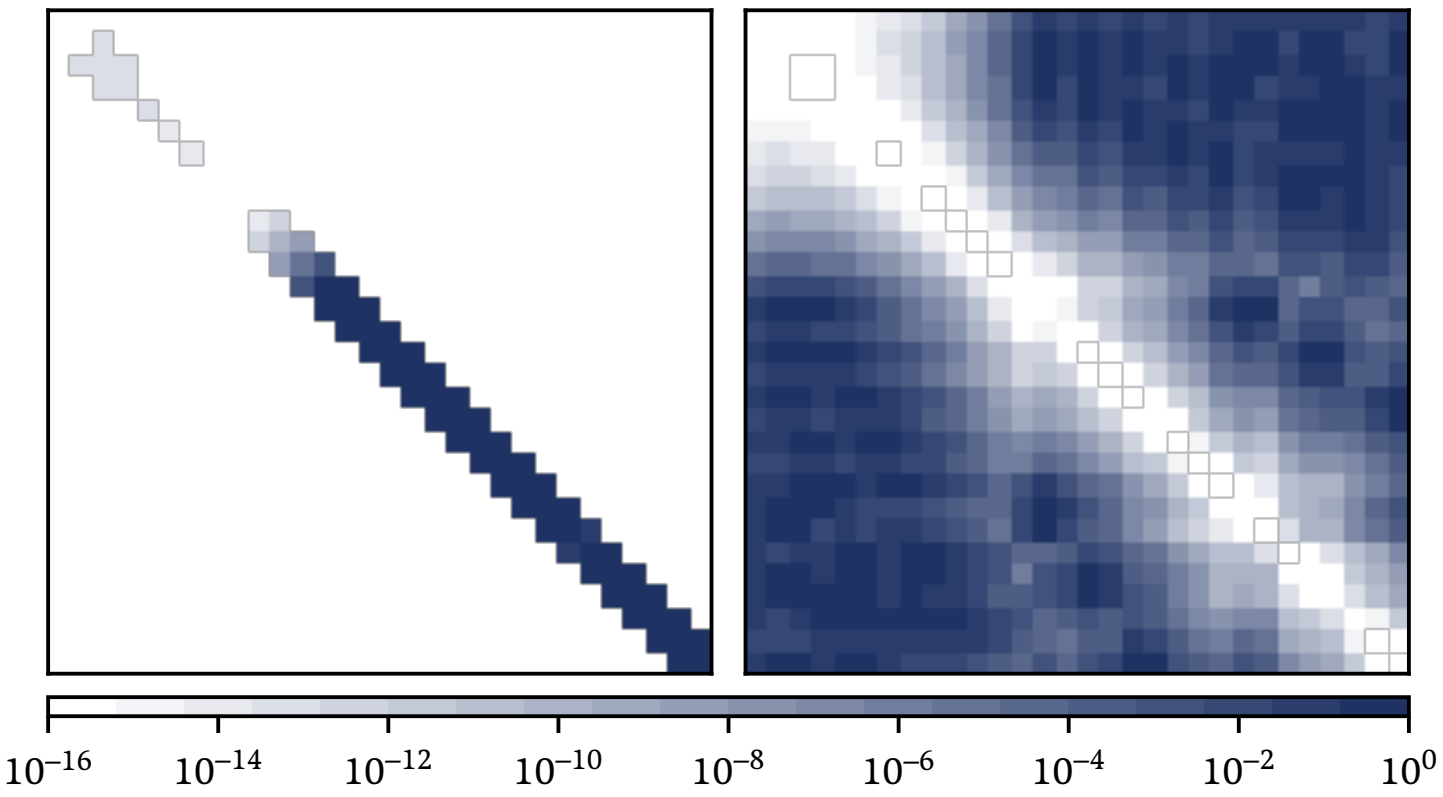

**Figure 4.1:** Instability of the Lanczos algorithm. *Left*: magnitude of entries of $\mathbf{T}_k - \overline{\mathbf{T}}_k$. *Right*: magnitude of entries of $\overline{\mathbf{Q}}_k^\top \overline{\mathbf{Q}}_k - \mathbf{I}$. *Takeaway*: The coefficients $\overline{\mathbf{T}}_k$ produced by the finite precision computation are far from what they would be in exact arithmetic, and there is a complete loss of orthogonality among the Lanczos vectors $\overline{\mathbf{Q}}_k$. In addition, while the top eigenvalues of $\mathbf{T}_k$ (and $\mathbf{A}$) are $1000.00, 588.17, 345.97, \ldots.$, the top eigenvalues of $\overline{\mathbf{T}}_k$ are $1000.00, 1000.00, 1000.00, 588.17, 588.17, \ldots$. There are multiple ghost eigenvalues of $\overline{\mathbf{T}}_k$ near the top eigenvalues of $\mathbf{A}$.

Section 4.1, by introducing the work of Chris Paige, which provides the foundation for most other analyses of the Lanczos algorithm in finite precision arithmetic. Then, in Section 4.2, we summarize the backwards stability analysis of Anne Greenbaum. More detailed treatments of the works of Paige and Greenbaum are readily available in textbooks and surveys [Par98; MS06; Meu06]. In Section 4.3 we discuss in detail a forward stability result of Leonid Knizhnerman. This result seems to have been mostly overlooked, even by the numerical analysis community, but elegantly explains why many Lanczos-based methods for matrix functions typically work fine in finite-precision arithmetic.

As with many results in numerical analysis, it is the intuition conveyed, rather than the actual bound itself, which is most valuable. Thus, we have tried to present the ideas of this section as clearly as possible, without getting bogged down in technical details. In particular, while we often rely on "$\approx$", all of the analyses are significantly more fine-grained, providing explicit dependencies on the dimension $d$, the number of iterations $k$, and the machine precision $\epsilon_{\text{mach}}$.

## 4.1 Paige's theory

Paige's work on the Lanczos algorithm is foundational; see e.g. [Pai70; Pai72; Pai76; Pai80]. Here we summarize only the aspects most relevant to methods for matrix functions. We refer readers to [Meu06] for a more detailed, but still accessible, treatment.

While errors can (and do) accumulate in the Lanczos algorithm, we expect *local* errors (the errors made at each iteration) to be small. In particular, the algorithm will still nearly satisfy a three term recurrence, and each Lanczos vector will be nearly orthogonal to the previous Lanczos vector and nearly unit length:

$$\mathbf{A}\overline{\mathbf{q}}_n \approx \overline{\beta}_{n-1}\overline{\mathbf{q}}_{n-1} + \overline{\alpha}_n\overline{\mathbf{q}}_n + \overline{\beta}_n\overline{\mathbf{q}}_{n+1}, \qquad \overline{\mathbf{q}}_n^\mathsf{T}\overline{\mathbf{q}}_{n+1} \approx 0, \qquad \|\overline{\mathbf{q}}_n\| \approx 1. \tag{4.1}$$

It will be useful to write

$$\mathbf{A}\overline{\mathbf{Q}}_k = \overline{\mathbf{Q}}_k\overline{\mathbf{T}}_k + \overline{\beta}_{k-1}\overline{\mathbf{q}}_k\mathbf{e}_k^\mathsf{T} + \mathbf{F}_k, \tag{4.2}$$

where the perturbation term $\mathbf{F}_k$ accounts for local rounding errors made by the algorithm. Paige's analysis implies $\mathbf{F}_k \approx \mathbf{0}$. In addition, Paige's analysis shows that,

$$\Lambda(\overline{\mathbf{T}}_k) \subset [\lambda_{\min}(\mathbf{A}) - \eta_k, \lambda_{\max}(\mathbf{A}) + \eta_k], \qquad \eta_k \approx 0. \tag{4.3}$$

We are unfortunately unaware of any intuition for why one might expect $\eta_k$ to be small, and a simple explanation would be of great interest to the author.

As it relates to the analyses in this monograph, the most important results of Paige are the following.

**Theorem 4.1** (informal; see [Pai70; Pai72; Pai76; Pai80]). Treat the input $(\mathbf{A}, \mathbf{b})$ and hence the dimension $d$ as a constant. Suppose the Lanczos algorithm is run for $k$ iterations on a computer with relative machine precision $\epsilon_{\text{mach}} = O(\text{poly}(k)^{-1})$. Then, the approximate equalities in (4.2) and (4.3) hold up to additive error $O(\text{poly}(k)\,\epsilon_{\text{mach}})$.

**Remark 4.2.** The error in the approximate equalities in (4.2) and (4.3) can be computed directly after the Lanczos algorithm has been run (e.g. using higher precision arithmetic).

Paige's work extends beyond what is summarized above. In particular, it reveals that a loss of orthogonality in $\overline{\mathbf{Q}}_k$ coincides with the convergence of a Ritz value to an eigenvalue of $\mathbf{A}$. This analysis motivated a number of schemes, such as selective reorthogonalization [PS79] and partial reorthogonalization [Sim84] , which aim to maintain orthogonality in the Lanczos basis $\overline{\mathbf{Q}}_k$ without orthogonalizing against all previous Lanczos vectors (full reorthogonalization). Such methods are beyond the scope of our discussion.

## 4.2 Greenbaum's theory

Greenbaum proved a backwards stability[2] result for the Lanczos algorithm.

> **Theorem 4.3** (informal; see [Gre89]). Under the assumptions of Paige's analysis, suppose that Lanczos is run in finite precision arithmetic on $(\mathbf{A}, \mathbf{b})$ for $k$ iterations to produce $\overline{\mathbf{T}}_k$. There is a matrix $\tilde{\mathbf{A}}$ and vector $\tilde{\mathbf{b}}$ so that the (exact) Lanczos algorithm run on $(\tilde{\mathbf{A}}, \tilde{\mathbf{b}})$ for $k$ iterations produces $\overline{\mathbf{T}}_k$ and $\psi(x; \mathbf{A}, \mathbf{b}) \approx \psi(x; \tilde{\mathbf{A}}, \tilde{\mathbf{b}})$ in the following sense:
>
> - Let $\tilde{\lambda}$ be an eigenvalue of $\tilde{\mathbf{A}}$. Then there exists an eigenvalue $\lambda$ of $\mathbf{A}$ so that $\tilde{\lambda} \approx \lambda$.
> - Let $\lambda_i$ be an eigenvalue of $\mathbf{A}$. Define $S_i = \{j : |\tilde{\lambda}_j - \lambda_i| = \min_\ell |\tilde{\lambda}_j - \lambda_\ell|\}$ as the indices the eigenvalues of $\tilde{\mathbf{A}}$ nearer to $\lambda_i$ than any other eigenvalue of $\mathbf{A}$. Then $\sum_{j \in S_i} |\tilde{\mathbf{b}}^\mathsf{T} \tilde{\mathbf{u}}_j|^2 \approx |\mathbf{b}^\mathsf{T} \mathbf{u}_i|^2$.

Note that here the "nearby problem" $(\tilde{\mathbf{A}}, \tilde{\mathbf{b}})$ is of a different dimension than $(\mathbf{A}, \mathbf{b})$. Thus, the nearness of $(\tilde{\mathbf{A}}, \tilde{\mathbf{b}})$ to $(\mathbf{A}, \mathbf{b})$ must be interpreted in terms of the eigenvector densities $\psi(x; \tilde{\mathbf{A}}, \tilde{\mathbf{b}})$ and $\psi(x; \mathbf{A}, \mathbf{b})$. In particular, Theorem 4.3 asserts that the support of $\psi(x; \tilde{\mathbf{A}}, \tilde{\mathbf{b}})$ is contained in tiny intervals about the support of $\psi(x; \mathbf{A}, \mathbf{b})$, and that the mass of $\psi(x; \mathbf{A}, \mathbf{b})$ on each eigenvalue of $\mathbf{A}$ is approximated by the mass of $\psi(x; \tilde{\mathbf{A}}, \tilde{\mathbf{b}})$ on the corresponding eigenvalues of $\tilde{\mathbf{A}}$.

An example of the what $\psi(x; \tilde{\mathbf{A}}, \tilde{\mathbf{b}})$ could look like in relation to $\psi(x; \mathbf{A}, \mathbf{b})$ is illustrated in Figure 4.2. Note the presence of clusters of eigenvalues of $\tilde{\mathbf{A}}$ around those of $\mathbf{A}$, and that the total mass of the eigenvalue cluster of $\tilde{\mathbf{A}}$ matches closely that of the corresponding eigenvalue of $\mathbf{A}$. The presence of clusters of eigenvalues of $\tilde{\mathbf{A}}$ explains the Ghost Ritz values of $\overline{\mathbf{T}}_k$.

**A smoothed analysis?** Greenbaum's theory guarantees the existence of a nearby problem $(\tilde{\mathbf{A}}, \tilde{\mathbf{b}})$ for which the tridiagonal matrix output by (exact arithmetic) Lanczos on $(\tilde{\mathbf{A}}, \tilde{\mathbf{b}})$ *exactly* matches the output $\overline{\mathbf{T}}_k$ of finite precision arithmetic on the original problem $(\mathbf{A}, \mathbf{b})$. One could construct a problem $(\hat{\mathbf{A}}, \hat{\mathbf{b}})$ qualatatively similar to $(\tilde{\mathbf{A}}, \tilde{\mathbf{b}})$; i.e with the eigenvalues of $\hat{\mathbf{A}}$ clustered around those of $\mathbf{A}$ and such that the mass on each eigenvalue of $\mathbf{A}$ is approximated by the mass on the corresponding eigenvalues of $\hat{\mathbf{A}}$.[3] Interestingly, the tridiagonal matrix output by Lanczos on $(\hat{\mathbf{A}}, \hat{\mathbf{b}})$ is often qualatatively very similar to $\overline{\mathbf{T}}_k$. This was explored in [GS92b]; see also [CLS24].

[2] In numerical analysis, an algorithm is *backwards stable* if, for any given input, the finite precision output is exactly equal to what the algorithm would produce (in exact arithmetic) on some slightly perturbed input.

[3] This can be thought of as convolving $\psi(x; \mathbf{A}, \mathbf{b})$ with some kernel with support contained in a tiny interval about zero.

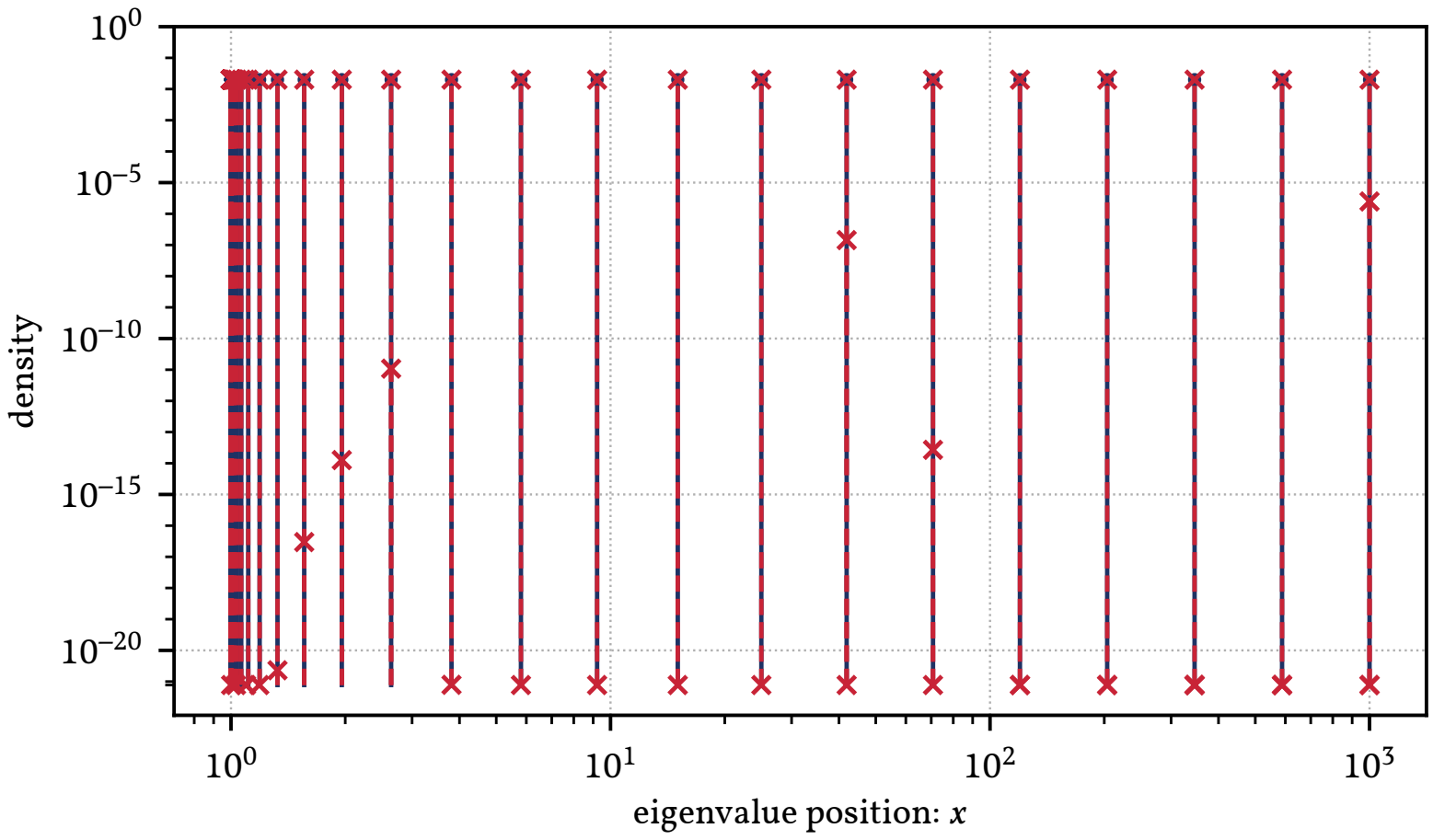

**Figure 4.2:** Illustration of the backward stability guarantee of Greenbaum. $\psi(x; \mathbf{A}, \mathbf{b})$ ( —●— ), and $\psi(x; \tilde{\mathbf{A}}, \tilde{\mathbf{b}})$ ( —×— ), where Lanczos run on $(\tilde{\mathbf{A}}, \tilde{\mathbf{b}})$ produces $\overline{\mathbf{T}}_k$, the exact quantity output by Lanczos run on $(\mathbf{A}, \mathbf{b})$ in finite precision arithmetic. *Takeaway*: $\tilde{\mathbf{A}}$ has more eigenvalues than $\mathbf{A}$, but these eigenvalues are closely clustered about those of $\mathbf{A}$. The backwards stability result of Greenbaum must be interpreted in terms of $\psi(x; \tilde{\mathbf{A}}, \tilde{\mathbf{b}})$ and $\psi(x; \mathbf{A}, \mathbf{b})$.

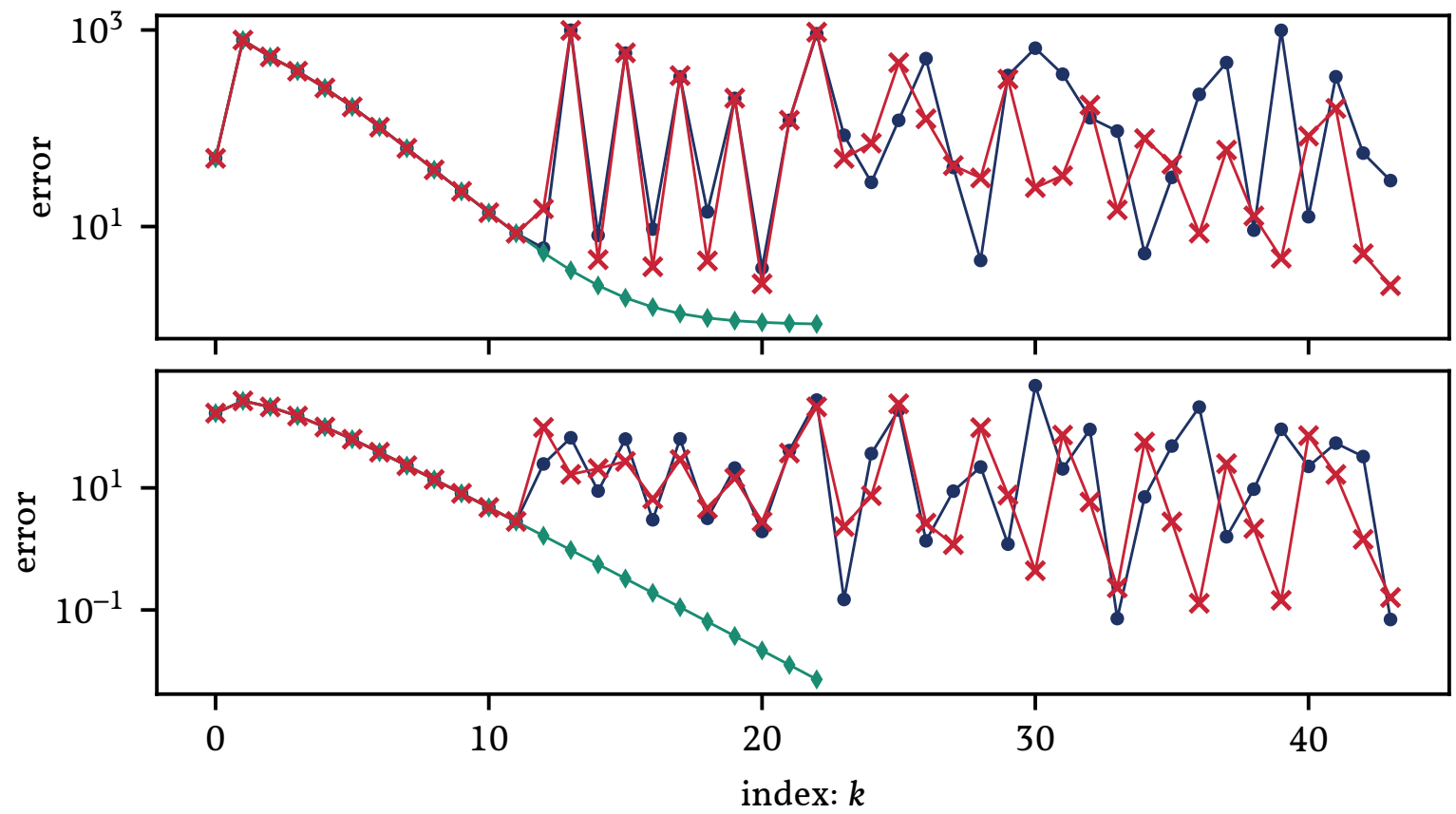

**Figure 4.3:** Diagonal (top) and off-diagonal (bottom) Lanczos coefficients with ( —×— ) and without ( —●— ) reorthgonalization, and on simulated a problem $(\hat{\mathbf{A}}, \hat{\mathbf{b}})$ with reorthgonalization ( —♦— ). *Takeaway*: Qualitative properties of $\overline{\mathbf{T}}_k$ can be realized in exact arithmetic by constructing a problem with eigenvalues clustered about those of the original problem.

We illustrate this phenomenon in Figure 4.3, where we show the entries of the tridiagonal matrices generated by Lanczos on $(\mathbf{A}, \mathbf{b})$ in exact and finite precision arithmetic, as well as those generated by Lanczos on $(\hat{\mathbf{A}}, \hat{\mathbf{b}})$ in exact arithmetic. Here $\hat{\mathbf{A}}$ is obtained by replacing each eigenvalue of $\mathbf{A}$ with 10 eigenvalue equally spaced throughout an interval of width $1.2 \times 10^{-13}$ centered at the original eigenvalue.

## 4.3 Knizhnerman's theory

Knizhnerman's work [Kni96], which unfortunately seems to have been mostly overlooked in the literature, asserts that the Chebyshev moments of the output of the Lanczos algorithm in finite precision arithmetic are close to the true Chebyshev moments. This remarkable fact is illustrated in Figure 4.4.

**Theorem 4.4** (informal; see [Kni96])**.** Suppose $\|\mathbf{A}\| = \|\mathbf{b}\| = 1$ and that Lanczos algorithm is run for $k$ iterations on a computer with relative machine precision $\epsilon_{\text{mach}} = O(\text{poly}(k)^{-1})$. Let $\psi(x) = \psi(x; \mathbf{A}, \mathbf{b})$ and $\overline{\psi}_k(x) = \psi(x; \overline{\mathbf{T}}_k, \mathbf{e}_1)$ be the Gaussian quadrature produced by Lanczos in finite precision arithmetic. Then, for all $n \leq 2k - 1$,

$$\left| \int T_n(x)\psi(x)\mathrm{d}x - \int T_n(x)\overline{\psi}_k(x)\mathrm{d}x \right| = O(\text{poly}(k)\,\epsilon_{\text{mach}}).$$

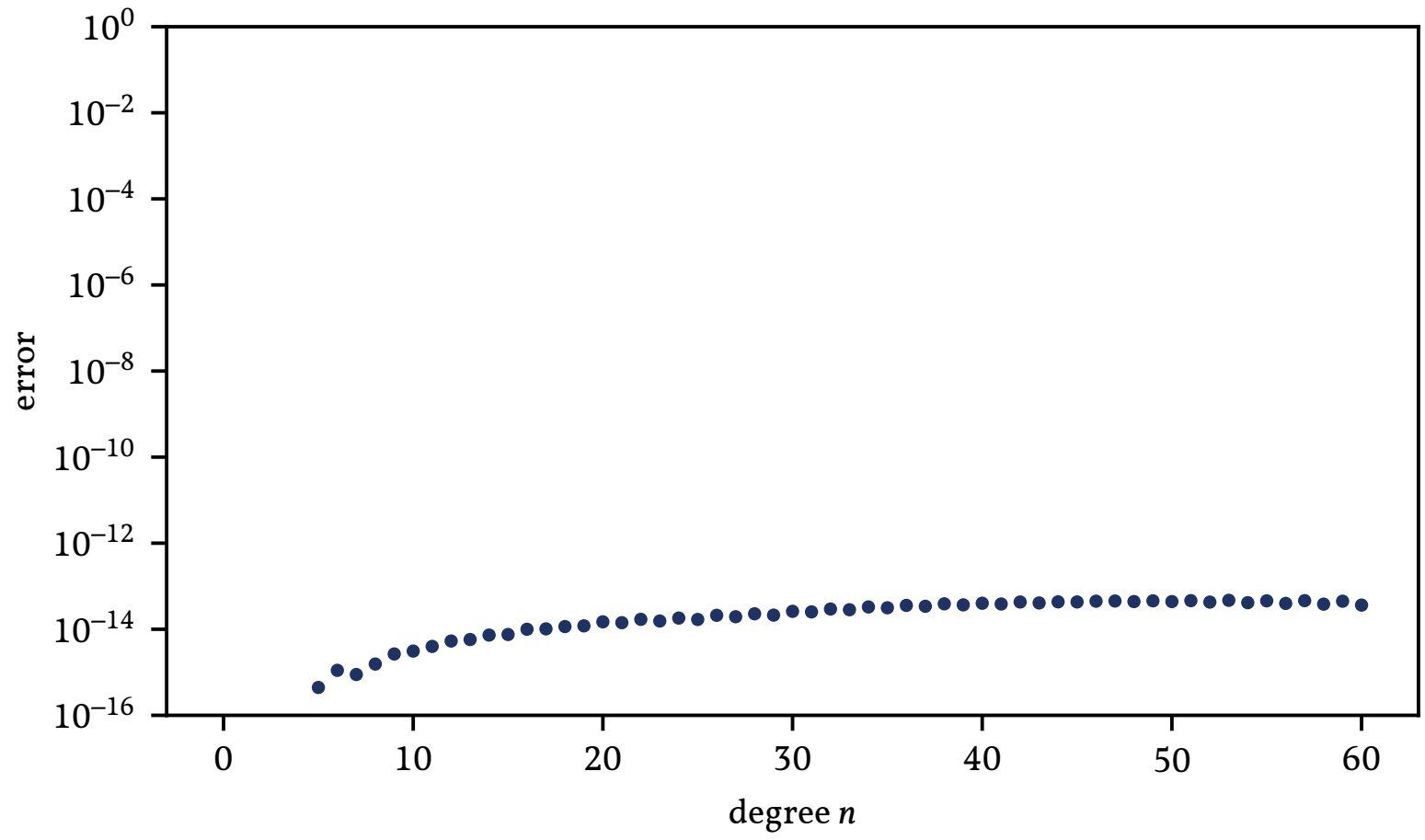

**Figure 4.4:** Error of the Chebyshev moments of $\psi_k(x)$ and $\overline{\psi}_k(x)$ for the same example as Figure 4.1. *Takeaway*: The Chebyshev moments of $\psi_k(x)$ and $\overline{\psi}_k(x)$ differ only on the order of the machine precision, even though the recurrence coefficients $\mathbf{T}_k$ and $\overline{\mathbf{T}}_k$ differ greatly. In this sense, the Lanczos algorithm is forward stable!

**Remark 4.5.** Despite many claims that Chebyshev-based methods are "more stable" than Lanczos-based methods, the precise version of Theorem 4.4 bounds is almost identical to what can be obtained for a method based directly on the Chebyshev recurrence (3.15).

In light of the fundamental equivalence between the moments of $\psi(x; \mathbf{A}, \mathbf{b})(x)$ and the tridiagonal matrix $\mathbf{T}_k$ output by the Lanczos algorithm, which we described in Section 3.2.1, Theorem 4.4 is a *forward stability* result for the Lanczos algorithm. As such, one might view the fact that $\overline{\mathbf{T}}_k$ is far from $\mathbf{T}_k$ as a result of poor conditioning of the map from moments to tridiagonal matrix (see Remark 3.9), rather than an instability of the Lanczos algorithm.

## TL;DR

The Lanczos algorithm behaves very differently in finite precision arithmetic than in exact arithmetic. However, much is known about how the algorithm behaves. In particular, the results of Greenbaum and Knizhnerman, based on the work of Paige, show that, in a certain sense, the Lanczos algorithm is actually forward/backward stable.

# 5 Linear Systems and least squares problems

One of the most important tasks in all of the computational sciences is solving the linear system of equations

$$\mathbf{A}\mathbf{x} = \mathbf{b}; \tag{5.1}$$

i.e. approximating $\mathbf{A}^{-1}\mathbf{b}$. Solving linear systems is also one of the main uses of KSMs, and there are many textbooks covering this topic. Our presentation differs from most textbooks in that we do not even write down the typical iteration for the conjugate gradient or MINRES algorithms. Instead, we define the algorithms using the quantities generated by the Lanczos algorithm, and focus primarily on developing key ideas which will be important in Chapter 6 when discussing more general matrix functions.

## 5.1 Positive definite systems: conjugate gradient

Throughout this subsection, we assume that $\mathbf{A}$ is positive definite; that is, that every eigenvalue of $\mathbf{A}$ is strictly positive. In this case, we can define the $\mathbf{A}$-norm $\|\cdot\|_{\mathbf{A}}$ by $\|\mathbf{x}\|_{\mathbf{A}} = (\mathbf{x}^{\mathsf{T}}\mathbf{A}\mathbf{x})^{1/2} = \|\mathbf{A}^{1/2}\mathbf{x}\|$.

**Definition 5.1.** The $k$-th conjugate gradient iterate $\mathbf{x}_k^{\mathrm{CG}}$ is given by

$$\mathbf{x}_k^{\mathrm{CG}} = \mathbf{x}_k^{\mathrm{CG}}(\mathbf{A},\mathbf{b}) := \|\mathbf{b}\|\mathbf{Q}_k\mathbf{T}_k^{-1}\mathbf{e}_1.$$

This definition of the CG iterate differs from the more typical implementation [HS52]. We explain the connection between Definition 5.1 and the more standard implementation in Section 5.1.3. However, to see that the quantity we call CG is indeed equivalent to the more standard definition, it suffices to verify the well-known $\mathbf{A}$-norm optimality guarantee.

**Theorem 5.2.** If $\mathbf{A}$ is positive definite, the CG iterate satisfies the formula

$$\mathbf{x}_k^{\mathrm{CG}} = \operatorname*{argmin}_{\mathbf{x}\in\mathcal{K}_k(\mathbf{A},\mathbf{b})} \|\mathbf{A}^{-1}\mathbf{b} - \mathbf{x}\|_{\mathbf{A}}.$$

*Proof.* Since $\mathbf{Q}_k$ is a basis for $\mathcal{K}_k(\mathbf{A},\mathbf{b})$, any vector $\mathbf{x} \in \mathcal{K}_k(\mathbf{A},\mathbf{b})$ can be written as $\mathbf{x} = \mathbf{Q}_k\mathbf{c}$ for some $\mathbf{c} \in \mathbb{R}^k$. Using this and the definition of the $\mathbf{A}$-norm we have

$$\underset{\mathbf{x}\in\mathcal{K}_k(\mathbf{A},\mathbf{b})}{\operatorname{argmin}} \|\mathbf{A}^{-1}\mathbf{b}-\mathbf{x}\|_{\mathbf{A}} = \mathbf{Q}_k \underset{\mathbf{c}\in\mathbb{R}^k}{\operatorname{argmin}} \|\mathbf{A}^{-1/2}\mathbf{b}-\mathbf{A}^{1/2}\mathbf{Q}_k\mathbf{c}\|. \tag{5.2}$$

The solution to this least squares problem is

$$\mathbf{Q}_k(\mathbf{Q}_k^\mathsf{T}\mathbf{A}^{1/2}\mathbf{A}^{1/2}\mathbf{Q}_k)^{-1}\mathbf{Q}_k^\mathsf{T}\mathbf{A}^{1/2}\mathbf{A}^{-1/2}\mathbf{b} = \|\mathbf{b}\|\mathbf{Q}_k\mathbf{T}_k^{-1}\mathbf{e}_1 = \mathbf{x}_k^{\mathrm{CG}}. \tag{5.3}$$

Here we have used that $\mathbf{Q}_k^\mathsf{T}\mathbf{A}\mathbf{Q}_k = \mathbf{T}_k$ and that $\mathbf{Q}_k^\mathsf{T}\mathbf{b} = \|\mathbf{b}\|\mathbf{e}_1$. ■

**5.1.1 Error bounds and estimates** We now review several standard error bounds and estimates. The main goal is to highlight the type results which are known for approximating $\mathbf{A}^{-1}\mathbf{b}$ in order to serve as reference for our discussion on algorithms for $f(\mathbf{A})\mathbf{b}$ in Chapter 6. Recall that $\Lambda$ is the spectrum of $\mathbf{A}$ and $\mathcal{I} = [\lambda_{\min}, \lambda_{\max}]$ is the smallest interval containing $\Lambda$.

**A priori error bounds** The optimality of the CG iterate allows the derivation of a priori error bounds in terms of spectral quantities of $\mathbf{A}$ such as the location of eigenvalues and condition number.

**Corollary 5.3.** The CG iterate satisfies the error bound

$$\frac{\|\mathbf{A}^{-1}\mathbf{b}-\mathbf{x}_k^{\mathrm{CG}}\|_{\mathbf{A}}}{\|\mathbf{A}^{-1}\mathbf{b}-\mathbf{x}_0^{\mathrm{CG}}\|_{\mathbf{A}}} \leq \min_{\substack{\deg(p)\leq k\\ p(0)=1}} \|p\|_\Lambda.$$

*Proof.* Observe that

$$\min_{\mathbf{x}\in\mathcal{K}_k(\mathbf{A},\mathbf{b})} \|\mathbf{A}^{-1}\mathbf{b}-\mathbf{x}\|_{\mathbf{A}} = \min_{\deg(p)<k} \|\mathbf{A}^{-1}\mathbf{b}-p(\mathbf{A})\mathbf{b}\|_{\mathbf{A}} \tag{5.4}$$

$$= \min_{\deg(p)<k} \|(\mathbf{I}-p(\mathbf{A})\mathbf{A})\mathbf{A}^{-1}\mathbf{b}\|_{\mathbf{A}} \tag{5.5}$$

$$= \min_{\substack{\deg(p)\leq k\\ p(0)=1}} \|p(\mathbf{A})\mathbf{A}^{-1}\mathbf{b}\|_{\mathbf{A}}. \tag{5.6}$$

Now, since $\mathbf{A}^{1/2}p(\mathbf{A}) = p(\mathbf{A})\mathbf{A}^{1/2}$, using the definition of the operator norm,

$$\|p(\mathbf{A})\mathbf{A}^{-1}\mathbf{b}\|_{\mathbf{A}} = \|p(\mathbf{A})\mathbf{A}^{1/2}\mathbf{A}^{-1}\mathbf{b}\| \leq \|p(\mathbf{A})\|\|\mathbf{A}^{1/2}\mathbf{A}^{-1}\mathbf{b}\|. \tag{5.7}$$

Finally, since $\|p(\mathbf{A})\|_2 = \|p\|_\Lambda$ and $\|\mathbf{A}^{1/2}\mathbf{A}^{-1}\mathbf{b}\| = \|\mathbf{A}^{-1}\mathbf{b}\|_\mathbf{A}$, combing these expressions and rearranging gives the desired result. ■

**Corollary 5.4.** The CG iterate satisfies the error bound

$$\frac{\|\mathbf{A}^{-1}\mathbf{b} - \mathbf{x}_k^{\mathrm{CG}}\|_\mathbf{A}}{\|\mathbf{A}^{-1}\mathbf{b} - \mathbf{x}_0^{\mathrm{CG}}\|_\mathbf{A}} \leq 2\exp\left(\frac{-2k}{\sqrt{\kappa(\mathbf{A})}}\right).$$

*Proof.* Since $\Lambda \subset \mathcal{I}$, $\|p\|_\Lambda \leq \|p\|_\mathcal{I}$ and hence

$$\min_{\substack{\deg(p)\leq k\\ p(0)=1}} \|p\|_\Lambda \leq \min_{\substack{\deg(p)\leq k\\ p(0)=1}} \|p\|_\mathcal{I}. \tag{5.8}$$

Observe that

$$p(x) := T_k\left(\frac{2x - \lambda_{\max} - \lambda_{\min}}{\lambda_{\max} - \lambda_{\min}}\right) \Big/ T_k\left(\frac{-\lambda_{\max} - \lambda_{\min}}{\lambda_{\max} - \lambda_{\min}}\right) \tag{5.9}$$

is the Chebyshev polynomial $T_k(x)$ shifted and scaled from $[-1, 1]$ to $\mathcal{I}$ and normalized to take value 1 at the origin. The max-growth property Lemma 3.25 of the Chebyshev polynomials implies this is the minimizer to the right hand side of (5.8). From the explicit formula for Chebyshev polynomials Lemma 3.21, one then obtains the well-known bound

$$\min_{\substack{\deg(p)\leq k\\ p(0)=1}} \|p\|_\mathcal{I} = 2\left(\left(\frac{\sqrt{\kappa}+1}{\sqrt{\kappa}-1}\right)^k + \left(\frac{\sqrt{\kappa}-1}{\sqrt{\kappa}+1}\right)^k\right)^{-1} \leq 2\left(\frac{\sqrt{\kappa}-1}{\sqrt{\kappa}+1}\right)^k, \tag{5.10}$$

where $\kappa = \kappa(\mathbf{A})$. While (5.10) is the standard bound in the numerical analysis literature, we find the stated bound, which follows immediately from basic properties of the exponential, easier to parse. ■

**Remark 5.5.** Corollary 5.4 asserts that a small condition number implies fast convergence. However, this is simply an upper bound; CG might converge very quickly on some ill-conditioned problems depending on how the eigenvalues are arranged. Thus, it is incorrect to claim that a smaller condition number means faster convergence of CG.

As an explicit example, the following bound shows that there are very ill-conditioned systems on which it is possible for CG to converge very rapidly.

**Corollary 5.6.** For any $\ell < k$ the CG iterate satisfies the bound

$$\frac{\|\mathbf{A}^{-1}\mathbf{b} - \mathbf{x}_k^{\mathrm{CG}}\|_{\mathbf{A}}}{\|\mathbf{A}^{-1}\mathbf{b} - \mathbf{x}_0^{\mathrm{CG}}\|_{\mathbf{A}}} \leq 2\exp\left(\frac{-2(k-\ell)}{\sqrt{\lambda_{\ell+1}/\lambda_{\min}}}\right).$$

Here the eigenvalues are ordered in descending order, $\lambda_1 \geq \lambda_2 \geq \cdots \geq \lambda_n$, so that $\lambda_1 = \lambda_{\max}$ and $\lambda_n = \lambda_{\min}$.

*Proof.* We will describe the $\ell = 1$ case. Consider the polynomial

$$p(x) := \left(1 - \frac{x}{\lambda_{\max}}\right) T_{k-1}\left(\frac{2x - \lambda_2 - \lambda_{\min}}{\lambda_2 - \lambda_{\min}}\right) \bigg/ T_{k-1}\left(\frac{-\lambda_2 - \lambda_{\min}}{\lambda_2 - \lambda_{\min}}\right). \tag{5.11}$$

This polynomial aims to be small on the interval $[\lambda_{\min}, \lambda_2] \cup \{\lambda_{\max}\}$ rather than on all of $\mathcal{I}$. This is done by constructing a shifted and scaled Chebyshev polynomial on $[\lambda_{\min}, \lambda_2]$ of *one degree lower*, and then adding a root at $\lambda_{\max}$. Since $|1 - x/\lambda_{\max}| < 1$ on $[\lambda_{\min}, \lambda_2]$, we can bound $|p(x)|$ on this interval using the same technique as we used to bound (5.9), and since $p(x)$ has a root at $\lambda_{\max}$ it does not matter that the Chebyshev polynomial might be larger there. The same approach can be used for larger $\ell$. ■

Comparing this bound to Corollary 5.4, we observe it is weaker in that the numerator depends on $k - \ell$ rather than $k$, but it is stronger in that the denominator depends on the "condition number" of all but the top-$\ell$ eigenvalues. If $\lambda_{\max} \gg \lambda_{\ell+1}$ then the better condition number in the denominator can make this bound much stronger than Corollary 5.4. This bound also makes it clear that the convergence of CG does not depend only on the condition number $\lambda_{\max}/\lambda_{\min}$, but the overall arrangement of the eigenvalues of $\mathbf{A}$.

**A posteriori error estimates** A priori bounds such as Corrolaries 5.3, 5.4 and 5.6 provide insight into how the convergence of CG is impacted by problem parameters such as the condition number or distribution of eigenvalues. However, such problem parameters are typically unknown, and so these bounds cannot be instantiated in practice. As such, other methods for estimating or bounding the error a posteriori are needed. We summarize a few simple approaches, and turn readers to the book [MT24] for more details and other more advanced approaches.

The simplest approach, which can be used for any method, is to compute the residual norm $\|\mathbf{b} - \mathbf{A}\mathbf{x}_k\|$. Note that

$$\|\mathbf{b} - \mathbf{A}\mathbf{x}_k\| = \|\mathbf{A}(\mathbf{A}^{-1}\mathbf{b} - \mathbf{x}_k)\| = \|\mathbf{A}^{-1}\mathbf{b} - \mathbf{x}_k\|_{\mathbf{A}^2}. \tag{5.12}$$

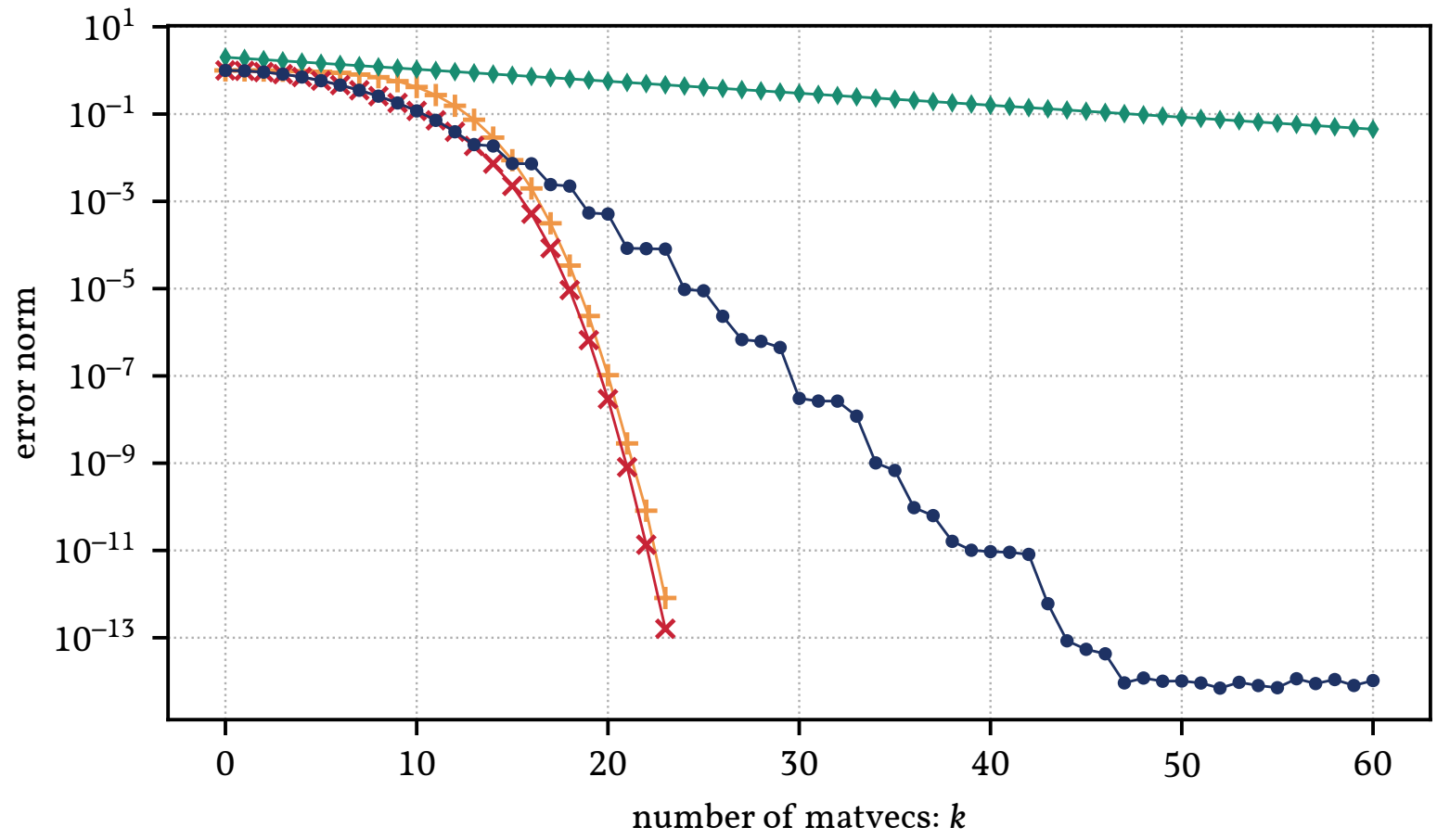

**Figure 5.1:** CG error $\|\mathbf{A}^{-1}\mathbf{b} - \mathbf{x}_k^{\mathrm{CG}}\|_{\mathbf{A}}/\|\mathbf{A}^{-1}\mathbf{b} - \mathbf{x}_0^{\mathrm{CG}}\|_{\mathbf{A}}$ with ( —×— ) and without ( —•— ) reorthgonalization, bound Corollary 5.3 on $\Lambda$ ( —+— ), and bound Corollary 5.4 on $\mathcal{I}$ ( —♦— ). *Takeaway*: The bound Corollary 5.3 typically accurately predicts the behavior of CG. The bound Corollary 5.4 can be very pessimistic, even when CG is run in finite precision arithmetic.

That is, the 2-norm of the residual is the $\mathbf{A}^2$-norm of the error. This in turn gives upper and lower bound for the $\mathbf{A}$-norm,

$$\sqrt{1/\lambda_{\max}}\|\mathbf{A}^{-1}\mathbf{b} - \mathbf{x}_k\|_{\mathbf{A}^2} \leq \|\mathbf{A}^{-1}\mathbf{b} - \mathbf{x}_k\|_{\mathbf{A}} \leq \sqrt{1/\lambda_{\min}}\|\mathbf{A}^{-1}\mathbf{b} - \mathbf{x}_k\|_{\mathbf{A}^2}. \tag{5.13}$$

A more sophisticated approach is to use the $\mathbf{A}$-norm optimality of the iterates. This optimality implies that for any $d \geq 0$, $\mathbf{A}^{-1}\mathbf{b} - \mathbf{x}_{k+d}$ and $\mathbf{x}_{k+d} - \mathbf{x}_k$ are $\mathbf{A}$-orthogonal and hence

$$\|\mathbf{A}^{-1}\mathbf{b} - \mathbf{x}_k\|_{\mathbf{A}}^2 = \|\mathbf{A}^{-1}\mathbf{b} - \mathbf{x}_{k+d}\|_{\mathbf{A}}^2 + \|\mathbf{x}_{k+d} - \mathbf{x}_k\|_{\mathbf{A}}^2. \tag{5.14}$$

If $\|\mathbf{A}^{-1}\mathbf{b} - \mathbf{x}_{k+d}\|_{\mathbf{A}} \ll \|\mathbf{A}^{-1}\mathbf{b} - \mathbf{x}_k\|_{\mathbf{A}}$ then

$$\|\mathbf{A}^{-1}\mathbf{b} - \mathbf{x}_k\|_{\mathbf{A}} \approx \|\mathbf{x}_{k+d} - \mathbf{x}_k\|_{\mathbf{A}}. \tag{5.15}$$

Since the $\mathbf{A}$-norm convergence of CG is monotonic, one might hope that this approximation is good even for small $d$. Note also that since (5.14) is an equality, $\|\mathbf{x}_{k+d} - \mathbf{x}_k\|_{\mathbf{A}}$ is in fact a lower-bound for the true error $\|\mathbf{A}^{-1}\mathbf{b} - \mathbf{x}_k\|_{\mathbf{A}}$.

**5.1.2 Finite precision arithmetic** The the main effects of finite precision arithmetic on CG are a delay of convergence and loss of maximal accuracy; see for instance [Meu06]. Both are described below and illustrated Figure 5.1.

- **Delay of convergence:** The error norms of the iterates in finite precision decreases more slowly than in exact arithmetic.

- **Lost of maximal accuracy:** The error norms eventually stops decreasing and plateau at some level of maximal accuracy.

Suppose the CG iterate Definition 5.1 is computed exactly from the quantities $\overline{\mathbf{Q}}_k$ and $\overline{\mathbf{T}}_k$ produced by the Lanczos algorithm in finite precision arithmetic. Greenbaum's theory (Section 4.2) implies that one can obtain a bound similar to Corollary 5.3 where $\Lambda$ is replaced with a union of small intervals each centered at an eigenvalue of $\mathbf{A}$. This in turn implies a bound similar to Corollary 5.4 still holds. In fact, this is also implied by the analysis in [DK91], and we will prove a more general version of this claim in Section 6.2. Note that, even in finite precision arithmetic, Corollary 5.4 is often pessimistic; see for instance Figure 5.1.

The exact behavior of CG depends on the implementation. Similar qualitative effects are observed for standard implementations of CG. However, if one uses a different implementation, the theoretical guarantees for the Lanczos algorithm in finite precision arithmetic described in Chapter 4 cannot be direct applied.

#### 5.1.3 Low-memory implementation

Although we have shown the $\mathbf{A}$-norm optimality of the CG iterate from Definition 5.1, we are missing any intuition for *why* this formulation is equivalent to the standard conjugate gradient method [HS52]. Perhaps the most notable difference is that Definition 5.1 involves the entire Lanczos basis $\mathbf{Q}_k$, while a standard CG implementation requires storing only several vectors of length $d$. In this section, we provide a brief argument showing the CG iterate can be computed without ever storing the $d \times k$ matrix $\mathbf{Q}_k$. A more detailed version of the argument, along with a full derivation of the standard CG algorithm can be found in [LS13, §2.5.1].

Note that since $\mathbf{A}$ is positive definite so is $\mathbf{T}_k$, and hence it has a Cholesky factorization $\mathbf{T}_k = \mathbf{L}_k\mathbf{L}_k^\mathsf{T}$ where $\mathbf{L}_k$ is some lower-triangular matrix. Define

$$\mathbf{P}_k := \begin{bmatrix} | & | & & | \\ \mathbf{p}_0 & \mathbf{p}_1 & \cdots & \mathbf{p}_{k-1} \\ | & | & & | \end{bmatrix} := \mathbf{Q}_k\mathbf{L}_k^{-\mathsf{T}}. \tag{5.16}$$

Observe that the vectors $\mathbf{p}_0, \ldots, \mathbf{p}_{k-1}$ are $\mathbf{A}$-orthogonal since

$$\mathbf{P}_k^\mathsf{T}\mathbf{A}\mathbf{P}_k = \mathbf{L}_k^{-1}\mathbf{Q}_k^\mathsf{T}\mathbf{A}\mathbf{Q}_k\mathbf{L}_k^{-\mathsf{T}} = \mathbf{L}_k^{-1}\mathbf{T}_k\mathbf{L}_k^{-\mathsf{T}} = \mathbf{I}. \tag{5.17}$$

Moreover, since $\mathbf{L}_k^{-\mathsf{T}}$ is upper triangular, we see that

$$\operatorname{span}\{\mathbf{p}_0, \ldots, \mathbf{p}_n\} = \operatorname{span}\{\mathbf{q}_0, \ldots, \mathbf{q}_n\} = \mathcal{K}_{n+1}(\mathbf{A}, \mathbf{b}). \tag{5.18}$$

Therefore, the $\mathbf{A}$-norm optimality of CG implies that $\mathbf{x}_k^{\mathrm{CG}}$ is obtained by projecting, in the $\mathbf{A}$-inner product, $\mathbf{A}^{-1}\mathbf{b}$ onto the $\mathbf{p}_0, \ldots, \mathbf{p}_{k-1}$. The iterate therefore

satisfies

$$\mathbf{x}_k^{\mathrm{CG}} = \sum_{n=0}^{k-1} (\mathbf{p}_n^\mathsf{T} \mathbf{A}\mathbf{A}^{-1}\mathbf{b})\mathbf{p}_n = \sum_{n=0}^{k-1} (\mathbf{p}_n^\mathsf{T}\mathbf{b})\mathbf{p}_n = \mathbf{x}_{k-1}^{\mathrm{CG}} + (\mathbf{p}_{k-1}^\mathsf{T}\mathbf{b})\mathbf{p}_{k-1}. \tag{5.19}$$

In particular, we observe that we can update the iterate $\mathbf{x}_{k-1}^{\mathrm{CG}}$ to obtain the iterate $\mathbf{x}_k^{\mathrm{CG}}$, so long as we have $\mathbf{p}_{k-1}$. We will now argue that we can obtain $\mathbf{p}_{k-1}$ from the Lanczos vectors $\mathbf{q}_{k-1}$ and $\mathbf{q}_{k-2}$. In particular, this means we do not need to store any of the previous Lanczos vectors.

The inverse Cholesky factor has the form

$$\mathbf{L}_k^{-1} = \mathbf{L}^{(k-1)} \cdots \mathbf{L}^{(2)}\mathbf{L}^{(1)}, \tag{5.20}$$

where $\mathbf{L}^{(n)}$ is chosen so as to a introduce a one on the diagonal by rescaling the $n$-th row and to introduce zeros below the diagonal of the $n$-th column of

$$\mathbf{L}^{(n-1)} \cdots \mathbf{L}^{(2)}\mathbf{L}^{(1)}\mathbf{T}_k \tag{5.21}$$

by subtracting some multiple of the $n$-th row. Since $\mathbf{T}_k$ is tridiagonal so is (5.21) and and $\mathbf{L}^{(n)}$ has the form

$$\mathbf{L}^{(m)} = \left[\begin{array}{c|cc|c} \mathbf{I}_{m-1} & & & \\ \hline & \times & 0 & \\ & \times & 1 & \\ \hline & & & \mathbf{I}_{m-k-2} \end{array}\right], \tag{5.22}$$

where we use "×" to indicate nonzero entries of $\mathbf{L}^{(m)}$. This implies that $\mathbf{L}_k^{-1}$ is in fact bidiagonal, and hence $\mathbf{L}_k^{-\mathsf{T}}$ is also bidiagonal. As such, we can compute $\mathbf{p}_{k-1}$ from only $\mathbf{T}_k$ and $\mathbf{q}_{k-2}$ and $\mathbf{q}_{k-1}$ as desired.

There are a number of additional cost savings which can be made. For instance, the Cholesky factorization of $\mathbf{T}_{k+1}$ is easily obtained from that of $\mathbf{T}_k$ and the new coefficients generated by the Lanczos algorithm. Moreover, the only part of the factorization that actually needs saved in order to proceed to the next step is the bottom right corner. We will not explain further.

## 5.2 Indefinite systems: MINRES

We now consider the case that the coefficient matrix $\mathbf{A}$ is possibly indefinite. When $\mathbf{A}$ is indefinite $\mathbf{T}_k$ can have an eigenvalue at zero, in which case the CG iterate Definition 5.1 undefined. Instead, it is common to use the minimum residual method (MINRES) [PS75]. In this section we derive some basic facts about MINRES and its relation to CG. We do not comment in detail on practical aspects such as efficient implementation or convergence in finite precision arithmetic, both of which are similar in spirit to CG and can be found in other texts.

**Definition 5.7.** The $k$-th MINRES iterate $\mathbf{x}_k^{\mathrm{M}}$ is given by

$$\mathbf{x}_k^{\mathrm{M}} = \mathbf{x}_k^{\mathrm{M}}(\mathbf{A}, \mathbf{b}) := \|\mathbf{b}\|\mathbf{Q}_k(\mathbf{T}_{k+1,k})^{\dagger}\mathbf{e}_1.$$

MINRES, as its name suggests, satisfies an residual optimality condition.

**Theorem 5.8.** If $\mathbf{A}$ is symmetric, the MINRES iterate satisfies the formula

$$\mathbf{x}_k^{\mathrm{M}} = \operatorname*{argmin}_{\mathbf{x}\in\mathcal{K}_k(\mathbf{A},\mathbf{b})} \|\mathbf{b} - \mathbf{A}\mathbf{x}\|.$$

*Proof.* Similar to the CG case,

$$\operatorname*{argmin}_{\mathbf{x}\in\mathcal{K}_k(\mathbf{A},\mathbf{b})} \|\mathbf{b} - \mathbf{A}\mathbf{x}\| = \mathbf{Q}_k \operatorname*{argmin}_{\mathbf{c}\in\mathbb{R}^k} \|\mathbf{b} - \mathbf{A}\mathbf{Q}_k\mathbf{c}\| = \mathbf{Q}_k \operatorname*{argmin}_{\mathbf{c}\in\mathbb{R}^k} \|\mathbf{b} - \mathbf{Q}_{k+1}\mathbf{T}_{k+1,k}\mathbf{c}\|. \tag{5.23}$$

The solution to this least squares problem is

$$\begin{aligned}\mathbf{Q}_k((\mathbf{Q}_{k+1}\mathbf{T}_{k+1,k})^{\mathsf{T}}\mathbf{Q}_{k+1}\mathbf{T}_{k+1,k})^{-1}(\mathbf{Q}_{k+1}\mathbf{T}_{k+1,k})^{\mathsf{T}}\mathbf{b}&\\ = \|\mathbf{b}\|\mathbf{Q}_k((\mathbf{T}_{k+1,k})^{\mathsf{T}}\mathbf{T}_{k+1,k})^{-1}(\mathbf{T}_{k+1,k})^{\mathsf{T}}\mathbf{e}_1.&\end{aligned} \tag{5.24}$$

Since $((\mathbf{T}_{k+1,k})^{\mathsf{T}}\mathbf{T}_{k+1,k})^{-1}\mathbf{T}_{k+1,k}^{\mathsf{T}} = (\mathbf{T}_{k+1,k})^{\dagger}$, we have the result. ■

We remark that, so long as $\mathbf{A}$ is invertible, the 2-norm of the residual is the same as the $\mathbf{A}^2$-norm of the error:

$$\|\mathbf{b} - \mathbf{A}\mathbf{x}\| = \|\mathbf{A}(\mathbf{A}^{-1}\mathbf{b} - \mathbf{x})\| = \|\mathbf{A}^{-1}\mathbf{b} - \mathbf{x}\|_{\mathbf{A}^2}. \tag{5.25}$$

So, while MINRES is typically stated in terms of a residual optimality condition, it also satisfies an equivalent error optimality condition, albeit in the $\mathbf{A}^2$-norm.

#### 5.2.1 Error bounds

The optimality of MINRES implies a bound on the spectrum of $\mathbf{A}$ analogous to Corollary 5.3.

**Corollary 5.9.** The MINRES iterate satisfies the error bounds

$$\frac{\|\mathbf{b} - \mathbf{A}\mathbf{x}_k^{\mathrm{M}}\|}{\|\mathbf{b} - \mathbf{A}\mathbf{x}_0^{\mathrm{M}}\|} \le \min_{\substack{\deg(p)\le k\\ p(0)=1}} \|p\|_{\Lambda}.$$

Corollary 5.9 allows us to obtain bounds for MINRES in terms of properties of polynomials. On positive definite system, MINRES satisfies bounds similar to Corrolaries 5.4 and 5.6. However, MINRES also satisfies bounds on indefinite systems. For instance, in [Gre97, §3.1] the following bounds are derived.

**Corollary 5.10.** Suppose $\Lambda \subset [a,b] \cup [c,d]$ where $b < 0 < c$ and $b-a = d-c$. The MINRES iterate satisfies the error bounds

$$\frac{\|\mathbf{b} - \mathbf{A}\mathbf{x}_k^{\mathrm{M}}\|}{\|\mathbf{b} - \mathbf{A}\mathbf{x}_0^{\mathrm{M}}\|} \leq 2\exp\left(\frac{-2\lfloor k/2 \rfloor}{\sqrt{|ad|/|bc|}}\right) \leq 2\exp\left(\frac{-2\lfloor k/2 \rfloor}{\kappa(\mathbf{A})}\right).$$

If the intervals are not situated symmetrically about zero, then the first bound can be significantly better than the second. When $b - a \neq d - c$ it is possible to obtain asymptotic rates for the best polynomial approximation in terms of Jacobi elliptic functions; see [Fis11, §3.3 and §3.4]. One can also derive a version of Corollary 5.10 analogous to Corollary 5.6.

### 5.2.2 Relations between MINRES and CG

From their formulas, it is clear the MINRES and CG iterates must be related, so one it is natural to ask whether the CG produces reasonable approximations to on indefinite systems. Unfortunately this is not generally possible since $\mathbf{T}_k$ can have eigenvalues at or near zero, in which case the CG iterate has infinite or large error. Even so, at iterations where the CG iterate is defined, one may wonder whether the error is similar to MINRES.

Define the residual vectors

$$\mathbf{r}_k^{\mathrm{M}} := \mathbf{b} - \mathbf{A}\mathbf{x}_k^{\mathrm{M}}, \qquad \mathbf{r}_k^{\mathrm{CG}} := \mathbf{b} - \mathbf{A}\mathbf{x}_k^{\mathrm{CG}}. \tag{5.26}$$

In Figure 5.2 we plot the residual norm convergence of CG and MINRES on an indefinite problem. While CG appears to have erratic convergence, its overall convergence seems to track that of MINRES.

Rearranging [MD20, Theorem 3.12], which is derived from a factorization of upper-Hessenberg matrices, the MINRES and CG residual norms satisfy the following relationship.

**Theorem 5.11.** The MINRES residual norms are related to the CG residual norms in that

$$\|\mathbf{r}_k^{\mathrm{M}}\| = \frac{1}{\sqrt{\sum_{n=0}^{k} 1/\|\mathbf{r}_n^{\mathrm{CG}}\|^2}}.$$

Theorem 5.11 implies what is perhaps the most well-known relation between the CG and MINRES residual norms. The following can also be derived more directly [Bro91; CG96].

**Theorem 5.12.** The CG residual norms are related to the MINRES residual norms in that

$$\|\mathbf{r}_k^{\mathrm{CG}}\| = \frac{\|\mathbf{r}_k^{\mathrm{M}}\|}{\sqrt{1 - \|\mathbf{r}_k^{\mathrm{M}}\|^2 / \|\mathbf{r}_{k-1}^{\mathrm{M}}\|^2}}.$$

Both Theorems 5.11 and 5.12 imply the well-known plateau-peak phenomenon: on iterations where MINRES makes good progress $\|\mathbf{r}_k^{\mathrm{M}}\|^2 \ll \|\mathbf{r}_{k-1}^{\mathrm{M}}\|^2$, the CG residual is similar to the MINRES residual, but on iterations where MINRES stagnates, the CG residual spikes. However, as stated, neither bound makes it clear whether it is necessary that CG encounter small residuals as $k$ increases.

In [Saa03], Theorem 5.11 is used to show that the *overall* convergence of CG on indefinite problems is similar to that of MINRES in that the smallest CG residual norm seen so far is never much bigger than the current MINRES residual norm.

**Theorem 5.13.** The overall residual norm convergence of CG is related to the MINRES residual norm in that

$$\min_{0 \le n \le k} \|\mathbf{r}_n^{\mathrm{CG}}\| \le \sqrt{k+1}\, \|\mathbf{r}_k^{\mathrm{M}}\|.$$

*Proof.* Bounding each term in the sum in Theorem 5.11 by the maximum we obtain

$$\|\mathbf{r}_k^{\mathrm{M}}\|_2 \ge \frac{1}{\sqrt{(k+1)\max_{0 \le j \le k} 1/\|\mathbf{r}_j^{\mathrm{CG}}\|_2^2}} = \frac{1}{\sqrt{k+1}} \min_{0 \le j \le k} \|\mathbf{r}_j^{\mathrm{CG}}\|_2, \tag{5.27}$$

Rearranging gives the desired result. ■

In [CM24] it is also proved that this bound is sharp; i.e. that there exist problems $(\mathbf{A}, \mathbf{b})$ for which the inequality is an equality.

## 5.3 Preconditioning

We conclude this section with a brief discussion on preconditioning, which in many ways has enabled the widespread use of KSMs for linear systems by transforming difficult problems to easier ones. While there are some “black-box" preconditioners, it is often the case that preconditioners are designed to work for systems arising from particular application areas [Saa03; Che05]. We do not discuss any specifics.

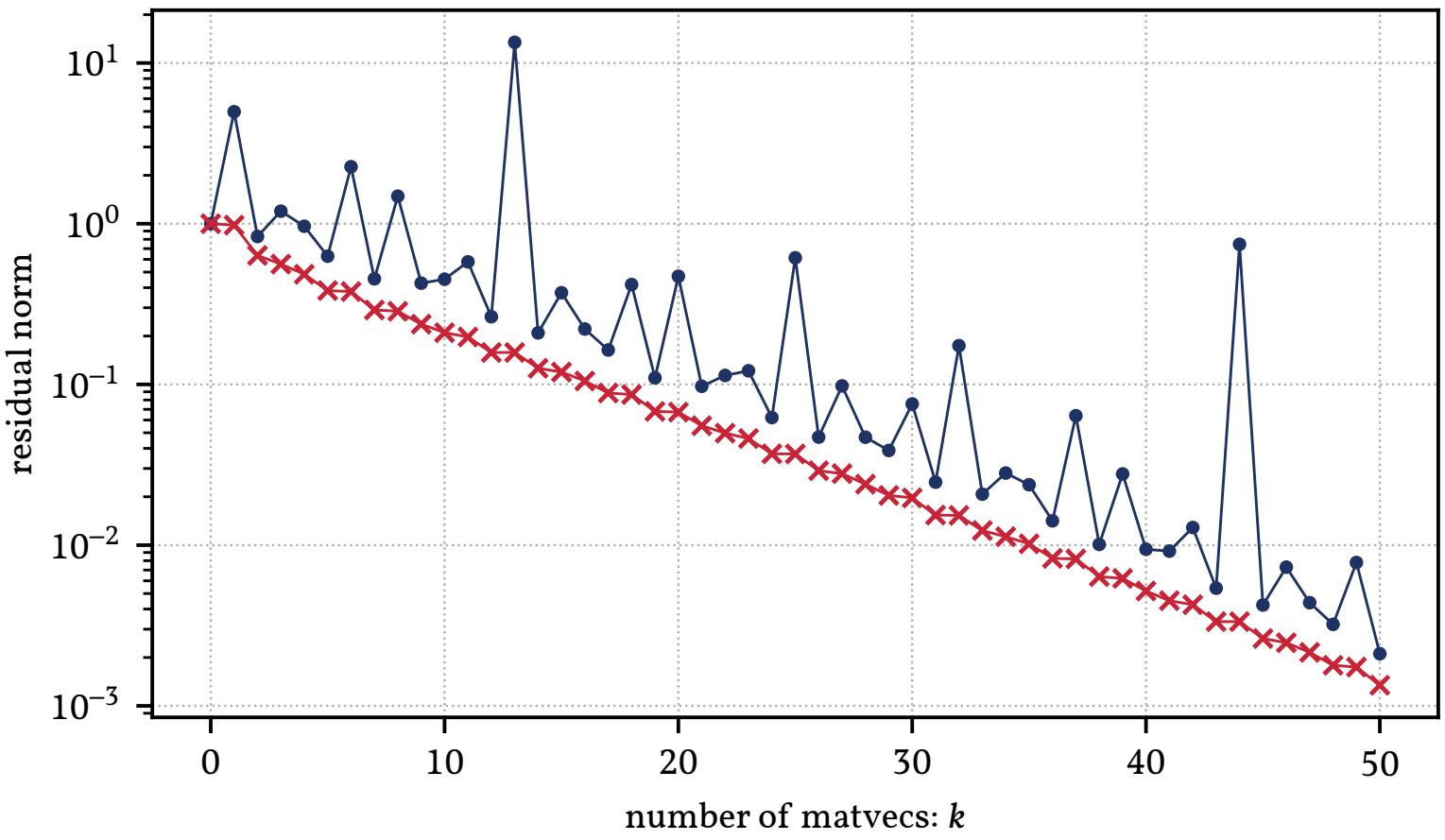

**Figure 5.2:** Convergence of CG ( —•— ) and MINRES ( —×— ) on an indefinite linear system. Observe the spikes of CG correspond to stagnation of MINRES, as described in Theorem 5.12. Moreover, despite intermittent spikes, the overall convergence of CG matches closely that of MINRES, as described in Theorem 5.13. *Takeaway*: While CG's convergence on indefinite systems may seem erratic, it is actually very closely related to the convergence of the optimal MINRES algorithm.

If $\mathbf{M}$ is invertible, then $\mathbf{A}^{-1}\mathbf{b} = \mathbf{M}^{\mathsf{T}}(\mathbf{M}\mathbf{A}\mathbf{M}^{\mathsf{T}})^{-1}\mathbf{M}\mathbf{b}$. Thus, we can obtain the solution to $\mathbf{A}\mathbf{x} = \mathbf{b}$ from the system of equations

$$\mathbf{M}\mathbf{A}\mathbf{M}^{\mathsf{T}}\mathbf{y} = \mathbf{M}\mathbf{b}, \qquad \mathbf{x} = \mathbf{M}^{\mathsf{T}}\mathbf{y}. \tag{5.28}$$

Here $\mathbf{M}$ is called a preconditioner, and (5.28) is referred to as the preconditioned system.

Note that $\mathbf{M}\mathbf{A}\mathbf{M}^{\mathsf{T}}$ is symmetric, and if $\mathbf{A}$ is positive definite, then so is $\mathbf{M}\mathbf{A}\mathbf{M}^{\mathsf{T}}$. If $\mathbf{M}\mathbf{A}\mathbf{M}^{\mathsf{T}}$ has better properties than $\mathbf{A}$, the the convergence of KSMs like CG or MINRES required to reach a given accuracy may be reduced. For instance, the best possible preconditioner would be $\mathbf{M} = \mathbf{A}^{-1/2}$ in which case the preconditioned system is solved in a single iteration by an iterative method like CG, and if $\mathbf{M}\mathbf{A}\mathbf{M}^{\mathsf{T}}$ has a small condition number, the fast convergence of CG is guaranteed.[1]

Of course, running an iterative method on the preconditioned system (5.28) is more expensive that on the original system. First, the preconditioned must be constructed, and second, the cost of an iteration is now increased. In particular, we must perform matrix-vector products with $\mathbf{M}\mathbf{A}\mathbf{M}^{\mathsf{T}}$ (which is typically done by first multiplying the vector with $\mathbf{M}^\mathsf{T}$, then with $\mathbf{A}$, and finally with $\mathbf{M}$). Thus, $\mathbf{M}$ must be chosen so as to be efficient to build and apply, while simultaneously improving the convergence properties of iterative methods.

[1] This does not mean that if the condition number of $\mathbf{M}\mathbf{A}\mathbf{M}^{\mathsf{T}}$ is smaller than that of $\mathbf{A}$ convergence is faster.

## TL;DR

The convergence of Lanczos-based methods like CG and MINRES for linear systems depends on fine-grained spectral properties of $\mathbf{A}$. Bounds based on the condition number are typically pessimistic, even in finite precision arithmetic.

# 6 Matrix functions times vectors

In this section we discuss methods for approximating the action of a matrix function on a vector,

$$f(\mathbf{A})\mathbf{b}. \tag{6.1}$$

In the special case $f(x) = 1/x$, this corresponds to solving a linear system $\mathbf{A}\mathbf{x} = \mathbf{b}$ which we discussed in Chapter 5. Beyond the multitude of applications of linear systems, matrix functions applied to vectors are used for computing the overlap operator in quantum chromodynamics [Esh+02], solving differential equations in applied math [Saa92; HL97], Gaussian process sampling in statistics [Ple+20], principal component projection and regression in data science [JS19], and a range of other applications [Hig08].

## 6.1 The Lanczos method for matrix function approximation

Perhaps the most successful general purpose KSM for approximating the action of symmetric matrix functions is the Lanczos method for matrix function approximation (Lanczos-FA). Early uses of Lanczos-FA were focused primarily on computing matrix exponentials applied to a vector; i.e. $f(x) = \exp(\beta x)$, and a number of papers studying the algorithm and its convergence properties were published around the same time [NW83; PL86; Vor87; DK88; DK89; GS92a; Saa92]. These early works were quickly followed by several papers demonstrating the effectiveness of Lanczos-FA in finite precision arithmetic [DK91; DK95; DGK98], a topic we discuss further in Section 6.2.

**Definition 6.1.** The $k$-th Lanczos-FA approximation to $f(\mathbf{A})\mathbf{b}$ is

$$\text{lan-FA}_k(f) = \text{lan-FA}_k(f; \mathbf{A}, \mathbf{b}) := \|\mathbf{b}\|\mathbf{Q}_k f(\mathbf{T}_k)\mathbf{e}_1.$$

Understanding why this simple algorithm works so well in practice is an ongoing topic of research [DK91; Dru08; DGK98; OSV12; MMS18; Che+22; Ams+24]. In the next several sections we will aim to provide some insight into the remarkable behavior of Lanczos-FA in exact and finite precision arithmetic.

### 6.1.1 Exactness of Lanczos-FA

It is well-known that the Lanczos-FA iterate is exact for low-degree polynomials.

**Theorem 6.2.** Suppose $\deg(p) < k$. Then,

$$\text{lan-FA}_k(p) = p(\mathbf{A})\mathbf{b}.$$

Theorem 6.2 is an immediate consequence of Lemma 3.19 on polynomial interpolation and Theorem 3.1 relating Lanczos to orthogonal polynomials. It is also an immediate consequence of Theorem 6.6 which, as we discuss in Section 6.2, provides insight into the finite precision behavior of the algorithm. We provide a third proof here.

*Proof.* By linearity, it suffices to verify that $\text{lan-FA}_k(x^n) = \mathbf{A}^n\mathbf{b}$ for all $n < k$. Observe that $\mathbf{Q}_k\mathbf{Q}_k^\mathsf{T}$ is the orthogonal projector onto $\mathcal{K}_k(\mathbf{A}, \mathbf{b})$. Hence, since $\mathbf{b}, \mathbf{A}\mathbf{b}, \ldots, \mathbf{A}^n\mathbf{b} \in \mathcal{K}_k(\mathbf{A}, \mathbf{b})$,

$$\begin{aligned}
\mathbf{A}^n\mathbf{b} &= \mathbf{Q}_k\mathbf{Q}_k^\mathsf{T}\mathbf{A}^n\mathbf{b} && (6.2)\\
&= \mathbf{Q}_k\mathbf{Q}_k^\mathsf{T}\mathbf{A}\mathbf{Q}_k\mathbf{Q}_k^\mathsf{T}\mathbf{A}^{n-1}\mathbf{b} && (6.3)\\
&= \mathbf{Q}_k\mathbf{T}_k\mathbf{Q}_k^\mathsf{T}\mathbf{A}^{n-1}\mathbf{b} && (6.4)\\
&\;\;\vdots && (6.5)\\
&= \mathbf{Q}_k\mathbf{T}_k^n\mathbf{Q}_k^\mathsf{T}\mathbf{b}. && (6.6)
\end{aligned}$$

The result then follows from the fact that $\mathbf{Q}_k^\mathsf{T}\mathbf{b} = \|\mathbf{b}\|\mathbf{e}_1$. ■

Theorem 6.2 also provides an interesting characterization of Lanczos-FA.

**Corollary 6.3.** Let $p(x)$ be the polynomial interpolant to $f(x)$ of degree less than $k$ at the eigenvalues of $\mathbf{T}_k$. Then

$$\text{lan-FA}_k(f) = p(\mathbf{A})\mathbf{b}.$$

*Proof.* Since $\mathbf{T}_k$ is of $k \times k$ and $p(x)$ is, by definition, the interpolant to $f(x)$ at the eigenvalues of $\mathbf{T}_k$, we have that $p(\mathbf{T}_k) = f(\mathbf{T}_k)$. Therefore, using Theorem 6.2,

$$\text{lan-FA}_k(f) = \|\mathbf{b}\|\mathbf{Q}_k f(\mathbf{T}_k)\mathbf{e}_1 = \|\mathbf{b}\|\mathbf{Q}_k p(\mathbf{T}_k)\mathbf{e}_1 = p(\mathbf{A})\mathbf{b}. \qquad (6.7)$$

■

**6.1.2 A simple error bound** We can use Theorem 6.2 to derive a simple error bound for Lanczos-FA in terms of polynomial approximation on $\mathcal{I} = [\lambda_{\min}, \lambda_{\max}]$. This bound essentially appears in [DK89; Saa92] and guarantees exponential convergence of Lanczos-FA for many matrix functions.

**Corollary 6.4.** The Lanczos-FA iterate satisfies

$$\|f(\mathbf{A})\mathbf{b} - \text{lan-FA}_k(f)\| \leq 2\|\mathbf{b}\| \min_{\deg(p)<k} \|f - p\|_{\mathcal{I}}.$$

*Proof.* For any polynomial $p(x)$ with $\deg(p) < k$, Theorem 6.2 asserts $p(\mathbf{A})\mathbf{b} = \text{lan-FA}_k(p)$. Then, by the triangle inequality

$$\|f(\mathbf{A})\mathbf{b} - \text{lan-FA}_k(f)\| \leq \|f(\mathbf{A})\mathbf{b} - p(\mathbf{A})\mathbf{b}\| + \|\text{lan-FA}_k(f) - \text{lan-FA}_k(p)\|. \quad (6.8)$$

Introduce the notation $e(x) := f(x) - p(x)$. Using the definition of operator norm and the fact that $\Lambda \subset \mathcal{I}$,

$$\|e(\mathbf{A})\mathbf{b}\| \leq \|e(\mathbf{A})\|\|\mathbf{b}\| = \|e\|_{\Lambda}\|\mathbf{b}\| \leq \|e\|_{\mathcal{I}}\|\mathbf{b}\|. \quad (6.9)$$

Similarly, since $\mathbf{Q}_k$ is orthogonal, and using Lemma 3.11 which implies that $\Lambda(\mathbf{T}_k) \subset \mathcal{I}$,

$$\|\mathbf{Q}_k e(\mathbf{T}_k)\mathbf{Q}_k^\mathsf{T}\mathbf{b}\| \leq \|\mathbf{Q}_k\|\|e(\mathbf{T}_k)\|\|\mathbf{Q}_k^\mathsf{T}\|\|\mathbf{b}\| \leq \|e\|_{\Lambda(\mathbf{T}_k)}\|\mathbf{b}\| \leq \|e\|_{\mathcal{I}}\|\mathbf{b}\|. \quad (6.10)$$

Plugging (6.9) and (6.10) in (6.8) gives the result. ■

Note that Corollary 6.4 is analogous to the bound Corollary 5.4, which guarantees CG converges exponentially at a rate depending on the square root of the condition number $\lambda_{\max}/\lambda_{\min}$. In particular, as with Corollary 5.4, Corollary 6.4 is typically not indicative of the true convergence behavior of Lanczos-FA.

**6.1.3 Spectrum adaptivity** When $f(x) = 1/x$, Lanczos-FA is mathematically equivalent to the CG iterate defined in Definition 5.1, and therefore satisfies strong spectrum dependent error guarantees such as Corollary 5.3 and Theorem 5.13. In practice, Lanczos-FA enjoys similar spectrum adaptivity for other functions $f(x)$, and as illustrated in Figure 6.1, often converges significantly faster than the bound Corollary 6.4. We now aim to provide some intuition for why this is the case, by relating Lanczos-FA to CG.

In many cases, $f(x)$ can be expressed as an integral of the form

$$f(x) = \int_{\Gamma} g(z)(x - z)^{-1}\,\mathrm{d}z, \quad (6.11)$$

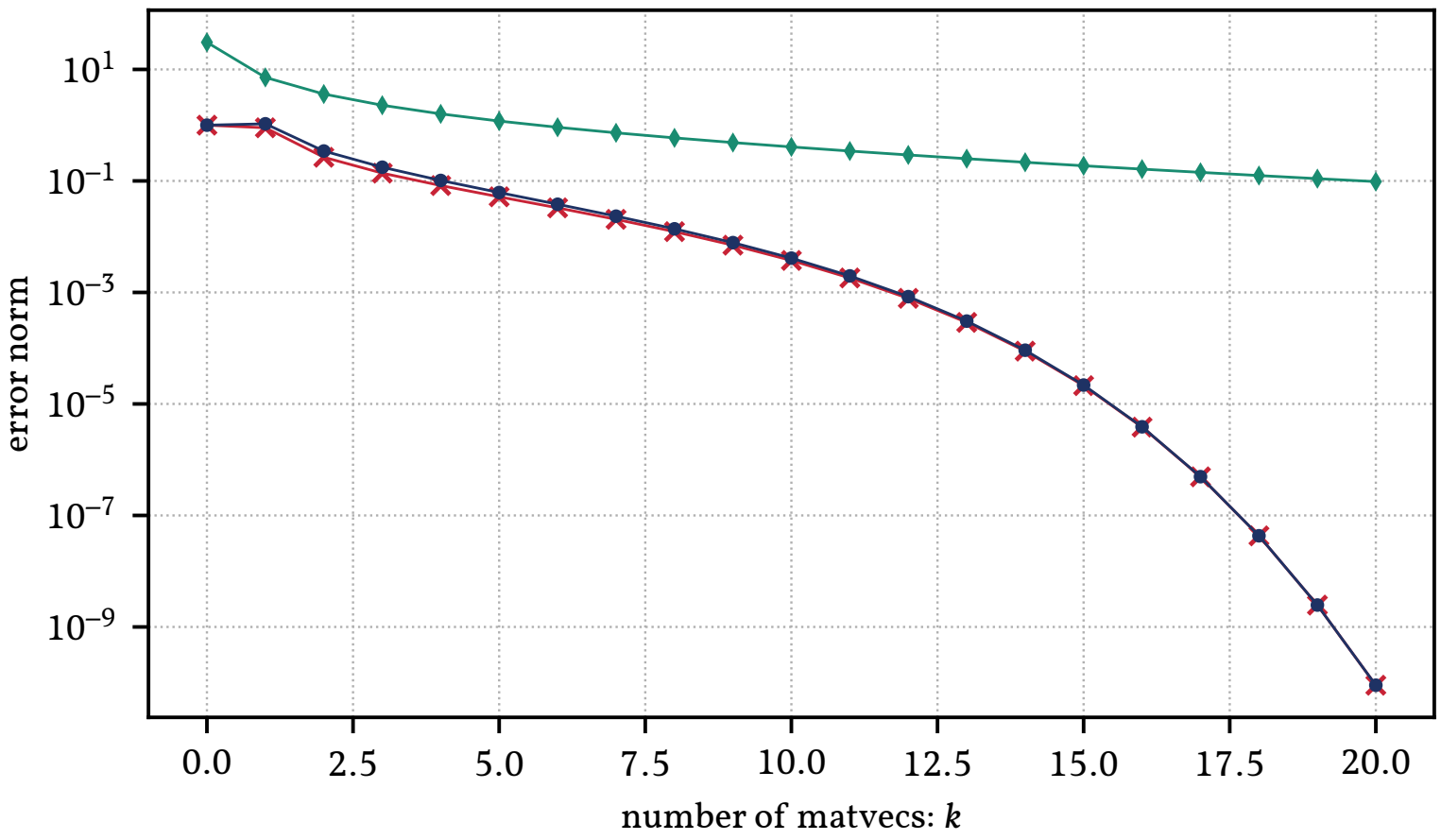

**Figure 6.1:** Comparison of the Euclidean norm convergence of Lanczos-FA ( ), optimal KSM ( ), and bound Corollary 6.4 ( ) for $f(x) = \sqrt{x}$. *Takeaway*: In many instances, Lanczos-FA exhibits nearly optimal convergence and the bound Corollary 6.4 is not at all indicative of the convergence of the algorithm.

where $g(z)$ is some function and $\Gamma$ a contour of integration. For instance, the inverse square root has the expression

$$\frac{1}{\sqrt{x}} = \frac{1}{\pi} \int_{(-\infty,0]} |z|^{-1/2}(x-z)^{-1}\mathrm{d}z. \tag{6.12}$$

This is an example of a Markov/Stieltjes function. Similarly, the indicator function that $x$ is inside a simple contour $\Gamma$ has the expression

$$f(x) = -\frac{1}{2\pi\mathrm{i}} \int_{\Gamma} (x-z)^{-1}\mathrm{d}z. \tag{6.13}$$

This is an example of the Cauchy Integral Formula. Finally, a rational function $f(x) = \sum_i w_i(x-z_i)^{-1}$ with simple poles $\{z_i\}$ can be expressed as

$$f(x) = \int_{-\infty}^{\infty} \Big(\textstyle\sum_i w_i\,\delta(z-z_i)\Big)(x-z)^{-1}\mathrm{d}z, \tag{6.14}$$

where $\delta(z)$ is a Dirac delta function centered at zero.

So long as (6.11) holds for each of the eigenvalues of $\mathbf{A}$,

$$f(\mathbf{A})\mathbf{b} = \int_{\Gamma} g(z)(\mathbf{A}-z\mathbf{I})^{-1}\mathbf{b}\,\mathrm{d}z. \tag{6.15}$$

Similarly, so long as (6.11) holds for each of the eigenvalues of $\mathbf{T}_k$,

$$\text{lan-FA}_k(f) = \int_{\Gamma} g(z)\|\mathbf{b}\|\mathbf{Q}_k(\mathbf{T}_k - z\mathbf{I})^{-1}\mathbf{e}_1\,\mathrm{d}z \tag{6.16}$$

$$= \int_\Gamma g(z)\,\text{lan-FA}_k((x-z)^{-1})\,\mathrm{d}z. \tag{6.17}$$

Assuming (6.15) and (6.17) hold, we can write the Lanczos-FA error as an integral

$$f(\mathbf{A})\mathbf{b} - \text{lan-FA}_k(f) = \int_\Gamma g(z)\big[(\mathbf{A}-z\mathbf{I})^{-1}\mathbf{b} - \text{lan-FA}_k((x-z)^{-1})\big]\,\mathrm{d}z. \tag{6.18}$$

This in turn leads to bounds for the norm of the Lanczos-FA error (6.18)

$$\|f(\mathbf{A})\mathbf{b} - \text{lan-FA}_k(f)\| \leq \int_\Gamma |g(z)|\|(\mathbf{A}-z\mathbf{I})^{-1}\mathbf{b} - \text{lan-FA}_k((x-z)^{-1})\||\mathrm{d}z|. \tag{6.19}$$

Thus, to understand the Lanczos-FA error, it suffices to understand the error of the integrand at each value $z$ in the contour $\Gamma$.

The key observation is that

$$\text{lan-FA}_k((x-z)^{-1}) \tag{6.20}$$

is equivalent to the CG iterate Definition 5.1 applied to the linear system $(\mathbf{A} - z\mathbf{I})\mathbf{x} = \mathbf{b}$. Indeed, the Lanczos factorization (2.12) can be shifted (even for complex $z$) to obtain

$$(\mathbf{A}-z\mathbf{I})\mathbf{Q}_k = \mathbf{Q}_k(\mathbf{T}_k - z\mathbf{I}) + \beta_{k-1}\mathbf{q}_k\mathbf{e}_k^\mathsf{T}. \tag{6.21}$$

That is, Lanczos applied to $(\mathbf{A}, \mathbf{b})$ for $k$ steps produces output $\mathbf{Q}_k$ and $\mathbf{T}_k$ satisfying (2.12) while Lanczos applied to $(\mathbf{A}-z\mathbf{I}, \mathbf{b})$ for $k$ steps produces output $\mathbf{Q}_k$ and $\mathbf{T}_k - z\mathbf{I}$ satisfying (6.21). Thus, the norm of the error $\|(\mathbf{A}-z\mathbf{I})^{-1}\mathbf{b} - \text{lan-FA}_k((x-z)^{-1})\|$ can be bounded (a priori or a posteriori) using theory about CG Chapter 5. This general approach is commonly used to derive bounds for Lanczos-FA and related algorithms [FS08b; FS09; ITS09; Fro+13; FGS14; FS15; Che+22].

**6.1.4 Is Lanczos-FA nearly optimal?** If $f(x) = 1/x$, the Lanczos-FA iterate coincides with the CG iterate defined in Definition 5.1 and is therefore an optimal approximation (in the $\mathbf{A}$-norm) if $\mathbf{A}$ is positive definite; see Theorem 5.2. Even when $\mathbf{A}$ is not positive definite, Lanczos-FA is nearly optimal in the sense of overall convergence; see Theorem 5.13.

One may wonder whether Lanczos-FA exhibits similar near-optimality guarantees for other functions. In Figure 6.1 we show the convergence of Lanczos-FA, as well as the optimal KSM error and the bound Corollary 6.4. In this example the convergence of Lanczos-FA closely tracks the convergence of the optimal KSM, converging significantly faster than the bound Corollary 6.4. In fact, it is hard to find examples where similar behavior is not observed.

For a class of functions called Stieltjes functions, [Sch25] proves that Lanczos-FA is nearly optimal; see [CP26] for follow-up work.

**Theorem 6.5.** Suppose $\mathbf{A}$ is positive definite and let $f$ be a Stieltjes function. Then

$$\|f(\mathbf{A})\mathbf{b} - \text{lan-FA}_k(f)\| \leq \left(1 + \kappa(\mathbf{A})^2\right) \min_{\mathbf{x} \in \mathcal{K}_k(\mathbf{A},\mathbf{b})} \|f(\mathbf{A})\mathbf{b} - \mathbf{x}\|.$$

## 6.2 Finite precision arithmetic

The main qualitative effects on Lanczos-FA are a delay of convergence and loss of maximal accuracy, similar to what we described in Section 5.1.2 for CG.

The main theoretical guarantee for Lanczos-FA in finite precision arithmetic asserts that Lanczos-FA accurately applies Chebyshev polynomials.

**Theorem 6.6** (informal; see [DK91]). Suppose $\|\mathbf{A}\| = \|\mathbf{b}\| = 1$ and that Lanczos algorithm is run for $k$ iterations on a computer with relative machine precision $\epsilon_{\text{mach}} = O(\text{poly}(k)^{-1})$. Then, for all $n \leq k-1$,

$$\|T_n(\mathbf{A})\mathbf{b} - \overline{\mathbf{Q}}_k T_n(\overline{\mathbf{T}}_k)\mathbf{e}_1\| = O(\text{poly}(k)\,\epsilon_{\text{mach}}).$$

We can then use Theorem 6.6 to derive a bound similar to Corollary 6.4 in finite precision arithmetic. This bound is implicit in [DK91]; see also [MMS18] for a more detailed analysis and some extensions.

**Corollary 6.7.** Suppose $\|\mathbf{A}\| = \|\mathbf{b}\| = 1$ and that Lanczos algorithm is run for $k$ iterations on a computer with relative machine precision $\epsilon_{\text{mach}} = O(\text{poly}(k)^{-1})$. Then, with $\overline{\mathcal{I}} = [-1-\eta_k, 1+\eta_k]$ where $\eta_k = O(\text{poly}(k)\,\epsilon_{\text{mach}})$,

$$\|f(\mathbf{A})\mathbf{b} - \overline{\mathbf{Q}}_k f(\overline{\mathbf{T}}_k)\mathbf{e}_1\| \leq O(k) \min_{\deg(p)<k} \|f-p\|_{\overline{\mathcal{I}}} + O(\|f\|_{\overline{\mathcal{I}}}\,\text{poly}(k)\,\epsilon_{\text{mach}}).$$

We will prove these results in Section 6.2.1.

**Remark 6.8.** Corollary 6.7 is for the quantity $\|\mathbf{b}\|\overline{\mathbf{Q}}_k f(\overline{\mathbf{T}}_k)\mathbf{e}_1$ and *not* the quantity $\overline{\mathbf{Q}}_k f(\overline{\mathbf{T}}_k)\overline{\mathbf{Q}}_k^\mathsf{T}\mathbf{b}$. Indeed, while these quantities are equivalent in exact arithmetic, they might not even be nearby in finite precision arithmetic. Lemma 3.19 provides intuition into why the analyzed quantity is the right one, without requiring the full analysis of Corollary 6.7. The convergence of these quantities is illustrated in Figure 6.2, where we observe the latter quantity does not even converge.

**Remark 6.9.** Despite widespread claims that Chebyshev-based methods are "more stable" than Lanczos-based methods, the precise version of Corollary 6.7 is almost identical to what can be obtained for a method based directly on the Chebyshev recurrence (3.15); see also Remark 4.5.

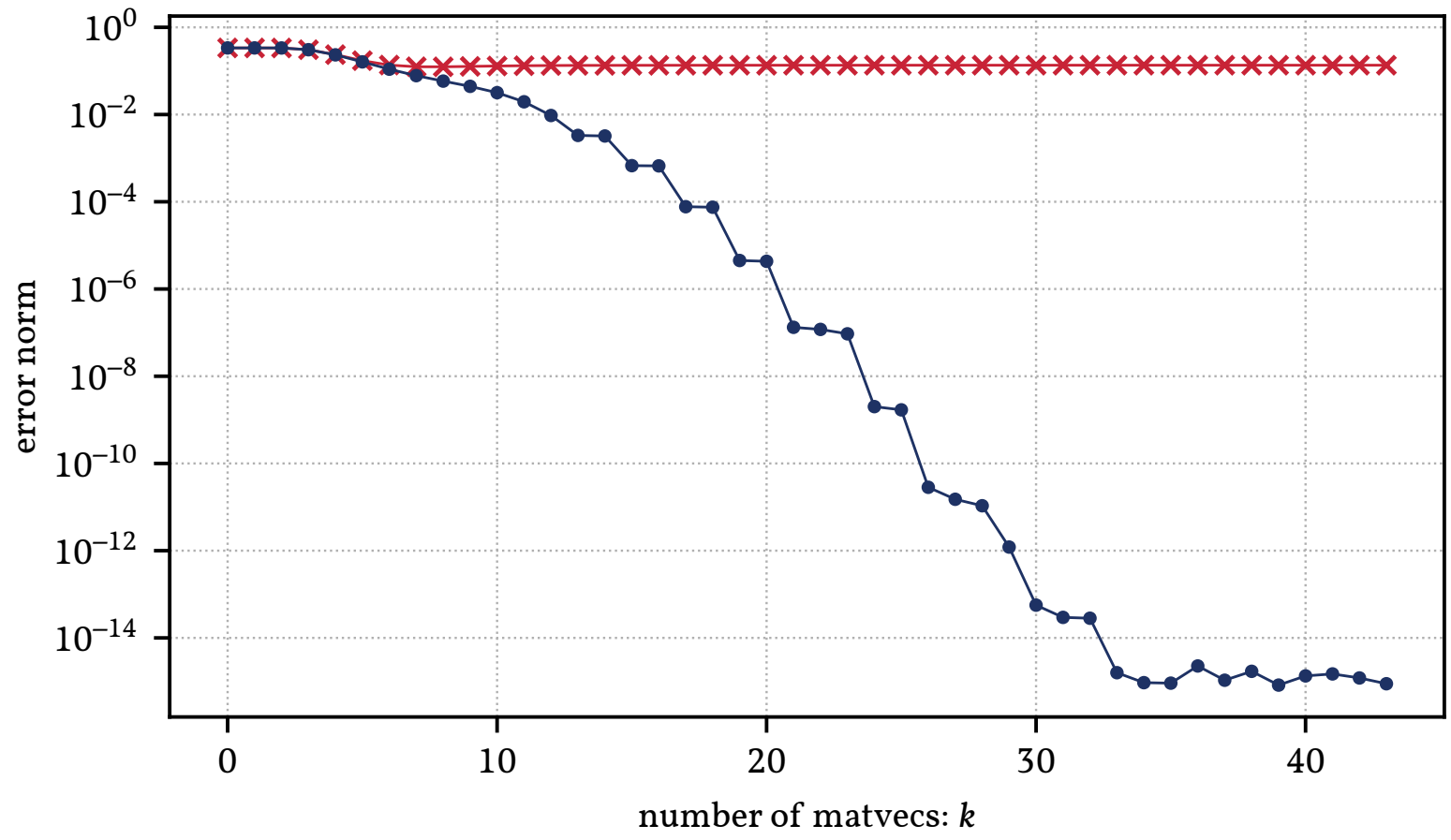

**Figure 6.2:** Comparison of the Euclidean error approximating $f(\mathbf{A})\mathbf{b}$ of $\|\mathbf{b}\|\overline{\mathbf{Q}}_k f(\overline{\mathbf{T}}_k)\mathbf{e}_1$ ( —•— ) and $\overline{\mathbf{Q}}_k f(\overline{\mathbf{T}}_k)\overline{\mathbf{Q}}_k^{\mathsf{T}}\mathbf{b}$ ( —×— ) for $f(x) = \exp(-x)$ in finite precision arithmetic. *Take-away*: Use the expression $\|\mathbf{b}\|\overline{\mathbf{Q}}_k f(\overline{\mathbf{T}}_k)\mathbf{e}_1$, otherwise Lanczos-FA might not converge!!!

### 6.2.1 Proofs

We will use the following standard fact about perturbed Chebyshev recurrences [Cle55]. This is a special case of a more general formula involving the associated polynomials of some family of orthogonal polynomials.

**Lemma 6.10.** Suppose that, for $n = 2, 3, \ldots,$

$$d_0 = 0, \qquad d_1 = f_0, \qquad d_n(x) = 2xd_{n-1}(x) - d_{n-2}(x) + 2f_{n-1}.$$

Then, introducing the notation $U_{-1}(x) = 0$, for $n = 0, 1, \ldots,$

$$d_n(x) = U_{n-1}(x)f_0 + 2\sum_{i=2}^{n} U_{n-i}(x)f_{i-1}.$$

*Proof.* Suppose the lemma holds for $i < n$. Then,

$$d_n(x) = 2xd_{n-1}(x) - d_{n-2}(x) + 2f_{n-1} \tag{6.22}$$

$$= 2x\left(U_{n-2}(x)f_0 + 2\sum_{i=2}^{n-1} U_{n-1-i}(x)f_{i-1}\right)$$

$$- \left(U_{n-3}(x)f_0 + 2\sum_{i=2}^{n-2} U_{n-2-i}(x)f_{i-1}\right) + 2f_{n-1} \quad (6.23)$$

$$= \left(2xU_{n-2}(x) - U_{n-3}(x)\right)f_0$$

$$+ 2\left(\sum_{i=2}^{n-2} \left(2xU_{n-1-i}(x) - U_{n-2-i}(x)\right)f_{i-1}\right) + 4xf_{n-2} + 2f_{n-1} \quad (6.24)$$

$$= U_{n-1}(x)f_0 + 2\left(\sum_{i=2}^{n-2} U_{n-i}(x)f_{i-1}\right) + 2U_1(x)f_{n-2} + 2U_0(x)f_{n-1} \quad (6.25)$$

$$= U_{n-1}(x)f_0 + 2\sum_{i=2}^{n} U_{n-i}(x)f_{i-1}. \quad (6.26)$$

The result follows as the base case is assumed. ■

The proofs of Theorem 6.6 and Corollary 6.7 are now straightforward.

*Proof of Theorem 6.6.* For notational brevity, define the vectors

$$\mathbf{t}_n := T_n(\mathbf{A})\mathbf{b}, \qquad \bar{\mathbf{t}}_n := T_n(\overline{\mathbf{T}}_k)\mathbf{e}_1, \qquad \mathbf{d}_n := \mathbf{t}_n - \overline{\mathbf{Q}}_k\bar{\mathbf{t}}_n. \quad (6.27)$$

Since $k > 1$, using the notation in (6.27) and recalling the perturbed recurrence (4.2), we have

$$\mathbf{d}_0 = \mathbf{b} - \overline{\mathbf{Q}}_k\mathbf{e}_1 = \mathbf{0}, \qquad \mathbf{d}_1 = \mathbf{A}\mathbf{b} - \overline{\mathbf{Q}}_k\overline{\mathbf{T}}_k\mathbf{e}_1 = (\bar{\beta}_{k-1}\bar{\mathbf{q}}_k\mathbf{e}_k^\mathsf{T} + \mathbf{F}_k)\mathbf{e}_1 = \mathbf{F}_k\bar{\mathbf{t}}_0. \quad (6.28)$$

For $n = 2, \ldots, k-1$, we can use the definitions of $\mathbf{t}_n$ and $\bar{\mathbf{t}}_n$, the definition of the Chebyshev polynomials, and the perturbed recurrence (4.2) to write

$$\mathbf{d}_n = (2\mathbf{A}\mathbf{t}_{n-1} - \mathbf{t}_{n-2}) - (2\overline{\mathbf{Q}}_k\overline{\mathbf{T}}_k\bar{\mathbf{t}}_{n-1} - \overline{\mathbf{Q}}_k\bar{\mathbf{t}}_{n-2}) \quad (6.29)$$

$$= 2(\mathbf{A}\mathbf{t}_{n-1} - (\mathbf{A}\overline{\mathbf{Q}}_k - \bar{\beta}_{k-1}\bar{\mathbf{q}}_k\mathbf{e}_k^\mathsf{T} - \mathbf{F}_k)\bar{\mathbf{t}}_{n-1}) - (\mathbf{t}_{n-2} - \overline{\mathbf{Q}}_k\bar{\mathbf{t}}_{n-2}) \quad (6.30)$$

$$= 2(\mathbf{A}\mathbf{t}_{n-1} - (\mathbf{A}\overline{\mathbf{Q}}_k\bar{\mathbf{t}}_{n-1} - \bar{\beta}_{k-1}\bar{\mathbf{q}}_k\mathbf{e}_k^\mathsf{T}\bar{\mathbf{t}}_{n-1} - \mathbf{F}_k\bar{\mathbf{t}}_{n-1})) - \mathbf{d}_{n-2} \quad (6.31)$$

Note that $(\overline{\mathbf{T}}_k)^i$ has half bandwidth $i$, so $(\overline{\mathbf{T}}_k)^i$ is zero in the bottom left entry provided $i < k-1$. Since $T_i$ is a degree $i$ polynomial, this implies that that $\mathbf{e}_k^\mathsf{T}\bar{\mathbf{t}}_i = \mathbf{e}_k^\mathsf{T}T_i(\overline{\mathbf{T}}_k)\mathbf{e}_1 = 0$ for any $i < k-1$. Since $n < k$, applying this with $i = n-1$ we find

$$\mathbf{d}_n = 2\mathbf{A}\mathbf{d}_{n-1} - \mathbf{d}_{n-2} + 2\mathbf{F}_k\bar{\mathbf{t}}_{n-1}. \quad (6.32)$$

Lemma 6.10 with $x \to \mathbf{A}$, $d_n(x) \to \mathbf{d}_n$, and $f_n \to \mathbf{F}_k\bar{\mathbf{t}}_n$ allows us to obtain an explicit expression

$$\mathbf{d}_n = U_{n-1}(\mathbf{A})\mathbf{F}_k\bar{\mathbf{t}}_0 + 2\sum_{i=2}^{n} U_{n-i}(\mathbf{A})\mathbf{F}_k\bar{\mathbf{t}}_{i-1}. \tag{6.33}$$

Intuitively (6.33) is small because $\mathbf{F}_k$ is small (by Paige's analysis in Section 4.1) and the Chebyshev polynomials $U_\ell(x)$ and $T_\ell(x)$ are not very large on $[-1, 1]$. It is slightly annoying to prove this though, since it is possible that $\|\overline{\mathbf{T}}_k\|$ can be slightly larger than one.

Our assumption on $\epsilon_{\mathrm{mach}}$ and Theorem 4.1 guarantees that $\|\overline{\mathbf{T}}_k\| \leq 1 + 1/(2k^2)$. Hence, Lemmas 3.26 and 3.27 guarantees

$$\|T_\ell(\overline{\mathbf{T}}_k)\| \leq 2, \qquad \|U_\ell(\mathbf{A})\| \leq k. \tag{6.34}$$

Theorem 4.1 also gives the bound $\|\mathbf{F}_k\| = O(\mathrm{poly}(k)\,\epsilon_{\mathrm{mach}})$.

Therefore, we can apply the triangle inequality to (6.33), apply the above bounds, and use the fact $n \leq k$ to obtain the bound

$$\|\mathbf{d}_n\| \leq 2\sum_{i=1}^{n} \|U_{n-i}(\mathbf{A})\|\|\mathbf{F}_k\|\|\bar{\mathbf{t}}_{i-1}\| = O(\mathrm{poly}(k)\,\epsilon_{\mathrm{mach}}). \tag{6.35}$$

This is the desired result. ■

*Proof of Corollary 6.7.* This result intuitively follows because any function bounded in $[-1, 1]$ has a Chebyshev series with bounded coefficients. As with the proof of Theorem 6.6, the result is slightly complicated by the fact that $\|\overline{\mathbf{T}}_k\|$ can be greater than one.

Theorem 4.1 ensures $\Lambda(\overline{\mathbf{T}}_k) \subset \bar{\mathcal{I}}$. Choose

$$p(x) := \operatorname*{argmin}_{\deg(p)<k} \|f - p\|_{\bar{\mathcal{I}}}. \tag{6.36}$$

Define $e(x) := f(x) - p(x)$. The triangle inequality gives

$$\|f(\mathbf{A})\mathbf{b} - \overline{\mathbf{Q}}_k f(\overline{\mathbf{T}}_k)\mathbf{e}_1\| \leq \|p(\mathbf{A})\mathbf{b} - \overline{\mathbf{Q}}_k p(\overline{\mathbf{T}}_k)\mathbf{e}_1\| + \|e(\mathbf{A})\mathbf{b}\| + \|\overline{\mathbf{Q}}_k e(\overline{\mathbf{T}}_k)\mathbf{e}_1\|. \tag{6.37}$$

According to Lemma 3.24 we can expand $p(x)$ in the Chebyshev basis

$$p(x) = c_0 T_0(x) + 2\sum_{n=1}^{k-1} c_n T_n(x), \qquad c_n = \int p(x)T_n(x)\mu_T(x)\mathrm{d}x. \tag{6.38}$$

Applying the triangle inequality and Theorem 6.6 we find that

$$\|p(\mathbf{A})\mathbf{b} - \overline{\mathbf{Q}}_k p(\overline{\mathbf{T}}_k)\mathbf{e}_1\| \tag{6.39}$$

$$= \left\| c_0\big(T_0(\mathbf{A})\mathbf{b} - \overline{\mathbf{Q}}_k T_0(\overline{\mathbf{T}}_k)\mathbf{e}_1\big) + 2\sum_{n=1}^{k-1} c_n\big(T_n(\mathbf{A})\mathbf{b} - \overline{\mathbf{Q}}_k T_n(\overline{\mathbf{T}}_k)\mathbf{e}_1\big) \right\| \quad (6.40)$$

$$\leq |c_0|\|T_0(\mathbf{A})\mathbf{b} - \overline{\mathbf{Q}}_k T_0(\overline{\mathbf{T}}_k)\mathbf{e}_1\| + 2\sum_{n=1}^{k-1} |c_n|\|T_n(\mathbf{A})\mathbf{b} - \overline{\mathbf{Q}}_k T_n(\overline{\mathbf{T}}_k)\mathbf{e}_1\|. \quad (6.41)$$

By the triangle inequality and our choice of $p(x)$,

$$\|p\|_{\bar{\mathcal{I}}} \leq \|f\|_{\bar{\mathcal{I}}} + \|f-p\|_{\bar{\mathcal{I}}} \leq \|f\|_{\bar{\mathcal{I}}} + \|f-0\|_{\bar{\mathcal{I}}} = 2\|f\|_{\bar{\mathcal{I}}} \quad (6.42)$$

Hence, using Lemma 3.26 which asserts $\|T_n\|_{\mathcal{I}} \leq 1$,

$$|c_n| \leq \int |p(x)||T_n(x)|\mu_T(x)\mathrm{d}x \leq \|p\|_{[-1,1]} \int \mu_T(x)\mathrm{d}x \leq 2\|f\|_{\bar{\mathcal{I}}}. \quad (6.43)$$

Plugging this bound into (6.41) we get a bound

$$\|p(\mathbf{A})\mathbf{b} - \overline{\mathbf{Q}}_k p(\overline{\mathbf{T}}_k)\mathbf{e}_1\| \leq \|f\|_{\bar{\mathcal{I}}} O(\mathrm{poly}(k)\,\epsilon_{\mathrm{mach}}). \quad (6.44)$$

In addition, we have that

$$\|e(\mathbf{A})\mathbf{b}\| \leq \|e\|_{\Lambda} \leq \|e\|_{\mathcal{I}} \leq \|e\|_{\bar{\mathcal{I}}}. \quad (6.45)$$

From Theorem 4.1 we have $\|\overline{\mathbf{q}}_n\| \leq 1 + O(\mathrm{poly}(k)\epsilon_{\mathrm{mach}})$. Our assumption on $\epsilon_{\mathrm{mach}}$ ensures $\|\overline{\mathbf{q}}_n\| = O(1)$ and hence $\|\overline{\mathbf{Q}}_k\| = O(k)$. Therefore

$$\|\overline{\mathbf{Q}}_k e(\overline{\mathbf{T}}_k)\mathbf{e}_1\| \leq \|\overline{\mathbf{Q}}_k\|\|e\|_{\Lambda(\overline{\mathbf{T}}_k)} \leq O(k)\|e\|_{\bar{\mathcal{I}}}. \quad (6.46)$$

Finally, plugging (6.44)–(6.46) into (6.37) gives the result. ■

## 6.3 Low-memory algorithms

An apparent downside of Lanczos-FA compared with explicit polynomial approaches or CG is that the Krylov basis $\mathbf{Q}_k$ must be stored. However, there are several ways to avoid this storage cost.

**6.3.1 Two-pass Lanczos-FA** This storage cost can be avoided by incurring additional computational cost. In particular, we can use an implementation called two pass Lanczos-FA [Bor00; FS08a]. On the first pass, the tridiagonal matrix $\mathbf{T}_k$ is computed using the short-recurrence version of Lanczos; i.e., without storing all of $\mathbf{Q}_k$. Once $\mathbf{T}_k$ has been computed, $f(\mathbf{T}_k)\mathbf{e}_1$ can be evaluated. Lanczos is then run again and the vector $\mathbf{Q}_k f(\mathbf{T}_k)\mathbf{e}_1$ is computed as the columns of $\mathbf{Q}_k$ become available. Note that on the second run, the *exact same*[1] Lanczos vectors can be computed without any inner products by using the values computed in the first run and stored in $\mathbf{T}_k$.

Such an approach can be generalized by re-generating the Lanczos recurrence from multiple points simultaneously on the second pass [Li22]. Specifically, on the first pass, vectors $\mathbf{q}_j$ and $\mathbf{q}_{j-1}$ can be saved for $j = 0, q, 2q, \ldots$. Then, on the second pass, the rest of the Lanczos vectors can be constructed by continuing the three-term Lanczos recurrence (2.12) from each of the roughly $d/q$ start points in parallel. The case $q = d$ gives the original two-pass approach.

**6.3.2 Multi-shift CG/MINRES** Often $f(x)$ is well approximated by a rational function

$$f(x) \approx \sum_{i=1}^{q} \frac{w_i}{x - z_i}. \tag{6.47}$$

For instance, the integral representation (6.11) can be discretized using a quadrature rule. In this case,

$$f(\mathbf{A})\mathbf{b} \approx \sum_{i=1}^{q} w_i(\mathbf{A} - z_i\mathbf{I})^{-1}\mathbf{b}. \tag{6.48}$$

Using the shift-invariance property of Krylov subspaces (6.21), we can *simultaneously* apply a low-memory implementation of CG or MINRES to the shifted linear systems $(\mathbf{A} - z_i\mathbf{I})\mathbf{x}_i = \mathbf{b}$ without using additional matrix-vector products with $\mathbf{A}$. This requires $q$ times the memory of a single run of CG or MINRES, but does not increase the number of matrix-vector products. For a survey on methods of this flavor, see [GS21].

> **Remark 6.11.** The accuracy of (6.48) is limited by not only how well each term in the rational function is approximated by CG/MINRES, but also by the accuracy of the rational approximation (6.47).

## TL;DR

Lanczos-FA is a general purpose method for approximating $f(\mathbf{A})\mathbf{b}$. The method is closely related to CG, and therefore enjoys similar spectrum adaptivity, often exhibiting similar behavior to the best possible KSM. Bounds for Lanczos-FA based on best polynomial approximation on an interval are often pessimistic, but still hold in finite precision arithmetic.

[1] This is true, even in finite precision arithmetic, so long as computations are done deterministically. However, in many computing environments computations may be nondeterministic. There are several efforts to ensure reproducibility in such environments [Iak+15; ADN20].

# 7 Quadratic forms and trace approximation

In this section we discuss quadrature-based methods for approximating the quadratic form

$$\mathbf{b}^{\mathsf{T}} f(\mathbf{A})\mathbf{b}. \tag{7.1}$$

We are particularly interested in such approximations because, as discussed in Section 7.3, they are useful in approximating the trace of a matrix function

$$\operatorname{tr}(f(\mathbf{A})) = \sum_{i=1}^{d} f(\lambda_i), \tag{7.2}$$

also referred to as a spectral sum.

Applications of spectral sums include characterizing the degree of protein folding in biology [Est00], studying the thermodynamics of spin systems in quantum physics and chemistry [Wei+06; SS10; SRS20; Jin+21], benchmarking quantum devices in quantum information theory [Joz94], maximum likelihood estimation in statistics [BP99; PL04], designing better public transit in urban planning [BS22; Wan+21], and finding triangle counts and other structure in network science [Avr10; DBB19; BB20].

## 7.1 Lanczos quadrature

In Chapter 6, we described a Lanczos-based method for approximating $f(\mathbf{A})\mathbf{b}$ called the Lanczos method for matrix function approximation (Lanczos-FA). A similar approach, which we call the Lanczos method for quadratic form approximation (Lanczos-QF), can be used to approximate the quadratic form $\mathbf{b}^{\mathsf{T}} f(\mathbf{A})\mathbf{b}$.

**Definition 7.1.** The $k$-th Lanczos-QF approximation to $\mathbf{b}^{\mathsf{T}} f(\mathbf{A})\mathbf{b}$ is

$$\operatorname{lan-QF}_k(f) = \operatorname{lan-QF}_k(f; \mathbf{A}, \mathbf{b}) := \|\mathbf{b}\|^2 \mathbf{e}_1^{\mathsf{T}} f(\mathbf{T}_k)\mathbf{e}_1.$$

Note that $\operatorname{lan-QF}_k(f) = \mathbf{b}^{\mathsf{T}} \operatorname{lan-FA}_k(f)$. However, it is typically more useful to think of $\operatorname{lan-QF}_k(f)$ in terms of quadrature. In particular, unlike the Lanczos-FA iterate, the Lanczos-QF iterate can be computed without every looking at the Lanczos basis $\mathbf{Q}_k$.

**7.1.1 Exactness and error bounds** Similar to Lanczos-FA, Lanczos-QF is exact for low-degree polynomials.

**Theorem 7.2.** Suppose $\|\mathbf{b}\|^2 = 1$. Let $\psi_k(x; \mathbf{A}, \mathbf{b})$ be the $k$-point Gaussian quadrature for $\psi(x; \mathbf{A}, \mathbf{b})$. Then

$$\text{lan-QF}_k(f) = \int f(x)\psi_k(x; \mathbf{A}, \mathbf{b})\mathrm{d}x.$$

*Proof.* This is an immediate consequence of Theorem 3.1 the relationship between the Lanczos algorithm and orthogonal polynomials and the definition of Gaussian quadrature (Definition 3.4). ∎

As an immediate consequence of Theorem 7.2 and the exactness of Gaussian quadrature rules (Theorem 3.5), we have the following.[1]

**Corollary 7.3.** Let $p(x)$ be a polynomial with $\deg(p) < 2k$. Then

$$\text{lan-QF}_k(p) = \mathbf{b}^\mathsf{T} p(\mathbf{A})\mathbf{b}.$$

Therefore, analogous to Corollary 6.4 we have an error bound

**Corollary 7.4.** The Lanczos-QF iterate satisfies

$$|\mathbf{b}^\mathsf{T} f(\mathbf{A})\mathbf{b} - \text{lan-QF}_k(f)| \leq 2\|\mathbf{b}\|^2 \min_{\deg(p)<2k} \|f - p\|_{\mathcal{I}}.$$

*Proof.* For any polynomial $p(x)$ with $\deg(p) < 2k$, Corollary 7.3 asserts $\text{lan-QF}_k(p) = \mathbf{b}^\mathsf{T} p(\mathbf{A})\mathbf{b}$. Recalling that $\text{lan-QF}_k(f) = \mathbf{b}^\mathsf{T}\mathbf{Q}_k f(\mathbf{T}_k)\mathbf{Q}_k^\mathsf{T}\mathbf{b}$, the triangle inequality gives

$$|\mathbf{b}^\mathsf{T} f(\mathbf{A})\mathbf{b} - \text{lan-QF}_k(f)| \leq |\mathbf{b}^\mathsf{T} e(\mathbf{A})\mathbf{b}| + |\mathbf{b}^\mathsf{T}\mathbf{Q}_k e(\mathbf{T}_k)\mathbf{Q}_k^\mathsf{T}\mathbf{b}|, \tag{7.3}$$

where $e(x) := f(x) - p(x)$. Using the definition of the operator norm and the fact that $\Lambda \subset \mathcal{I}$,

$$|\mathbf{b}^\mathsf{T} e(\mathbf{A})\mathbf{b}| \leq \|\mathbf{b}\|^2 \|e(\mathbf{A})\| = \|\mathbf{b}\|^2 \|e\|_\Lambda \leq \|\mathbf{b}\|^2 \|e\|_{\mathcal{I}}. \tag{7.4}$$

Similarly, since $\mathbf{Q}_k$ is orthonormal, and using Lemma 3.11 which implies that $\Lambda(\mathbf{T}_k) \subset \mathcal{I}$,

$$|\mathbf{b}^\mathsf{T}\mathbf{Q}_k e(\mathbf{T}_k)\mathbf{Q}_k^\mathsf{T}\mathbf{b}| \leq \|\mathbf{b}\|^2 \|e(\mathbf{T}_k)\| = \|\mathbf{b}\|^2 \|e\|_{\Lambda(\mathbf{T}_k)} \leq \|\mathbf{b}\|^2 \|e\|_{\mathcal{I}}. \tag{7.5}$$

[1] This can also be proved directly using Theorem 6.2; see Theorem 9.10.

Plugging (7.4) and (7.5) into (7.3) and minimizing over $p$ gives the result. ∎

## 7.2 Finite precision arithmetic

In Section 4.3 we described Theorem 4.4, which is effectively a bound for how well Lanczos-QF applies Chebyshev polynomials.

**Theorem 4.4** (informal; see [Kni96]). Suppose $\|\mathbf{A}\| = \|\mathbf{b}\| = 1$ and that Lanczos algorithm is run for $k$ iterations on a computer with relative machine precision $\epsilon_{\text{mach}} = O(\text{poly}(k)^{-1})$. Let $\psi(x) = \psi(x; \mathbf{A}, \mathbf{b})$ and $\overline{\psi}_k(x) = \psi(x; \overline{\mathbf{T}}_k, \mathbf{e}_1)$ be the Gaussian quadrature produced by Lanczos in finite precision arithmetic. Then, for all $n \leq 2k-1$,

$$\left| \int T_n(x)\psi(x)\,\mathrm{d}x - \int T_n(x)\overline{\psi}_k(x)\,\mathrm{d}x \right| = O(\text{poly}(k)\,\epsilon_{\text{mach}}).$$

The proof of Theorem 4.4 in [Kni96] is based on the bound Theorem 6.6 for Lanczos-FA described in Chapter 6 and the analysis of Paige described in Section 4.1. The analysis of Knizhnerman is more-or-less rederived using notation similar to this monograph in [CT24].

We can obtain a bound for Lanczos-QF in finite precision arithmetic analagous to Corollary 7.4. The argument is the same as we used to prove Corollary 6.7 for Lanczos-FA,

**Corollary 7.5** (informal; see [Kni96]). Suppose $\|\mathbf{A}\| = \|\mathbf{b}\| = 1$ and that Lanczos algorithm is run for $k$ iterations on a computer with relative machine precision $\epsilon_{\text{mach}} = O(\text{poly}(k)^{-1})$. Then, with $\overline{\mathcal{I}} = [-1-\eta_k, 1+\eta_k]$ where $\eta_k = O(\text{poly}(k)\,\epsilon_{\text{mach}})$,

$$|\mathbf{b}^\mathsf{T} f(\mathbf{A})\mathbf{b} - \mathbf{e}_1^\mathsf{T} f(\overline{\mathbf{T}}_k)\mathbf{e}_1| = O(k) \min_{\deg(p)<2k} \|f-p\|_{\overline{\mathcal{I}}} + O(\|f\|_{\overline{\mathcal{I}}}\,\text{poly}(k)\,\epsilon_{\text{mach}}).$$

## 7.3 Stochastic trace estimation

Fix a matrix $\mathbf{M}$ and independently draw a random vector $\mathbf{b}$ for which $\mathbb{E}[\mathbf{b}\mathbf{b}^\mathsf{T}] = d^{-1}\mathbf{I}$ for instance, each entry of $\mathbf{b}$ could be an independent normal random variable with variance $d^{-1}$ or an independent $\pm d^{-1/2}$ random variable. Then, using the cyclic property of the trace and the linearity of the expectation, we have that

$$\mathbb{E}\big[\mathbf{b}^\mathsf{T}\mathbf{M}\mathbf{b}\big] = \mathbb{E}\big[\text{tr}(\mathbf{b}^\mathsf{T}\mathbf{M}\mathbf{b})\big] = \mathbb{E}\big[\text{tr}(\mathbf{M}\mathbf{b}\mathbf{b}^\mathsf{T})\big] = \text{tr}(\mathbf{M}\mathbb{E}[\mathbf{b}\mathbf{b}^\mathsf{T}]) = d^{-1}\,\text{tr}(\mathbf{M}). \quad (7.6)$$

The history of estimators of this form is somewhat murky. In numerical analysis and theoretical computer science, such estimators are often attributed to Hutchinson [Hut89], particularly when $\mathbf{b}$ has iid $\pm d^{-1/2}$ entries [AT11; RA14; Mey+21]. However, the abstract of [Hut89] itself cites Girard [Gir87], who also used the method for estimating the trace of a matrix and provided a theoretical analysis. In computational physics, the use of random vectors to produce empirical averages of quantities of interest has been used since at least the mid 1970s [Alb+75; WW76; WW77; RV89]; see [Jin+21] for a review. More broadly, expressions like (7.6) have close relation to the concept of *typicality* in quantum physics and appear as early as the late 1920s [Sch27; Neu29]; see [Gol+10] for a review.

Let $\mathbb{S}^{d-1} = \{\mathbf{x} \in \mathbb{R}^d : \|\mathbf{x}\| = 1\}$ be the unit hypersphere in $\mathbb{R}^d$. Suppose $\mathbf{b} = \mathbf{g}/\|\mathbf{g}\|$, where the entries of $\mathbf{g} \in \mathbb{R}^d$ are independent standard normal random variables. Then $\mathbf{b} \sim \mathrm{Unif}(\mathbb{S}^{d-1})$, the uniform distribution on $\mathbb{S}^{d-1}$. It is well-known that $\mathbf{b}^\mathsf{T}\mathbf{M}\mathbf{b}$ concentrates strongly about its expectation value $d^{-1}\operatorname{tr}(\mathbf{M})$.

**Theorem 7.6.** There is a universal constant $C > 0$ such that for any $\delta \in (0,1)$ and $d \times d$ matrix $\mathbf{M}$ the following holds. Let $\mathbf{b}_1, \dots, \mathbf{b}_m \sim \mathrm{Unif}(\mathbb{S}^{d-1})$ be independent samples. Then

$$\mathbb{P}\left[\left|\frac{1}{d}\operatorname{tr}(\mathbf{M}) - \frac{1}{m}\sum_{\ell=1}^{m} \mathbf{b}_\ell^\mathsf{T}\mathbf{M}\mathbf{b}_\ell\right| > \sqrt{\frac{\|\mathbf{M}\|^2 \log(2/\delta)}{Cdm}}\right] \le \delta.$$

*Proof.* We follow the approach of [PSW06]. Define $F(\mathbf{x}) := \mathbf{x}^\mathsf{T}\mathbf{M}\mathbf{x}$ and observe that

$$\begin{aligned} |F(\mathbf{x}) - F(\mathbf{y})| &= |\mathbf{x}^\mathsf{T}\mathbf{M}\mathbf{x} - \mathbf{y}^\mathsf{T}\mathbf{M}\mathbf{y}| && (7.7)\\ &= \frac{1}{2}|(\mathbf{x}+\mathbf{y})^\mathsf{T}\mathbf{M}(\mathbf{x}-\mathbf{y}) + (\mathbf{x}-\mathbf{y})^\mathsf{T}\mathbf{M}(\mathbf{x}+\mathbf{y})| && (7.8)\\ &\le \frac{1}{2}|(\mathbf{x}+\mathbf{y})^\mathsf{T}\mathbf{M}(\mathbf{x}-\mathbf{y})| + \frac{1}{2}|(\mathbf{x}-\mathbf{y})^\mathsf{T}\mathbf{M}(\mathbf{x}+\mathbf{y})|. && (7.9) \end{aligned}$$

We also have that

$$|(\mathbf{x}+\mathbf{y})^\mathsf{T}\mathbf{M}(\mathbf{x}-\mathbf{y})| \le \|\mathbf{x}+\mathbf{y}\|\|\mathbf{M}\|\|\mathbf{x}-\mathbf{y}\| \le (\|\mathbf{x}\| + \|\mathbf{y}\|)\|\mathbf{M}\|\|\mathbf{x}-\mathbf{y}\|, \quad (7.10)$$

and an analogous bound for $|(\mathbf{x}-\mathbf{y})^\mathsf{T}\mathbf{M}(\mathbf{x}+\mathbf{y})|$. Since $\|\mathbf{x}\| = \|\mathbf{y}\| = 1$ and $\|\mathbf{M}^\mathsf{T}\| = \|\mathbf{M}\|$ we therefore find that $F(\mathbf{x})$ is $2\|\mathbf{M}\|$-Lipshitz on $\mathbb{S}^{d-1}$.

By Lévy's lemma (see for instance [Ver18, Theorem 5.1.4]), there is a universal constant $C$ such that, if $\mathbf{b} \sim \mathrm{Unif}(\mathbb{S}^{d-1})$,

$$\mathbb{P}\big[|F(\mathbf{b}) - \mathbb{E}[F(\mathbf{b})]| > \varepsilon\big] \le 2\exp\left(\frac{-Cd\varepsilon^2}{(2\|\mathbf{M}\|)^2}\right).$$

The final result holds by a standard bound for the average of iid copies of a sub-Gaussian random variable (see for instance [Ver18, §2.5 and Theorem 2.6.2]) and relabeling $C$. ■

For simplicity we consider only the case of uniform unit vectors, but other choices of distribution, such as random sign vectors, are also common [MT20, §4]. In addition, we note that potentially sharper bounds, with $\|\mathbf{M}\|^2$ replaced by $\|\mathbf{M}\|_\mathsf{F}^2/d$, are known for many distributions [CK21; Mey+21; Epp24]. Such details are beyond the scope of this monograph.

## 7.4 Stochastic Lanczos quadrature

It is natural to combine Lanczos quadrature with stochastic trace estimation to obtain a method for approximating $\operatorname{tr}(f(\mathbf{A}))$. The resulting algorithm is commonly called stochastic Lanczos quadrature (SLQ) [BFG96; BG96; UCS17].

**Definition 7.7.** Let $\mathbf{b}_1, \ldots, \mathbf{b}_m \sim \operatorname{Unif}(\mathbb{S}^{d-1})$ be independent samples. The $(k, m)$-th SLQ approxiamtion to $d^{-1}\operatorname{tr}(f(\mathbf{A}))$ is

$$\operatorname{SLQ}_{k,m}(f) = \operatorname{SLQ}_{k,m}(f; \mathbf{A}) := \frac{1}{m}\sum_{\ell=1}^{m} \operatorname{lan-QF}_k(f; \mathbf{A}, \mathbf{b}_\ell).$$

We can use the triangle inequality to decompose the SLQ error into a term which accounts for the statistical noise from the stochastic trace estimator and a term which accounts for errors made by Lanczos quadrature. Specifically, we have that

$$\begin{aligned}
&\left|\frac{1}{d}\operatorname{tr}(f(\mathbf{A})) - \operatorname{SLQ}_{k,m}(f)\right| \\
&\qquad \leq \underbrace{\left|\frac{1}{d}\operatorname{tr}(f(\mathbf{A})) - \frac{1}{m}\sum_{\ell=1}^{m}\mathbf{b}_\ell^\mathsf{T} f(\mathbf{A})\mathbf{b}_\ell\right|}_{\text{trace estimator noise}} + \frac{1}{m}\sum_{\ell=1}^{m}\underbrace{\left|\mathbf{b}_\ell^\mathsf{T} f(\mathbf{A})\mathbf{b}_\ell - \operatorname{lan-QF}_k(f; \mathbf{A}, \mathbf{b}_\ell)\right|}_{\text{Lanczos quadrature error}}. \quad (7.11)
\end{aligned}$$

The first term is made small by increasing $m$, while the second term is made small by increasing $k$. This allows us to obtain theoretical guarantees various values of $k$ and $m$. A prototypical theoretical guarantee is the following.

**Theorem 7.8.** There is a universal constant $C > 0$ such that, for any function $f(x)$ and failure probability $\delta$, the $(k, m)$-th SLQ approximation to $d^{-1}\operatorname{tr}(f(\mathbf{A}))$ satisfies

$$\mathbb{P}\left[\left|d^{-1}\operatorname{tr}(f(\mathbf{A})) - \mathrm{SLQ}_{k,m}(f)\right| > \sqrt{\frac{\|f\|_{\mathcal{I}}^2 \log(2/\delta)}{Cdm}} + 2 \min_{\deg(p)<2k} \|f-p\|_{\mathcal{I}}\right] \leq \delta.$$

*Proof.* First, using Theorem 7.6 and that $\|f(\mathbf{A})\| \leq \|f\|_{\mathcal{I}}$ we have

$$\mathbb{P}\left[\left|\frac{1}{d}\operatorname{tr}(f(\mathbf{A})) - \frac{1}{m}\sum_{\ell=1}^{m} \mathbf{b}_\ell^\mathsf{T} f(\mathbf{A})\mathbf{b}_\ell\right| > \sqrt{\frac{\|f\|_{\mathcal{I}}^2 \log(2/\delta)}{Cdm}}\right] \leq \delta. \tag{7.12}$$

Next, using Corollary 7.4 and the fact that $\|\mathbf{b}_\ell\| = 1$, we have that

$$\left|\mathbf{b}_\ell^\mathsf{T} f(\mathbf{A})\mathbf{b}_\ell - \text{lan-QF}_k(f;\mathbf{A},\mathbf{b}_\ell)\right| \leq 2 \min_{\deg(p)<2k} \|f-p\|_{\mathcal{I}}. \tag{7.13}$$

Finally, using (7.11) gives the result. ■

Bounds like Theorem 7.8 were first derived in [UCS17]. In particular, [UCS17] gives explicit bounds for $m$ and $k$ required to get $\varepsilon$ error with probability $1-\delta$. Theorem 7.8 already reveals how we must set $m$ (and in particular $m = O(\varepsilon^{-2})$. To get an explicit bound for $k$, [UCS17] considers particular cases of functions $f(x)$; e.g. $f(x)$ analytic on a Bernstein ellipse containing $\mathcal{I}$. In these cases one can derive specific bounds for the the number of Lanczos steps $k$ required using classic approximation theory. For instance, for bounded functions analytic on some Bernstein ellipse containing $\mathcal{I}$, it suffices to take $k = O(\log(1/\varepsilon))$ [Tre19].

### 7.4.1 Variance reduction

By the linearity of the trace we have that

$$\operatorname{tr}(\mathbf{A}) = \operatorname{tr}(\tilde{\mathbf{A}}) + \operatorname{tr}(\mathbf{A}-\tilde{\mathbf{A}}), \tag{7.14}$$

for any matrix $\tilde{\mathbf{A}}$. If one can obtain a matrix $\tilde{\mathbf{A}}$ for which the exact trace can be computed efficiently, then we can apply the stochastic trace estimator to approximate the residual term $\operatorname{tr}(\mathbf{A}-\tilde{\mathbf{A}})$. If $\|\mathbf{A}-\tilde{\mathbf{A}}\| \ll \|\mathbf{A}\|$, then this can result in significantly lower variance than trying to apply stochastic trace estimation to $\operatorname{tr}(\mathbf{A})$ directly. For the regular trace estimation problem, such an observation has been made several times [Gir87; GSO17; Mey+21]. Notably, [Mey+21] introduces Hutch++, which uses sketching-based low-rank approximation to obtain $\tilde{\mathbf{A}}$. This provably reduces the number of matrix-vector products with $\mathbf{A}$ from $O(\varepsilon^{-2})$ to $O(\varepsilon^{-1})$; see also [PCK22; ETW24].

One can naturally implement Hutch++ to approximate $\operatorname{tr}(f(\mathbf{A}))$ by using black-box methods like Lanczos-FA from Chapter 6 to approximate matrix-vector products with $f(\mathbf{A})$. However, it is often better to open up the black-box, and look more closely at how the variance reduction and Krylov subspace methods interact [CH23; PK23; PMM25; PCM26].

## TL;DR

Methods for approximating traces and quadratic forms of matrix functions are closely related to quadrature. As with Lanczos-FA, Lanczos-QF works fine in finite precision arithmetic.

# 8 Spectrum approximation

Computing the spectral density $\varphi(x;\mathbf{A})$ is equivalent to computing the spectrum of $\mathbf{A}$ and is therefore infeasible in many settings. However, approximations to $\varphi(x)$ are also useful in that they provide a global picture of the spectrum of $\mathbf{A}$. Such coarse grained approximations are used in electronic structure computations and other tasks in physics[1] [Wei+06; Jin+21], probing the behavior of neural networks in machine learning [GKX19; Pap19; GWG19; Yao+20], understanding the structure of networks in spectral graph theory [KV17; BKM22], load balancing modern parallel eigensolvers in numerical linear algebra [Pol09; Li+19], and computing the product of matrix functions with vectors [Fan+19].

We remark that spectrum approximation is closely related to spectral sum approximation, discussed in Chapter 7. In particular,

$$d^{-1}\operatorname{tr}(f(\mathbf{A})) = \int f(x)\varphi(x;\mathbf{A})\mathrm{d}x. \tag{8.1}$$

This connection is will be essential to analyzing methods for spectrum approximation.

## 8.1 Stochastic Lanczos Quadrature

Using the same arguments as in Section 7.3, is not hard to see that if $\mathbb{E}[\mathbf{b}\mathbf{b}^\mathsf{T}] = d^{-1}\mathbf{I}$, then for each value $x$

$$\mathbb{E}[\psi(x;\mathbf{A},\mathbf{b})] = \varphi(x;\mathbf{A}). \tag{8.2}$$

While computing $\psi(x;\mathbf{A},\mathbf{b})$ is not easy, we can compute the Gaussian quadrature approximation $\psi_k(x;\mathbf{A},\mathbf{b})$ to $\psi(x;\mathbf{A},\mathbf{b})$. Averaging repeated copies of $\psi_k(x;\mathbf{A},\mathbf{b})$ to reduce statistical noise yields the SLQ approximation to the spectral density function.

**Definition 8.1.** Let $\mathbf{b}_1,\ldots,\mathbf{b}_m \sim \operatorname{Unif}(\mathbb{S}^{d-1})$ be independent vectors drawn from the unit hypersphere. For each $\ell = 1,\ldots,m$, denote by $\psi_k(x;\mathbf{A},\mathbf{b}_\ell)$ be the $k$-point Guassian quadrature approximation to $\psi(x;\mathbf{A},\mathbf{b}_\ell)$. The $(k,m)$

[1]In physics, $\varphi(x)$ is often called the density of states (DOS).

SLQ approximation $\psi_{k,m}^{\mathrm{SLQ}}(x)$ to $\varphi(x;\mathbf{A})$ is defined as

$$\psi_{k,m}^{\mathrm{SLQ}}(x) = \psi_{k,m}^{\mathrm{SLQ}}(x;\mathbf{A}) := \frac{1}{m}\sum_{\ell=1}^{m}\psi_k(x;\mathbf{A},\mathbf{b}_\ell)$$

Using Theorem 7.2 we see that the SLQ trace estimator is compatible with the SLQ spectral density estimator.

**Lemma 8.2.** The $(k,m)$ SLQ approximation $\mathrm{SLQ}_{k,m}(f)$ to $\mathrm{tr}(f(\mathbf{A}))$ is related to the $(k,m)$-th SLQ approximation $\psi_{k,m}^{\mathrm{SLQ}}(x)$ to $\varphi(x;\mathbf{A})$ by

$$\mathrm{SLQ}_{k,m}(f) = \int f(x)\psi_{k,m}^{\mathrm{SLQ}}(x)\mathrm{d}x. \tag{8.3}$$

In order to obtain theoretical guarantees for SLQ spectrum approximation, we must establish a measure of distance between two densities. A common way of defining the distance between to densities is the Wasserstein distance. A visual explication is given in Figure 8.1.

**Definition 8.3.** Let $\mu_1(x)$ and $\mu_2(x)$ be two probability density functions and let $M_1(x)$ and $M_2(x)$ be their respective cumulative distribution functions. The Wasserstein distance between $\mu_1(x)$ and $\mu_2(x)$, denoted $\mathrm{d_W}(\mu_1,\mu_2)$, is defined by

$$\mathrm{d_W}(\mu_1,\mu_2) = \int |M_1(x) - M_2(x)|\mathrm{d}x.$$

SLQ satisfies the following Wasserstein convergence guarantee.

**Theorem 8.4.** Suppose

$$m \geq \frac{\log^2(\varepsilon^{-1})\log(\varepsilon^{-1}\delta^{-1})}{d\varepsilon^2}, \qquad k \geq \frac{1}{\varepsilon}.$$

Then, the $(k,m)$-SLQ approximation to $\varphi(\,\cdot\,;\mathbf{A})$ satisfies

$$\mathbb{P}\left[\mathrm{d_W}\left(\varphi(\,\cdot\,;\mathbf{A}),\psi_{k,m}^{\mathrm{SLQ}}\right) > |\lambda_{\max} - \lambda_{\min}|\,O(\varepsilon)\right] \leq \delta.$$

In particular, note that if $\varepsilon \gg d^{-0.499}$, then we can instantiate Theorem 8.4 with $m = 1$ (i.e. a single random test vector), in which case one can obtain a $O(\epsilon)$-accurate Wasserstein approximation with $\varepsilon^{-1}$ matrix-vector products.

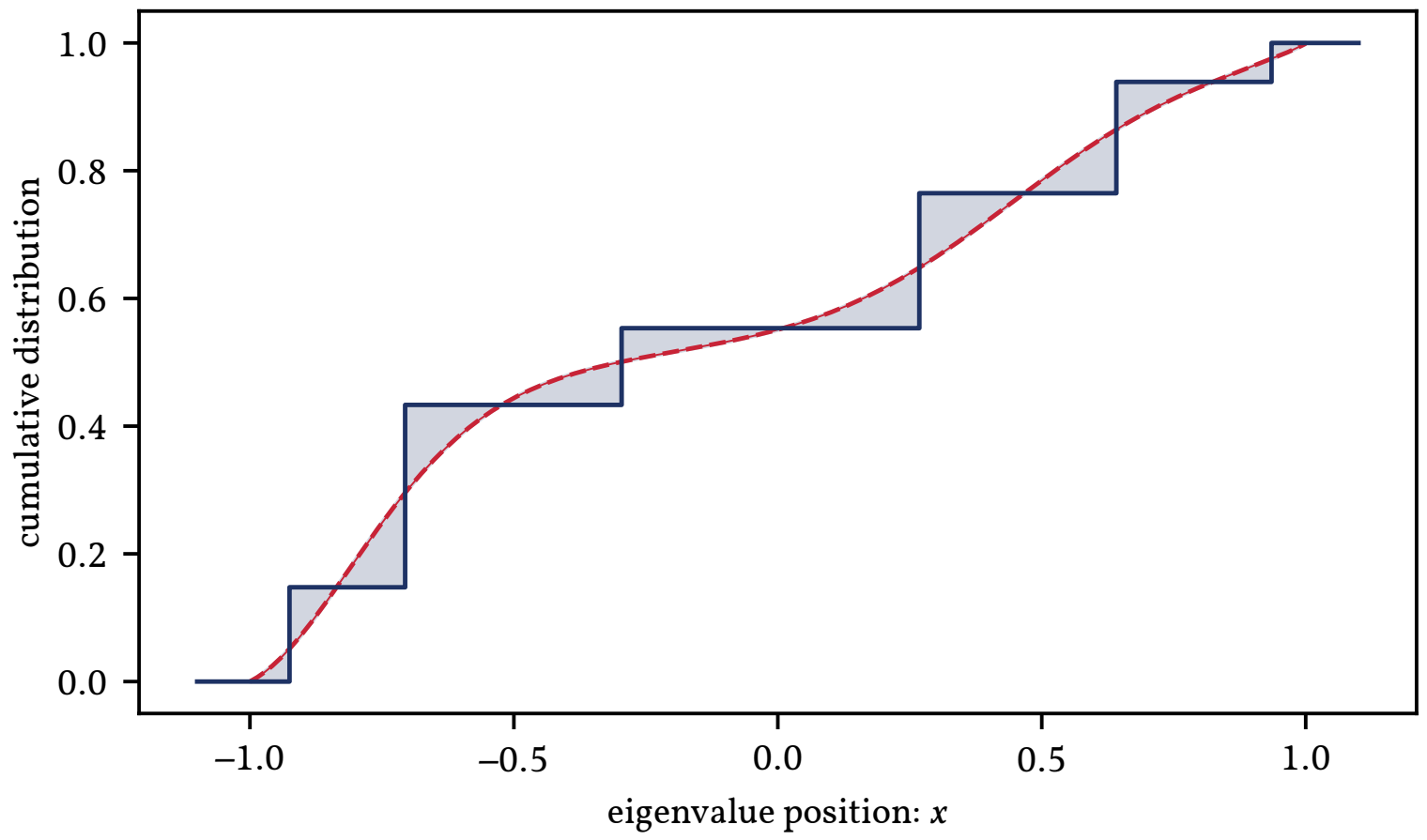

**Figure 8.1:** Cumulative distribution functions $M_1(x)$ ( ) and $M_2(x)$ ( ) corresponding to densities $\mu_1(x)$ and $\mu_2(x)$. The Wasserstein distance $\mathrm{d_W}(\mu_1, \mu_2)$ is the total area of the shaded region between the two curves.

To the best of our knowledge, the first theoretical bounds for SLQ were established in [CTU21]. Our proof technique here more closely follows the approach of [BKM22] for the kernel polynomial method described in Section 8.2. A related, but more complicated, algorithm which gives an $\varepsilon$-accurate Wasserstein approximation with roughly $\varepsilon^{-1}$ matrix-vector products for any $\varepsilon \gg d^{-1}$ is described in [Mus+25].

### 8.1.1 Proof of Theorem 8.4

It is well-known (see for instance [Vil09, Remark 6.5]) that the Wasserstein distance is equivalently defined by

$$\mathrm{d_W}(\mu_1, \mu_2) = \sup\left\{ \int f(x)\big(\mu_1(x) - \mu_2(x)\big)\,\mathrm{d}x : f(x) \text{ 1-Lipschitz} \right\}. \tag{8.4}$$

This characterization will be useful in our analysis.

We begin by stating Jackson's theorem [Jac12], which asserts that all Lipschitz functions can be approximated by a low-degree polynomial with nice Chebyshev coefficients; see [Riv81, Theorem 1.4] and [Tre19, Theorem 7.1] for details.

**Theorem 8.5.** Suppose $f(x)$ is 1-Lipschitz. There exists[2] a polynomial

$$p(x) = c_0 T_0(x) + 2\sum_{n=1}^{2k-1} c_n T_n(x)$$

such that

$$|c_n| = O(n^{-1}), \qquad \|f - p\|_{[-1,1]} = O(k^{-1}).$$

Next, we show that if two densities have similar Chebyshev moments, then they are close in Wasserstein distance.[3]

**Corollary 8.6.** Suppose $\mu_1(x)$ and $\mu_2(x)$ are two probability density functions supported on $[-1,1]$ satisfying, for all $n < 1/\varepsilon$,

$$\left|\int T_n(x)(\mu_1(x)-\mu_2(x))\,\mathrm{d}x\right| = O\Big(\frac{\varepsilon}{\log(\varepsilon^{-1})}\Big).$$

Then

$$\mathrm{d_W}(\mu_1,\mu_2) = O(\varepsilon).$$

*Proof.* Set $k = 1/\varepsilon$. Let $f(x)$ be an arbitrary 1-Lipschitz function and $p(x)$ as in Theorem 8.5. By the triangle inequality,

$$\left|\int f(x)(\mu_1(x)-\mu_2(x))\,\mathrm{d}x\right| \\ \leq \left|\int (f(x)-p(x))(\mu_1(x)-\mu_2(x))\,\mathrm{d}x\right| + \left|\int p(x)(\mu_1(x)-\mu_2(x))\,\mathrm{d}x\right|. \quad (8.5)$$

Since $\mu_1(x)$ and $\mu_2(x)$ are both probability density functions, Theorem 8.5 implies

$$\left|\int (f(x)-p(x))(\mu_1(x)-\mu_2(x))\,\mathrm{d}x\right| \leq \int |f(x)-p(x)||\mu_1(x)-\mu_2(x)|\,\mathrm{d}x \quad (8.6)$$

$$\leq \|f-p\|_{[-1,1]}\int |\mu_1(x)|+|\mu_2(x)|\,\mathrm{d}x \quad (8.7)$$

$$\leq 2\|f-p\|_{[-1,1]} = O(k^{-1}) = O(\varepsilon). \quad (8.8)$$

Next, using the expansion of $p(x)$, and again noting that $\mu_1(x)$ and $\mu_2(x)$ are both probability density functions, then using Theorem 8.5 and the fact that $1+1/2+\cdots+1/n = O(\log(n))$, we have that

$$\left|\int p(x)(\mu_1(x)-\mu_2(x))\,\mathrm{d}x\right| = \left|2\sum_{n=1}^{2k-1} c_n \int T_n(x)(\mu_1(x)-\mu_2(x))\,\mathrm{d}x\right|. \quad (8.9)$$

$$\leq 2\sum_{n=1}^{2k-1}|c_n|\left|\int T_n(x)(\mu_1(x)-\mu_2(x))\,\mathrm{d}x\right| \quad (8.10)$$

$$\leq 2\sum_{n=1}^{2k-1}\frac{1}{n}\frac{\varepsilon}{\log(k)} = O(\varepsilon). \quad (8.11)$$

[2] In fact, the polynomial can be obtained constructively.

[3] This bound is more-or-less sharp; it's not hard to show there exist distributions with matching moments through degree $1/\varepsilon$ with Wasserstein distances $\Omega(\varepsilon)$ [KV17; Jin+23].

Plugging (8.8) and (8.11) into (8.5) we have that

$$\left| \int f(x)(\mu_1(x) - \mu_2(x))\,\mathrm{d}x \right| = O(\varepsilon). \tag{8.12}$$

Hence, since $f(x)$ is arbitrary, using (8.4) gives the result. ■

The proof of proof of Theorem 8.4 then follows as an immediate result of the previous to facts.

*Proof of Theorem 8.4.* W.l.o.g., scale $\mathbf{A}$ so that $\mathcal{I} = [-1, 1]$. By Lemma 3.26 we have that $\|T_n\|_{\mathcal{I}} \leq 1$. Therefore, using Theorem 7.8, we have that

$$\forall n < 2k : \mathbb{P}\left[ |d^{-1}\operatorname{tr}(T_n(\mathbf{A})) - \mathrm{SLQ}_{k,m}(T_n)| > \sqrt{\frac{\log(4k/\delta)}{Cdm}} \right] \leq \frac{\delta}{2k}.$$

Thus, applying a union bound we therefore have that

$$\mathbb{P}\left[ \exists n < 2k : |d^{-1}\operatorname{tr}(T_n(\mathbf{A})) - \mathrm{SLQ}_{k,m}(T_n)| > \sqrt{\frac{\log(4k/\delta)}{Cdm}} \right] \leq \delta. \tag{8.13}$$

Our choice of $m$ allows us to instantiate Corollary 8.6, and our choice of $k$ (and removing the condition on $[\lambda_{\min}, \lambda_{\max}]$) then gives the result. ■

## 8.2 The Kernel Polynomial Method

Lanczos quadrature can be viewed as approximating $\psi(x; \mathbf{A}, \mathbf{b})$ by Gaussian quadrature. A common alternative is the kernel polynomial method (KPM) [Ski89; SR94; Sil+96; Wei+06].

**Definition 8.7.** Fix a reference density $\nu(x)$ with orthonormal polynomials $\{q_n(x)\}$. The $k$-th KPM approximation $\mu_k^{\mathrm{KPM}}(x)$ to $\mu(x)$ is the "density function"

$$\mu_k^{\mathrm{KPM}}(x) := \nu(x) \sum_{n=0}^{2k-1} \left( \int q_n(z)\mu(z)\,\mathrm{d}z \right) q_n(x).$$

A motivation for the definition of the KPM is the following. Expand $\mu(x)/\nu(x)$ in an orthogonal polynomial series

$$\frac{\mu(x)}{\nu(x)} = \sum_{n=0}^{\infty} \left( \int q_n(z) \frac{\mu(z)}{\nu(z)} \nu(z)\,\mathrm{d}z \right) q_n(x) = \sum_{n=0}^{\infty} \left( \int q_n(z)\mu(z)\,\mathrm{d}z \right) q_n(x). \tag{8.14}$$

So we obtain the approximating by truncating the series at $2k-1$ and multiplying up the $\nu(x)$. This also makes it clear that we need that the support of $\mu(x)$ is contained in the support of $\nu(x)$, otherwise $\mu(x)/\nu(x)$ will not be defined.

We note that when $\mu(x) = \psi(x; \mathbf{A}, \mathbf{b})$, then the coefficients in the KPM approximation can be computed by the identity

$$\int q_n(x)\mu(x)\mathrm{d}x = \mathbf{b}^\mathsf{T} q_n(\mathbf{A})\mathbf{b}. \tag{8.15}$$

Most commonly, the reference density $\nu(x)$ is taken to be a shifted and scaled version of the Chebyshev density $\mu_T(x)$ defined in (3.17). In this case, the coefficients (8.15) are often computed using an explicit Chebyshev recurrence (3.15); see for instance [Wei+06]. As we discuss in Section 8.2, Lanczos-based approaches are also possible.

As with the Gaussian quadrature approximation, the KPM approximation satisfies an exactness property for integrating polynomials.

**Lemma 8.8.** Suppose the support of $\mu(x)$ is contained in the support of $\nu(x)$. Then for any polynomial $p(x)$ of degree less than $2k$,

$$\int p(x)\mu_k^{\mathrm{KPM}}(x)\mathrm{d}x = \int p(x)\mu(x)\mathrm{d}x. \tag{8.16}$$

*Proof.* Fix $m < 2k$. By the orthonormality of $\{q_n(x)\}$ with respect to $\nu(x)$,

$$\begin{aligned}
\int q_m(x)\mu_k^{\mathrm{KPM}}(x)\mathrm{d}x &= \int \sum_{n=0}^{2k-1} \left( \int q_n(z)\mu(z)\mathrm{d}z \right) q_n(x) q_m(x)\nu(x)\mathrm{d}x && (8.17)\\
&= \sum_{n=0}^{2k-1} \left( \int q_n(z)\mu(z)\mathrm{d}z \right) \int q_n(x)q_m(x)\nu(x)\mathrm{d}x && (8.18)\\
&= \int q_m(z)\mu(z)\mathrm{d}z. && (8.19)
\end{aligned}$$

The result follows since $\{q_n(x)\}$ is a basis for polynomials of degree less than $2k$. ■

**Damping** Unlike the SLQ approximation, the KPM approximation could be negative at some values. This can sometimes be mitigated through the use of damping. Given damping coefficient $\{\rho_{n,k}\}$, the damped KPM approximation is defined as

$$\mu_k^{\mathrm{dKPM}}(x) := \nu(x) \sum_{n=0}^{2k-1} \rho_{n,k} \left( \int q_n(z)\mu(z)\mathrm{d}z \right) q_n(x). \tag{8.20}$$

Clearly the case $\rho_{n,k} = 1$ recovers the KPM approximation.

When $\nu(x) = \mu_T(x)$ and the support of $\mu(x)$ is contained in $[-1, 1]$, then Jackson's damping coefficients

$$\rho_{n,k} = (2k+1)^{-1}\left((2k-n+1)\cos\left(\frac{n\pi}{2k+1}\right) + \sin\left(\frac{n\pi}{2k+1}\right)\cot\left(\frac{\pi}{2k+1}\right)\right) \quad (8.21)$$

can be used to ensure that the resulting approximation is a actually a density function [Wei+06]. In this case, a bound similar to Corollary 7.4 can be obtained [BKM22; CTU25].

**Spectrum adaptivity** The reference density $\nu(x)$ must be chosen before the coefficients (8.15) can be computed. However, this does not require fixing the reference density before products with $\mathbf{A}$ are performed. In particular, as noted in [Che23], the Lanczos quadrature method from Section 7.4 can be used to compute the coefficients (8.15). This allows the reference density to be chosen adaptively based on the spectrum of $\mathbf{A}$, which can lead to better approximations in a number of settings.

## TL;DR

We can directly approximate the spectral density of a matrix using quadrature methods. This is closely related to trace estimation.

# 9 Block Methods

There are a number of situations in which we are simultaneously interested in the Krylov subspaces corresponding to a matrix $\mathbf{A}$ and multiple vectors $\mathbf{b}_1, \ldots, \mathbf{b}_m$. This motivates the definition of a block Krylov subspace.

**Definition 9.1.** The dimension $k$ block Krylov subspace generated by $\mathbf{A}$ and a matrix $\mathbf{B} = [\mathbf{b}_1, \ldots, \mathbf{b}_m]$ is defined as

$$\mathcal{K}_k(\mathbf{A}, \mathbf{B}) := \mathcal{K}_k(\mathbf{A}, \mathbf{b}_1) + \cdots + \mathcal{K}_k(\mathbf{A}, \mathbf{b}_m).$$

In particular, note that each Krylov subspace $\mathcal{K}_k(\mathbf{A}, \mathbf{b}_i)$ is contained in the block Krylov subspace $\mathcal{K}_k(\mathbf{A}, \mathbf{B})$.

The block Krylov subspace can equivalently be defined as

$$\mathcal{K}_k(\mathbf{A}, \mathbf{B}) = \operatorname{span}\{\mathbf{B}, \mathbf{A}\mathbf{B}, \ldots, \mathbf{A}^{k-1}\mathbf{B}\}, \tag{9.1}$$

where the span of a collection of matrices is interpreted as the span of all of the constituent columns of the matrices. Such a definition is more reminiscent of Definition 1.2 for standard Krylov subspaces.

## 9.1 Block Lanczos

The Lanczos algorithm naturally generalizes to the block setting [GU77]. An implementation of the block-Lanczos algorithm for obtaining an orthonormal basis for $\mathcal{K}_k(\mathbf{A}, \mathbf{B})$ is described in Algorithm 9.2.

**Algorithm 9.2** (Block Lanczos).

1: $\textsc{block-Lanczos}(\mathbf{A}, \mathbf{B}, k)$
2: $\quad \hat{\mathbf{Q}}_0, \hat{\mathbf{B}}_{-1} = \mathrm{QR}(\mathbf{B}), \hat{\mathbf{Q}}_{-1} = \mathbf{0}$
3: $\quad$ **for** $n = 0, 1, \ldots, k-1$ **do**
4: $\quad\quad \mathbf{Y}_{n+1} = \mathbf{A}\hat{\mathbf{Q}}_n - \hat{\mathbf{Q}}_{n-1}\hat{\mathbf{B}}_{n-1}^\mathsf{T}$
5: $\quad\quad \hat{\mathbf{A}}_n = \hat{\mathbf{Q}}_n^\mathsf{T}\mathbf{Y}_{n+1}$
6: $\quad\quad \mathbf{Z}_{n+1} = \mathbf{Y}_{n+1} - \hat{\mathbf{Q}}_n\hat{\mathbf{A}}_n$
7: $\quad\quad$ orthogonalize against $\hat{\mathbf{Q}}_0, \ldots, \hat{\mathbf{Q}}_n$ $\quad \triangleright$ optional

```
8:        Q̂_{n+1}, B̂_n = QR(Z_{n+1})        ▷ deflate to ensure orthogonality
9:    return {Q̂_n}, {Â_n}, {B̂_n}
```

Assuming the algorithm terminates successfully, the output satisfies a symmetric three term recurrence

$$\mathbf{A}\hat{\mathbf{Q}}_n = \hat{\mathbf{Q}}_{n-1}\hat{\mathbf{B}}_{n-1}^{\mathsf{T}} + \hat{\mathbf{Q}}_n\hat{\mathbf{A}}_n + \hat{\mathbf{Q}}_{n+1}\hat{\mathbf{B}}_n. \tag{9.2}$$

This recurrence can be written in matrix form as

$$\mathbf{A}\mathbf{Q}_k = \mathbf{Q}_k\mathbf{T}_k + \hat{\mathbf{Q}}_k\hat{\mathbf{B}}_{k-1}\mathbf{E}_k^{\mathsf{T}} \tag{9.3}$$

where $\mathbf{E}_k = \mathbf{e}_k \otimes \mathbf{I}_m$, where "$\otimes$" denotes the Kronecker product, and

$$\mathbf{Q}_k := \begin{bmatrix} | & | & & | \\ \hat{\mathbf{Q}}_0 & \hat{\mathbf{Q}}_1 & \cdots & \hat{\mathbf{Q}}_{k-1} \\ | & | & & | \end{bmatrix}, \quad \mathbf{T}_k := \begin{bmatrix} \hat{\mathbf{A}}_0 & \hat{\mathbf{B}}_0^{\mathsf{T}} & & \\ \hat{\mathbf{B}}_0 & \hat{\mathbf{A}}_1 & \ddots & \\ & \ddots & \ddots & \hat{\mathbf{B}}_{k-2}^{\mathsf{T}} \\ & & \hat{\mathbf{B}}_{k-2} & \hat{\mathbf{A}}_{k-1} \end{bmatrix}. \tag{9.4}$$

Here the $\hat{\mathbf{B}}_n$ are upper triangular, so $\mathbf{T}_k$ is of bandwidth at most $2m+1$.

Block Lanczos is substantially more difficult to implement in practice than than Lanczos. As with Lanczos, omitting the orthogonalization in Line 7 of Algorithm 9.2 will result in a loss of orthogonality in the columns of $\mathbf{Q}_k$. However, even if full reorthogonalization is used, the blocks $\mathbf{Z}_{n+1}$ may themselves become rank-deficient. This is because when $m > 1$, $\dim \mathcal{K}_k(\mathbf{A}, \mathbf{B}) < km$ does not imply that $\mathcal{K}_k(\mathbf{A}, \mathbf{B})$ is an invariant subspace of $\mathbf{A}$. As such, we do not necessarily want to terminate the block-Lanczos algorithm as soon as $\dim \mathcal{K}_k(\mathbf{A}, \mathbf{B}) < km$. However, this means that *even in exact arithmetic*, QR algorithm used in in Line 8 of Algorithm 9.2 must return $\hat{\mathbf{Q}}_{n+1}$ with rank equal to that of $\mathbf{Z}_{n+1}$. This is complicated further in finite precision arithmetic, where the algorithm must identify the numerical rank of $\mathbf{Z}_{n+1}$ and deflate appropriately.

**9.1.1 Finite precision arithmetic** The behavior of block Lanczos in finite precision arithmetic is much less understood than the behavior of the standard Lanczos algorithm. For instance, as far as we are aware, there is no Paige style analysis (see Section 4.1) which guarantees that a perturbed version of (9.3) holds in finite precision arithmetic. Even so, good implementations (which deflate when necessary) seem to work well in practice, producing outputs which nearly satisfy (9.3).

## 9.2 Block CG, block Lanczos-FA, and block Lanczos-QF

The conjugate gradient, Lanczos-FA, and Lanczos-QF algorithms have natural generalizations to the block setting [OLe80; Cab97; Fro+17, etc.].

### 9.2.1 Block CG

We can define a block conjugate gradient iterate analogously to the standard conjugate gradient iterate.

**Definition 9.3.** The $k$-th block conjugate gradient iterate $\mathbf{X}_k^{\mathrm{CG}}$ is given by

$$\mathbf{X}_k^{\mathrm{CG}} = \mathbf{X}_k^{\mathrm{CG}}(\mathbf{A},\mathbf{B}) := \mathbf{Q}_k\mathbf{T}_k^{-1}\mathbf{E}_1\hat{\mathbf{B}}_{-1}.$$

Define the matrix $\mathbf{A}$-norm $\|\cdot\|_{\mathbf{A}}$ by $\|\mathbf{X}\|_{\mathbf{A}} = \operatorname{tr}(\mathbf{X}^\mathsf{T}\mathbf{A}\mathbf{X})^{1/2} = \|\mathbf{A}^{1/2}\mathbf{X}\|_\mathsf{F}$. A similar proof gives the block version of Theorem 5.2.

**Theorem 9.4.** If $\mathbf{A}$ is positive definite, the block CG iterate satisfies the formula

$$\mathbf{X}_k^{\mathrm{CG}} = \operatorname*{argmin}_{\operatorname{range}(\mathbf{X})\subseteq\mathcal{K}_k(\mathbf{A},\mathbf{B})} \|\mathbf{A}^{-1}\mathbf{B}-\mathbf{X}\|_{\mathbf{A}}.$$

The optimality of the block-CG iterate over the block Krylov subspace implies that block-CG performs at least as well as CG.

**Corollary 9.5.** Let $[\mathbf{X}_k^{\mathrm{CG}}]_\ell$ denote the $\ell$-th column of $\mathbf{X}_k^{\mathrm{CG}}$. Then

$$\sum_{\ell=1}^{m}\|\mathbf{A}^{-1}\mathbf{b}_\ell-[\mathbf{X}_k^{\mathrm{CG}}]_\ell\|_{\mathbf{A}}^2 \leq \sum_{\ell=1}^{m}\|\mathbf{A}^{-1}\mathbf{b}_\ell-\mathbf{x}_k^{\mathrm{CG}}(\mathbf{A},\mathbf{b}_\ell)\|_{\mathbf{A}}^2$$

*Proof.* Using Theorems 5.2 and 9.4, we have that

$$\begin{aligned}
\|\mathbf{A}^{-1}\mathbf{b}_\ell-[\mathbf{X}_k^{\mathrm{CG}}]_\ell\|_{\mathbf{A}} &= \min_{\mathbf{x}\in\mathcal{K}_k(\mathbf{A},\mathbf{B})}\|\mathbf{A}^{-1}\mathbf{b}_\ell-\mathbf{x}\|_{\mathbf{A}} &(9.5)\\
&\leq \min_{\mathbf{x}\in\mathcal{K}_k(\mathbf{A},\mathbf{b}_\ell)}\|\mathbf{A}^{-1}\mathbf{b}_\ell-\mathbf{x}\|_{\mathbf{A}} &(9.6)\\
&= \|\mathbf{A}^{-1}\mathbf{b}_\ell-\mathbf{x}_k^{\mathrm{CG}}\|_{\mathbf{A}}. &(9.7)
\end{aligned}$$

The result then follows. ■

One can use Corrolaries 5.3 and 5.4 to derive bounds in terms of polynomial approximation. One can also extend Corollary 9.5 to related functions such as Stieltjes functions at the cost of condition-number factors.

### 9.2.2 Block Lanczos-FA

The block Lanczos-FA iterate can be defined analogously to the standard Lanczos-FA iterate (Definition 6.1).

**Definition 9.6.** The $k$-th block Lanczos-FA approximation to $f(\mathbf{A})\mathbf{B}$ is

$$\text{lan-FA}_k(f) = \text{lan-FA}_k(f;\mathbf{A},\mathbf{B}) := \mathbf{Q}_k f(\mathbf{T}_k)\mathbf{E}_1\hat{\mathbf{B}}_{-1}.$$

Analogously to Theorem 6.2 for Lanczos-FA, block Lanczos-FA is exact for low-degree polynomials.[1]

**Theorem 9.7.** Suppose $\deg(p) < k$. Then,

$$\text{lan-FA}_k(p) = p(\mathbf{A})\mathbf{B}.$$

*Proof.* By linearity, it suffices to verify that $\text{lan-FA}_k(x^n; \mathbf{A}, \mathbf{B}) = \mathbf{A}^n\mathbf{B}$ for all $n < k$. Observe that $\mathbf{Q}_k\mathbf{Q}_k^\mathsf{T}$ is the orthogonal projector onto $\mathcal{K}_k(\mathbf{A}, \mathbf{B})$. Hence, since $\mathbf{B}, \mathbf{A}\mathbf{B}, \ldots, \mathbf{A}^n\mathbf{B} \in \mathcal{K}_k(\mathbf{A}, \mathbf{B})$,

$$\begin{aligned} \mathbf{A}^n\mathbf{B} &= \mathbf{Q}_k\mathbf{Q}_k^\mathsf{T}\mathbf{A}^n\mathbf{B} && (9.8)\\ &= \mathbf{Q}_k\mathbf{Q}_k^\mathsf{T}\mathbf{A}\mathbf{Q}_k\mathbf{Q}_k^\mathsf{T}\mathbf{A}^{n-1}\mathbf{B} && (9.9)\\ &= \mathbf{Q}_k\mathbf{T}_k\mathbf{Q}_k^\mathsf{T}\mathbf{A}^{n-1}\mathbf{B} && (9.10)\\ &\ \vdots && (9.11)\\ &= \mathbf{Q}_k\mathbf{T}_k^n\mathbf{Q}_k^\mathsf{T}\mathbf{B}. && (9.12) \end{aligned}$$

Now, recall $\mathbf{B} = \hat{\mathbf{Q}}_0\hat{\mathbf{B}}_{-1}$, where $\hat{\mathbf{Q}}_0$ is the first block of $\mathbf{Q}_k$. The result then follows from the fact that $\mathbf{Q}_k^\mathsf{T}\mathbf{B} = \mathbf{E}_1\hat{\mathbf{B}}_{-1}$. ∎

As such, block Lanczos-FA satisfies bounds similar to Corollary 6.4.

**Corollary 9.8.** For any unitarily invariant norm $|||\cdot|||$, the Lanczos-FA iterate satisfies

$$|||f(\mathbf{A})\mathbf{B} - \text{lan-FA}_k(f)||| \leq 2|||\mathbf{B}||| \min_{\deg(p)<k} \|f - p\|_{\mathcal{I}}.$$

Moreover, for each $i = 1, \ldots, m$, the Lanczos-FA iterate satisfies

$$\|f(\mathbf{A})\mathbf{b}_i - \text{lan-FA}_k(f)\mathbf{e}_i\| \leq 2\|\mathbf{b}_i\| \min_{\deg(p)<k} \|f - p\|_{\mathcal{I}}$$

*Proof.* For any polynomial $p(x)$ with $\deg(p) < k$, Theorem 9.7 asserts $p(\mathbf{A})\mathbf{B} = \text{lan-FA}_k(p)$. Moreover, since $\mathbf{Q}_k^\mathsf{T}\mathbf{B} = \mathbf{E}_1\hat{\mathbf{B}}_{-1}$, the block Lanczos-FA iterate can be written as $\text{lan-FA}_k(f) = \mathbf{Q}_k f(\mathbf{T}_k)\mathbf{Q}_k^\mathsf{T}\mathbf{B}$. Introduce the notation $e(x) := f(x) - p(x)$, so that $f(\mathbf{A})\mathbf{B} - \text{lan-FA}_k(f) = e(\mathbf{A})\mathbf{B} - \mathbf{Q}_k e(\mathbf{T}_k)\mathbf{Q}_k^\mathsf{T}\mathbf{B}$. Then, by the triangle inequality,

$$|||f(\mathbf{A})\mathbf{B} - \text{lan-FA}_k(f)||| \leq |||e(\mathbf{A})\mathbf{B}||| + |||\mathbf{Q}_k e(\mathbf{T}_k)\mathbf{Q}_k^\mathsf{T}\mathbf{B}||| \tag{9.13}$$

[1] In fact, a stronger result about “matrix polynomials” is true [FLS20, Theorem 2.7].

$$\leq \|e(\mathbf{A})\|\,|||\mathbf{B}||| + \|e(\mathbf{T}_k)\|\,|||\mathbf{B}|||, \tag{9.14}$$

where we have used that every unitarily invariant norm is sub-multiplicative with respect to the operator norm (i.e. $|||\mathbf{XY}||| \leq \|\mathbf{X}\|\,|||\mathbf{Y}|||$), together with $\|\mathbf{Q}_k\| = \|\mathbf{Q}_k^\mathsf{T}\| = 1$. Since $\Lambda \subset \mathcal{I}$, we have $\|e(\mathbf{A})\| = \|e\|_\Lambda \leq \|e\|_\mathcal{I}$; and since $\mathbf{T}_k = \mathbf{Q}_k^\mathsf{T}\mathbf{A}\mathbf{Q}_k$ with $\mathbf{Q}_k$ orthonormal, the eigenvalues of $\mathbf{T}_k$ are Rayleigh quotients of $\mathbf{A}$, so that $\Lambda(\mathbf{T}_k) \subset \mathcal{I}$ and $\|e(\mathbf{T}_k)\| = \|e\|_{\Lambda(\mathbf{T}_k)} \leq \|e\|_\mathcal{I}$. Hence

$$|||f(\mathbf{A})\mathbf{B} - \text{lan-FA}_k(f)||| \leq 2\|e\|_\mathcal{I}|||\mathbf{B}|||,$$

and minimizing over $p$ gives the first bound. For the second bound, the same argument applied to any vector $\mathbf{z}$ gives

$$\|f(\mathbf{A})\mathbf{B}\mathbf{z} - \text{lan-FA}_k(f)\mathbf{z}\| \leq \|e(\mathbf{A})\mathbf{B}\mathbf{z}\| + \|\mathbf{Q}_k e(\mathbf{T}_k)\mathbf{Q}_k^\mathsf{T}\mathbf{B}\mathbf{z}\| \tag{9.15}$$

$$\leq \big(\|e(\mathbf{A})\| + \|e(\mathbf{T}_k)\|\big)\|\mathbf{B}\mathbf{z}\| \leq 2\|e\|_\mathcal{I}\|\mathbf{B}\mathbf{z}\|. \tag{9.16}$$

Taking $\mathbf{z} = \mathbf{e}_i$, so that $\mathbf{B}\mathbf{z} = \mathbf{b}_i$, and minimizing over $p$ gives the result. ■

**9.2.3 Block Lanczos-QF** The block Lanczos-QF iterate can be defined analogously to the standard Lanczos-QF iterate (Definition 7.1).

**Definition 9.9.** The $k$-th Lanczos-QF approximation to $\mathbf{B}^\mathsf{T} f(\mathbf{A})\mathbf{B}$ is

$$\text{lan-QF}_k(f) = \text{lan-QF}_k(f; \mathbf{A}, \mathbf{B}) := \hat{\mathbf{B}}_{-1}^\mathsf{T}\mathbf{E}_1^\mathsf{T} f(\mathbf{T}_k)\mathbf{E}_1\hat{\mathbf{B}}_{-1}.$$

**Theorem 9.10.** Let $p(x)$ be a polynomial with $\deg(p) < 2k$. Then

$$\text{lan-QF}_k(p) = \mathbf{B}^\mathsf{T} p(\mathbf{A})\mathbf{B}.$$

*Proof.* It suffices to verify the result for $p(x) = x^\ell$ for $\ell = 0, 1, \ldots, 2k-1$. For $m, n < k$, using Theorem 9.7 and the fact that $\mathbf{Q}_k^\mathsf{T}\mathbf{Q}_k = \mathbf{I}$, we have that

$$\mathbf{B}^\mathsf{T}\mathbf{A}^{m+n}\mathbf{B} = (\mathbf{A}^m\mathbf{B})^\mathsf{T}\mathbf{A}^n\mathbf{B} \tag{9.17}$$

$$= (\mathbf{Q}_k(\mathbf{T}_k)^m\mathbf{E}_1\hat{\mathbf{B}}_{-1})^\mathsf{T}(\mathbf{Q}_k(\mathbf{T}_k)^n\mathbf{E}_1\hat{\mathbf{B}}_{-1}) \tag{9.18}$$

$$= \text{lan-QF}_k(x^{m+n}). \tag{9.19}$$

Next, using that $\mathbf{Q}_k^\mathsf{T}\mathbf{A}\mathbf{Q}_k = \mathbf{T}_k$,

$$\mathbf{B}^\mathsf{T}\mathbf{A}^{2k-1}\mathbf{B} = (\mathbf{A}^{k-1}\mathbf{B})^\mathsf{T}\mathbf{A}\mathbf{A}^{k-1}\mathbf{B} \tag{9.20}$$

$$= (\mathbf{Q}_k(\mathbf{T}_k)^{k-1}\mathbf{E}_1\hat{\mathbf{B}}_{-1})^\mathsf{T}\mathbf{A}(\mathbf{Q}_k(\mathbf{T}_k)^{k-1}\mathbf{E}_1\hat{\mathbf{B}}_{-1}) \tag{9.21}$$

$$= \text{lan-QF}_k(x^{2k-1}). \tag{9.22}$$

Hence, the result follows by linearity. ■

We can then obtain an analog of Corollary 7.4.

**Corollary 9.11.** For any unitarily invariant norm $|||\cdot|||$, the Lanczos-QF iterate satisfies

$$|||\mathbf{B}^\mathsf{T} f(\mathbf{A})\mathbf{B} - \text{lan-QF}_k(f)||| \leq 2\|\mathbf{B}\||||\mathbf{B}||| \min_{\deg(p)<2k} \|f-p\|_{\mathcal{I}}.$$

Moreover, for each $i, j = 1, \ldots, m$, the Lanczos-QF iterate satisfies

$$|\mathbf{b}_i^\mathsf{T} f(\mathbf{A})\mathbf{b}_j - \mathbf{e}_i^\mathsf{T} \text{lan-QF}_k(f)\mathbf{e}_j| \leq 2\|\mathbf{b}_i\|\|\mathbf{b}_j\| \min_{\deg(p)<2k} \|f-p\|_{\mathcal{I}}.$$

*Proof.* For any polynomial $p(x)$ with $\deg(p) < 2k$, Theorem 9.10 asserts $\mathbf{B}^\mathsf{T} p(\mathbf{A})\mathbf{B} = \text{lan-QF}_k(p)$. Moreover, since $\mathbf{Q}_k^\mathsf{T}\mathbf{B} = \mathbf{E}_1\hat{\mathbf{B}}_{-1}$, the block Lanczos-QF iterate can be written as $\text{lan-QF}_k(f) = \mathbf{B}^\mathsf{T}\mathbf{Q}_k f(\mathbf{T}_k)\mathbf{Q}_k^\mathsf{T}\mathbf{B}$. Introduce the notation $e(x) := f(x) - p(x)$, so that $\mathbf{B}^\mathsf{T} f(\mathbf{A})\mathbf{B} - \text{lan-QF}_k(f) = \mathbf{B}^\mathsf{T} e(\mathbf{A})\mathbf{B} - \mathbf{B}^\mathsf{T}\mathbf{Q}_k e(\mathbf{T}_k)\mathbf{Q}_k^\mathsf{T}\mathbf{B}$. Then, by the triangle inequality,

$$\begin{aligned} |||\mathbf{B}^\mathsf{T} f(\mathbf{A})\mathbf{B} - \text{lan-QF}_k(f)||| &\leq |||\mathbf{B}^\mathsf{T} e(\mathbf{A})\mathbf{B}||| + |||\mathbf{B}^\mathsf{T}\mathbf{Q}_k e(\mathbf{T}_k)\mathbf{Q}_k^\mathsf{T}\mathbf{B}||| && (9.23) \\ &\leq \|\mathbf{B}\|\|e(\mathbf{A})\||||\mathbf{B}||| + \|\mathbf{B}\|\|e(\mathbf{T}_k)\||||\mathbf{B}|||, && (9.24) \end{aligned}$$

where we have again used the sub-multiplicativity of $|||\cdot|||$ with respect to the spectral norm, together with $\|\mathbf{Q}_k\| = \|\mathbf{Q}_k^\mathsf{T}\| = 1$. As in the proof of Corollary 9.8, $\Lambda \subset \mathcal{I}$ gives $\|e(\mathbf{A})\| = \|e\|_\Lambda \leq \|e\|_{\mathcal{I}}$, and $\Lambda(\mathbf{T}_k) \subset \mathcal{I}$ gives $\|e(\mathbf{T}_k)\| = \|e\|_{\Lambda(\mathbf{T}_k)} \leq \|e\|_{\mathcal{I}}$. Hence

$$|||\mathbf{B}^\mathsf{T} f(\mathbf{A})\mathbf{B} - \text{lan-QF}_k(f)||| \leq 2\|\mathbf{B}\|\|e\|_{\mathcal{I}}|||\mathbf{B}|||,$$

and minimizing over $p$ gives the first bound. For the second bound, the same argument applied to any vectors $\mathbf{y}, \mathbf{z}$ gives

$$\begin{aligned} &|\mathbf{y}^\mathsf{T}\mathbf{B}^\mathsf{T} f(\mathbf{A})\mathbf{B}\mathbf{z} - \mathbf{y}^\mathsf{T}\,\text{lan-QF}_k(f)\mathbf{z}| && (9.25) \\ &\quad\leq |\mathbf{y}^\mathsf{T}\mathbf{B}^\mathsf{T} e(\mathbf{A})\mathbf{B}\mathbf{z}| + |\mathbf{y}^\mathsf{T}\mathbf{B}^\mathsf{T}\mathbf{Q}_k e(\mathbf{T}_k)\mathbf{Q}_k^\mathsf{T}\mathbf{B}\mathbf{z}| && (9.26) \\ &\quad\leq \big(\|e(\mathbf{A})\| + \|e(\mathbf{T}_k)\|\big)\|\mathbf{B}\mathbf{y}\|\|\mathbf{B}\mathbf{z}\| \leq 2\|e\|_{\mathcal{I}}\|\mathbf{B}\mathbf{y}\|\|\mathbf{B}\mathbf{z}\|. && (9.27) \end{aligned}$$

Taking $\mathbf{y} = \mathbf{e}_i$ and $\mathbf{z} = \mathbf{e}_j$, so that $\mathbf{B}\mathbf{y} = \mathbf{b}_i$ and $\mathbf{B}\mathbf{z} = \mathbf{b}_j$, and minimizing over $p$ gives the result. ■

**9.2.4 Are block methods better?** The bounds above are strikingly similar to the single-vector counterparts, and do not provide any real evidence that block-methods are preferable. For instance, while Theorem 9.4 shows that block-CG performs no worse than running CG independently, it does not show that it performs better.

One may argue that block methods make better use of parallelism, since matrix-vector products are replaced by a matrix-matrix product. However, if one is going to run multiple single-vector Krylov methods with the same matrix, then the computation can be arranged such that the matrix-vector products are done in a blocked manner. Thus, the computational profile is not a differentiating factor.

There is some work which suggests that block methods may be preferable in the presence of isolated eigenvalues [OLe80; Che+26b; Mat+26]. This makes intuitive sense: while each instance of a single-vector method will find its own copies of the isolated eigenvalues, a block method will find just one copy, which leaves additional degrees of freedom to do useful stuff. Unfortunately, any comprehensive characterization of the theoretical benefits of block Krylov methods for tasks involving matrix functions is a seemingly distant goal.

## 9.3 Some facts and observations

**9.3.1 The Krylov subspace of Krylov subspace is a Krylov subspace** The following observation about block Krylov subspaces has found use in the design and analysis of algorithms [CH23; MMM24; Che+26b].

**Theorem 9.12.** Suppose the columns of $\mathbf{Q}$ are each contained in $\mathcal{K}_{k+1}(\mathbf{A},\mathbf{\Omega})$. Then

$$\mathcal{K}_t(\mathbf{A},\mathbf{Q}) \subseteq \mathcal{K}_{k+t}(\mathbf{A},\mathbf{\Omega}),$$

with equality if the columns of $\mathbf{Q}$ span $\mathcal{K}_{k+1}(\mathbf{A},\mathbf{\Omega})$.

*Proof.* We have

$$\begin{aligned}\mathcal{K}_t(\mathbf{A},\mathbf{Q}) &= \operatorname{span}\{\mathbf{Q},\mathbf{A}\mathbf{Q},\ldots,\mathbf{A}^{t-1}\mathbf{Q}\}\\ &\subseteq \operatorname{span}\{\mathbf{\Omega},\mathbf{A}\mathbf{\Omega},\ldots,\mathbf{A}^{k}\mathbf{\Omega},\\ &\qquad\qquad \mathbf{A}\mathbf{\Omega},\mathbf{A}^2\mathbf{\Omega},\ldots,\mathbf{A}^{k+1}\mathbf{\Omega},\\ &\qquad\qquad\qquad \ddots\\ &\qquad\qquad\qquad \mathbf{A}^{t-1}\mathbf{\Omega},\mathbf{A}^{t}\mathbf{\Omega},\ldots,\mathbf{A}^{k+t-1}\mathbf{\Omega}\} = \mathcal{K}_{k+t}(\mathbf{A},\mathbf{\Omega}).\end{aligned} \tag{9.28}$$

The second result follows similarly, with the subset inclusion replaced by equality. ■

We note that a good block-Lanczos implementation would multiply $\mathbf{A}$ with $\mathbf{Q}$, realize that $\mathbf{AQ}$ has a small rank, and automatically do deflation. Thus, while it would unnecessarily do some products with $\mathbf{A}$, the overall cost would be similar to an algorithm designed with Theorem 9.12 in mind. Even so, Theorem 9.12 is of conceptual value.

**9.3.2 Random block Krylov subspaces** In many settings where block KSMs are used, $\mathbf{B}$ is a random Gaussian matrix. In this case, we can provide a sufficient conditions for the Krylov subspace $\mathcal{K}_k(\mathbf{A}, \mathbf{B})$ to be of rank $km$. Quantitative bounds on the conditioning of Krylov matrices are derived in [Che+26a].

**Theorem 9.13.** Suppose each distinct eigenvalue of $\mathbf{A}$ has multiplicity at most $m$, and let $\mathbf{B}$ be a $d \times m$ random Gaussian matrix. Then, with probability one,

$$\dim(\mathcal{K}_k(\mathbf{A}, \mathbf{B})) = \min\{km, d\}.$$

Inspired by work analyzing block Krylov subspaces over finite fields [Kal95] we make use of the Schwartz–Zippel Lemma.

**Lemma 9.14** (Schwartz–Zippel)**.** Suppose $p : \mathbb{R}^q \to \mathbb{R}$ is a a nonzero polynomial of finite total degree and $\mathbf{x}$ is a random Gaussian vector. Then, with probability one, $p(\mathbf{x}) \neq 0$.

*Proof.* For convenience, suppose $km \leq d$. We will subsequently describe how to relax this condition.

We begin by relating the block Krylov subspace to a multivariate polynomial of the entries of the starting block as follows. Define the matrix $\mathbf{K}_k(\mathbf{B})$ by

$$\mathbf{K}_k(\mathbf{B}) = \begin{bmatrix} \mathbf{B} & \mathbf{AB} & \cdots & \mathbf{A}^{k-1}\mathbf{B} \end{bmatrix} \in \mathbb{R}^{d \times km}, \tag{9.29}$$

and define the polynomial $p : \mathbb{R}^{d \times m} \to \mathbb{R}$ by

$$p(\mathbf{B}) = \det(\mathbf{K}_k(\mathbf{B})^\mathsf{T}\mathbf{K}_k(\mathbf{B})). \tag{9.30}$$

Recall that $km \leq d$ by assumption. Therefore, the matrix $\mathbf{K}_k(\mathbf{B})$ is of full-rank $km$ exactly when $p(\mathbf{B}) \neq 0$.

If we can show that $p : \mathbb{R}^{d \times m} \to \mathbb{R}$ is not the zero polynomial, then the result follows immediately from the Schwartz–Zippel Lemma (Lemma 9.14); i.e we must exhibit a matrix $\hat{\mathbf{B}}$ such that $p(\hat{\mathbf{B}}) \neq 0$.

Since each distinct eigenvalue has multiplicity at most $m$, we can partition the $d$ eigenvalues of $\mathbf{A}$ into $m$ groups such that (i) each group has at least $k$ eigenvalues, and (ii) no eigenvalues within in a group are repeated. Work in a basis such that

$$\mathbf{A} = \begin{bmatrix} \mathbf{A}_1 & & & \\ & \mathbf{A}_2 & & \\ & & \ddots & \\ & & & \mathbf{A}_m \end{bmatrix}, \tag{9.31}$$

where each $\mathbf{A}_i$ is diagonal and corresponds to one of the aforementioned groups of eigenvalues. Denote the all ones vector (of appropriate size) by $\mathbf{1}$ and define $\hat{\mathbf{B}}$ by

$$\hat{\mathbf{B}} = \begin{bmatrix} \mathbf{1} & & & \\ & \mathbf{1} & & \\ & & \ddots & \\ & & & \mathbf{1} \end{bmatrix}. \tag{9.32}$$

Observe then that, for some permutation matrix $\mathbf{P}$

$$\mathbf{K}_k(\hat{\mathbf{B}}) = \begin{bmatrix} \hat{\mathbf{K}}_1 & & & \\ & \hat{\mathbf{K}}_2 & & \\ & & \ddots & \\ & & & \hat{\mathbf{K}}_m \end{bmatrix} \mathbf{P}, \qquad \hat{\mathbf{K}}_i = \begin{bmatrix} \mathbf{1} & \mathbf{A}_i\mathbf{1} & \cdots & \mathbf{A}_i^{k-1}\mathbf{1} \end{bmatrix}. \tag{9.33}$$

For each $i$, $\mathbf{1}$ has nonzero projection onto each of the $\mathbf{A}_i$ and hence by Lemma 2.2, $\operatorname{rank}(\hat{\mathbf{K}}_i) = k$. So then

$$\operatorname{rank}(\mathbf{K}_k(\hat{\mathbf{B}})) = \operatorname{rank}(\hat{\mathbf{K}}_1) + \cdots + \operatorname{rank}(\hat{\mathbf{K}}_m) = km. \tag{9.34}$$

Therefore, $\mathbf{K}_k(\hat{\mathbf{B}})$ is full rank and $p(\hat{\mathbf{B}}) \neq 0$ as desired.

Now, assume $km > d$. The proof is similar, but now we define $p(\mathbf{B}) = \det(\mathbf{K}_k(\mathbf{B})\mathbf{K}_k(\mathbf{B})^\mathsf{T})$ which has nonzero determinant when $\mathbf{K}_k(\mathbf{B})$ is of full rank $d$. Since $d < km$, we can group the eigenvalues so that each group has at most $k$ eigenvalues. Then the rank of $\hat{\mathbf{K}}_i$ will be exactly the number of eigenvalues in that group, and the sum is $d$. ■

## TL;DR

Block methods can be used for many tasks involving matrix functions, and satisfy many of the same bounds as single-vector methods. While theory for their behavior (in exact or finite precision arithmetic) is less developed than their single-vector counterparts, they are often used in practice with good results.

# 10 References

[ADN20] W. Ahrens, J. Demmel, and H. D. Nguyen. "Algorithms for Efficient Reproducible Floating Point Summation". In: *ACM Transactions on Mathematical Software* 46.3 (July 2020), pp. 1–49. ISSN: 1557-7295. DOI: 10.1145/3389360 (cited on page 50).

[Aic+03] M. Aichhorn, M. Daghofer, H. G. Evertz, and W. von der Linden. "Low-temperature Lanczos method for strongly correlated systems". In: *Physical Review B* 67.16 (Apr. 2003). DOI: 10.1103/physrevb.67.161103 (cited on page 21).

[Alb+75] R. Alben, M. Blume, H. Krakauer, and L. Schwartz. "Exact results for a three-dimensional alloy with site diagonal disorder: comparison with the coherent potential approximation". In: *Physical Review B* 12.10 (Nov. 1975), pp. 4090–4094. DOI: 10.1103/physrevb.12.4090 (cited on page 54).

[Ams+24] N. Amsel, T. Chen, A. Greenbaum, C. Musco, and C. Musco. "Near-Optimal Approximation of Matrix Functions by the Lanczos Method". In: *Advances in Neural Information Processing Systems*. 2024. arXiv: 2303.03358 [math.NA] (cited on page 40).

[Arn51] W. E. Arnoldi. "The principle of minimized iterations in the solution of the matrix eigenvalue problem". In: *Quarterly of Applied Mathematics* 9.1 (1951), pp. 17–29. ISSN: 1552-4485. DOI: 10.1090/qam/42792 (cited on page 6).

[AT11] H. Avron and S. Toledo. "Randomized algorithms for estimating the trace of an implicit symmetric positive semi-definite matrix". In: *Journal of the ACM* 58.2 (Apr. 2011), pp. 1–34. DOI: 10.1145/1944345.1944349 (cited on page 54).

[Avr10] H. Avron. "Counting Triangles in Large Graphs using Randomized Matrix Trace Estimation". In: *Proceedings of KDD-LDMTA*. 2010. eprint: http://www.math.tau.ac.il/~haimav/Avron-KDD-LDMTA10.pdf (cited on page 51).

[BB20] M. Benzi and P. Boito. "Matrix functions in network analysis". In: *GAMM-Mitteilungen* 43.3 (2020), e202000012. DOI: 10.1002/gamm.202000012. eprint: https://onlinelibrary.wiley.com/doi/pdf/10.1002/gamm.202000012 (cited on page 51).

[BB98] B. Beckermann and E. Bourreau. "How to choose modified moments?" In: *Journal of Computational and Applied Mathematics* 98.1 (Oct. 1998), pp. 81–98. ISSN: 0377-0427. DOI: 10.1016/s0377-0427(98)00116-2 (cited on page 11).

[BFG96] Z. Bai, G. Fahey, and G. Golub. "Some large-scale matrix computation problems". In: *Journal of Computational and Applied Mathematics* 74.1-2 (Nov. 1996), pp. 71–89. DOI: 10.1016/0377-0427(96)00018-0 (cited on page 55).

[BG96] Z. Bai and G. Golub. "Bounds for the trace of the inverse and the determinant of symmetric positive definite matrices". In: *Annals of Numerical Mathematic* 4 (Apr. 1996), pp. 29–38 (cited on page 55).

[BKM22] V. Braverman, A. Krishnan, and C. Musco. "Sublinear time spectral density estimation". In: *Proceedings of the 54th Annual ACM SIGACT Symposium on Theory of Computing*. arXiv cs.DS 2104.03461. ACM, June 2022. DOI: 10.1145/3519935.3520009. arXiv: 2104.03461 [cs.DS] (cited on pages 58, 60, 64).

[Bor00] A. Boriçi. "Fast Methods for Computing the Neuberger Operator". In: Springer Berlin Heidelberg, 2000, pp. 40–47. DOI: 10.1007/978-3-642-58333-9_4 (cited on page 49).

[BP99] R. P. Barry and R. K. Pace. "Monte Carlo estimates of the log determinant of large sparse matrices". In: *Linear Algebra and its Applications* 289.1-3 (Mar. 1999), pp. 41–54. DOI: 10.1016/s0024-3795(97)10009-x (cited on page 51).

[Bro91] P. N. Brown. "A Theoretical Comparison of the Arnoldi and GMRES Algorithms". In: *SIAM Journal on Scientific and Statistical Computing* 12.1 (Jan. 1991), pp. 58–78. ISSN: 2168-3417. DOI: 10.1137/0912003 (cited on page 37).

[BS22] K. Bergermann and M. Stoll. "Fast computation of matrix function-based centrality measures for layer-coupled multiplex networks". In: *Physical Review E* 105.3 (Mar. 2022). ISSN: 2470-0053. DOI: 10.1103/physreve.105.034305 (cited on page 51).

[Cab97] C. Cabos. "Evaluation of matrix functions with the block Lanczos algorithm". In: *Computers & Mathematics with Applications* 33.1–2 (Jan. 1997), pp. 45–57. ISSN: 0898-1221. DOI: 10.1016/s0898-1221(96)00218-0 (cited on page 66).

[CG96] J. Cullum and A. Greenbaum. "Relations between Galerkin and Norm-Minimizing Iterative Methods for Solving Linear Systems". In: *SIAM Journal on Matrix Analysis and Applications* 17.2 (1996), pp. 223–247. DOI: 10.1137/S0895479893246765. eprint: https://doi.org/10.1137/S0895479893246765 (cited on page 37).

[CH23] T. Chen and E. Hallman. "Krylov-Aware Stochastic Trace Estimation". In: *SIAM Journal on Matrix Analysis and Applications* 44.3 (Aug. 2023), pp. 1218–1244. DOI: 10.1137/22m1494257. arXiv: 2205.01736 [math.NA] (cited on pages 57, 71).

[Che+22] T. Chen, A. Greenbaum, C. Musco, and C. Musco. "Error Bounds for Lanczos-Based Matrix Function Approximation". In: *SIAM Journal on Matrix Analysis and Applications* 43.2 (May 2022), pp. 787–811. DOI: 10.1137/21m1427784. arXiv: 2106.09806 [math.NA] (cited on pages 40, 44).

[Che+26a] T. Chen, E. N. Epperly, R. A. Meyer, C. Musco, and A. Rao. "Does block size matter in randomized block Krylov low-rank approximation?" In: *Proceedings of the 2026 Annual ACM-SIAM Symposium on Discrete Algorithms (SODA)*. Society for Industrial and Applied Mathematics, Jan. 2026, pp. 1026–1046. ISBN: 9781611978971. DOI: 10.1137/1.9781611978971.42 (cited on page 72).

[Che+26b] T. Chen, C. Huber, E. Lin, and H. Zaid. "Preconditioning without a preconditioner using randomized block Krylov subspace methods". In: *ETNA - Electronic Transactions on Numerical Analysis* 65 (2026), pp. 63–92. ISSN: 1068-9613. DOI: 10.1553/etna_vol65s63. arXiv: 2501.18717 [math.NA] (cited on page 71).

[Che05] K. Chen. *Matrix preconditioning techniques and applications*. Cambridge monographs on applied and computational mathematics. Cambridge ; New York: Cambridge University Press, 2005. ISBN: 9780521838283 (cited on page 38).

[Che23] T. Chen. "A spectrum adaptive kernel polynomial method". In: *The Journal of Chemical Physics* 159.11 (Sept. 2023), p. 114101. DOI: 10.1063/5.0166678. arXiv: 2308.15683 [physics.comp-ph] (cited on page 64).

[CK21] A. Cortinovis and D. Kressner. "On Randomized Trace Estimates for Indefinite Matrices with an Application to Determinants". In: *Foundations of Computational Mathematics* (July 2021). DOI: 10.1007/s10208-021-09525-9 (cited on page 55).

[Cle55] C. W. Clenshaw. "A note on the summation of Chebyshev series". In: *Mathematics of Computation* 9.51 (1955), pp. 118–120. ISSN: 1088-6842. DOI: 10.1090/s0025-5718-1955-0071856-0 (cited on page 46).

[CLS24] E. Carson, J. Liesen, and Z. Strakoš. "Towards understanding CG and GMRES through examples". In: *Linear Algebra and its Applications* 692 (July 2024), pp. 241–291. ISSN: 0024-3795. DOI: 10.1016/j.laa.2024.04.003 (cited on pages 1, 24).

[CM24] T. Chen and G. Meurant. "Near-optimal convergence of the full orthogonalization method". In: *Electronic Transactions on Numerical Analysis* 60 (2024), pp. 421–427. DOI: 10.1553/etna_vol60s421. arXiv: 2403.07259 [math.NA] (cited on page 37).

[CP26] T. Chen and D. Persson. *Optimal near-optimality bounds for the Lanczos method for matrix functions*. 2026. arXiv: 2608.07160 [math.NA] (cited on page 44).

[CT24] T. Chen and T. Trogdon. "Stability of the Lanczos algorithm on matrices with regular spectral distributions". In: *Linear Algebra and its Applications* 682 (Feb. 2024), pp. 191–237. ISSN: 0024-3795. DOI: 10.1016/j.laa.2023.11.006. arXiv: 2302.14842 [math.NA] (cited on pages 19, 53).

[CTU21] T. Chen, T. Trogdon, and S. Ubaru. "Analysis of stochastic Lanczos quadrature for spectrum approximation". In: *Proceedings of the 38th International Conference on Machine Learning*. Vol. 139. Proceedings of Machine Learning Research. PMLR, July 2021, pp. 1728–1739. arXiv: 2105.06595 [cs.DS] (cited on page 60).

[CTU25] T. Chen, T. Trogdon, and S. Ubaru. "Randomized Matrix-Free Quadrature: Unified and Uniform Bounds for Stochastic Lanczos Quadrature and the Kernel Polynomial Method". In: *SIAM Journal on Scientific Computing* 47.3 (May 2025), A1733–A1757. ISSN: 1095-7197. DOI: 10.1137/23m1600414. arXiv: 2204.01941 [math.NA] (cited on page 64).

[CW12] J. K. Cullum and R. A. Willoughby. *Lanczos algorithms for large symmetric eigenvalue computations vol. II programs*. en. Progress in Scientific Computing. Secaucus, NJ: Birkhauser Boston, June 2012 (cited on page 1).

[DBB19] K. Dong, A. R. Benson, and D. Bindel. "Network Density of States". In: *Proceedings of the 25th ACM SIGKDD International Conference on Knowledge Discovery & Data Mining*. ACM, July 2019. DOI: 10.1145/3292500.3330891 (cited on page 51).

[DGK98] V. Druskin, A. Greenbaum, and L. Knizhnerman. "Using Nonorthogonal Lanczos Vectors in the Computation of Matrix Functions". In: *SIAM Journal on Scientific Computing* 19.1 (1998), pp. 38–54. DOI: 10.1137/S1064827596303661. eprint: https://doi.org/10.1137/S1064827596303661 (cited on page 40).

[DK88] V. Druskin and L. Knizhnerman. "Spectral Differential-Difference Method for Numeric Solution of Three-Dimensional Nonstationary Problems of Electric Prospecting". In: *Physics of the Solid Earth* 24 (Jan. 1988), pp. 641–648 (cited on page 40).

[DK89] V. Druskin and L. Knizhnerman. "Two polynomial methods of calculating functions of symmetric matrices". In: *USSR Computational Mathematics and Mathematical Physics* 29.6 (Jan. 1989), pp. 112–121. DOI: 10.1016/s0041-5553(89)80020-5 (cited on pages 40, 42).

[DK91] V. L. Druskin and L. A. Knizhnerman. "Error Bounds in the Simple Lanczos Procedure for Computing Functions of Symmetric Matrices and Eigenvalues". In: *Comput. Math. Math. Phys.* 31.7 (July 1991), pp. 20–30. ISSN: 0965-5425 (cited on pages 33, 40, 45).

[DK95] V. Druskin and L. Knizhnerman. "Krylov subspace approximation of eigenpairs and matrix functions in exact and computer arithmetic". In: *Numerical Linear Algebra with Applications* 2.3 (May 1995), pp. 205–217. DOI: 10.1002/nla.1680020303 (cited on page 40).

[Dru08] V. Druskin. "On monotonicity of the Lanczos approximation to the matrix exponential". In: *Linear Algebra and its Applications* 429.7 (Oct. 2008), pp. 1679–1683. ISSN: 0024-3795. DOI: 10.1016/j.laa.2008.04.046 (cited on page 40).

[Epp24] E. N. Epperly. *Note to Self: Hanson-Wright and Trace Estimation with Random Vectors on the Sphere*. 2024. eprint: https://www.ethanepperly.com/index.php/2024/11/04/note-to-self-hanson-wright-and-trace-estimation-with-random-vectors-on-the-sphere/ (cited on page 55).

[Esh+02] J. van den Eshof, A. Frommer, T. Lippert, K. Schilling, and H. van der Vorst. "Numerical methods for the QCDd overlap operator. I. Sign-function and error bounds". In: *Computer Physics Communications* 146.2 (2002), pp. 203–224. ISSN: 0010-4655. DOI: 10.1016/S0010-4655(02)00455-1 (cited on page 40).

[Est00] E. Estrada. "Characterization of 3D molecular structure". In: *Chemical Physics Letters* 319.5-6 (Mar. 2000), pp. 713–718. DOI: 10.1016/s0009-2614(00)00158-5 (cited on page 51).

[ETW24] E. N. Epperly, J. A. Tropp, and R. J. Webber. "XTrace: Making the Most of Every Sample in Stochastic Trace Estimation". In: *SIAM Journal on Matrix Analysis and Applications* 45.1 (Jan. 2024), pp. 1–23. ISSN: 1095-7162. DOI: 10.1137/23m1548323 (cited on page 56).

[Fan+19] L. Fan, D. I. Shuman, S. Ubaru, and Y. Saad. "Spectrum-adapted Polynomial Approximation for Matrix Functions". In: *ICASSP 2019 - 2019 IEEE International Conference on Acoustics, Speech and Signal Processing (ICASSP)*. IEEE, May 2019. DOI: 10.1109/icassp.2019.8683179 (cited on page 58).

[FGS14] A. Frommer, S. Güttel, and M. Schweitzer. "Convergence of Restarted Krylov Subspace Methods for Stieltjes Functions of Matrices". In: *SIAM Journal on Matrix Analysis and Applications* 35.4 (Jan. 2014), pp. 1602–1624. DOI: 10.1137/140973463 (cited on page 44).

[Fis11] B. Fischer. *Polynomial Based Iteration Methods for Symmetric Linear Systems*. Society for Industrial and Applied Mathematics, Jan. 2011. DOI: 10.1137/1.9781611971927 (cited on page 36).

[FLS20] A. Frommer, K. Lund, and D. B. Szyld. "Block Krylov Subspace Methods for Functions of Matrices II: Modified Block FOM". In: *SIAM Journal on Matrix Analysis and Applications* 41.2 (Jan. 2020), pp. 804–837. ISSN: 1095-7162. DOI: 10.1137/19m1255847 (cited on page 68).

[Fro+13] A. Frommer, K. Kahl, T. Lippert, and H. Rittich. "2-Norm Error Bounds and Estimates for Lanczos Approximations to Linear Systems and Rational Matrix Functions". In: *SIAM Journal on Matrix Analysis and Applications* 34.3 (2013), pp. 1046–1065. DOI: 10.1137/110859749 (cited on page 44).

[Fro+17] A. Frommer, K. Lund, M. Schweitzer, and D. B. Szyld. "The Radau–Lanczos Method for Matrix Functions". In: *SIAM Journal on Matrix Analysis and Applications* 38.3 (2017), pp. 710–732. DOI: 10.1137/16M1072565 (cited on page 66).

[FS08a] A. Frommer and V. Simoncini. "Matrix Functions". In: *Mathematics in Industry*. Springer Berlin Heidelberg, 2008, pp. 275–303. DOI: 10.1007/978-3-540-78841-6_13 (cited on page 49).

[FS08b] A. Frommer and V. Simoncini. "Stopping Criteria for Rational Matrix Functions of Hermitian and Symmetric Matrices". In: *SIAM Journal on Scientific Computing* 30.3 (Jan. 2008), pp. 1387–1412. DOI: 10.1137/070684598 (cited on page 44).

[FS09] A. Frommer and V. Simoncini. "Error Bounds for Lanczos Approximations of Rational Functions of Matrices". In: *Numerical Validation in Current Hardware Architectures*. Berlin, Heidelberg: Springer Berlin Heidelberg, 2009, pp. 203–216. ISBN: 978-3-642-01591-5 (cited on page 44).

[FS15] A. Frommer and M. Schweitzer. "Error bounds and estimates for Krylov subspace approximations of Stieltjes matrix functions". In: *BIT Numerical Mathematics* 56.3 (Dec. 2015), pp. 865–892. DOI: 10.1007/s10543-015-0596-3 (cited on page 44).

[Gau04] W. Gautschi. *Orthogonal polynomials: computation and approximation*. Numerical mathematics and scientific computation. Oxford University Press, 2004. ISBN: 9780198506720 (cited on pages 11, 13).

[Gau06] W. Gautschi. "Orthogonal Polynomials, Quadrature, and Approximation: Computational Methods and Software (in Matlab)". In: *Lecture Notes in Mathematics*. Springer Berlin Heidelberg, 2006, pp. 1–77. DOI: 10.1007/978-3-540-36716-1_1 (cited on page 10).

[Gir87] D. Girard. *Un algorithme simple et rapide pour la validation croisée généralisée sur des problèmes de grande taille*. 1987 (cited on pages 54, 56).

[GKX19] B. Ghorbani, S. Krishnan, and Y. Xiao. "An Investigation into Neural Net Optimization via Hessian Eigenvalue Density". In: *Proceedings of the 36th International Conference on Machine Learning*. Ed. by K. Chaudhuri and R. Salakhutdinov. Vol. 97. Proceedings of Machine Learning Research. PMLR, June 2019, pp. 2232–2241. arXiv: 1901.10159 [cs.LG] (cited on page 58).

[GM09] G. H. Golub and G. Meurant. *Matrices, moments and quadrature with applications*. Vol. 30. Princeton series in applied mathematics. Princeton University Press, 2009. ISBN: 9780691143415 (cited on pages 9, 10, 13).

[GM93] G. H. Golub and G. Meurant. "Matrices, moments and quadrature". In: *Numerical analysis 1993*. CRC Press, 1993, pp. 105–156 (cited on page 10).

[GM97] G. H. Golub and G. Meurant. "Matrices, moments and quadrature II; How to compute the norm of the error in iterative methods". In: *BIT Numerical Mathematics* 37.3 (Sept. 1997), pp. 687–705. DOI: 10.1007/bf02510247 (cited on page 10).

[Gol+10] S. Goldstein, J. L. Lebowitz, C. Mastrodonato, R. Tumulka, and N. Zanghì. "Normal typicality and von Neumann's quantum ergodic theorem". In: *Proceedings of the Royal Society A: Mathematical, Physical and Engineering Sciences* 466.2123 (May 2010), pp. 3203–3224. DOI: 10.1098/rspa.2009.0635 (cited on page 54).

[Gre89] A. Greenbaum. "Behavior of slightly perturbed Lanczos and conjugate-gradient recurrences". In: *Linear Algebra and its Applications* 113 (1989), pp. 7–63. ISSN: 0024-3795. DOI: 10.1016/0024-3795(89)90285-1 (cited on page 24).

[Gre97] A. Greenbaum. *Iterative Methods for Solving Linear Systems*. Philadelphia, PA, USA: Society for Industrial and Applied Mathematics, 1997. ISBN: 0-89871-396-X (cited on pages 1, 36).

[GS21] S. Güttel and M. Schweitzer. "A Comparison of Limited-memory Krylov Methods for Stieltjes Functions of Hermitian Matrices". In: *SIAM Journal on Matrix Analysis and Applications* 42.1 (Jan. 2021), pp. 83–107. DOI: 10.1137/20m1351072 (cited on page 50).

[GS92a] E. Gallopoulos and Y. Saad. "Efficient Solution of Parabolic Equations by Krylov Approximation Methods". In: *SIAM Journal on Scientific and Statistical Computing* 13.5 (Sept. 1992), pp. 1236–1264. DOI: 10.1137/0913071 (cited on page 40).

[GS92b] A. Greenbaum and Z. Strakos. "Predicting the Behavior of Finite Precision Lanczos and Conjugate Gradient Computations". In: *SIAM Journal on Matrix Analysis and Applications* 13.1 (Jan. 1992), pp. 121–137. ISSN: 1095-7162. DOI: 10.1137/0613011 (cited on page 24).

[GSO17] A. S. Gambhir, A. Stathopoulos, and K. Orginos. "Deflation as a Method of Variance Reduction for Estimating the Trace of a Matrix Inverse". In: *SIAM Journal on Scientific Computing* 39.2 (Jan. 2017), A532–A558. DOI: 10.1137/16m1066361 (cited on page 56).

[GU77] G. Golub and R. Underwood. "The Block Lanczos Method for Computing Eigenvalues". In: *Mathematical Software*. Elsevier, 1977, pp. 361–377. ISBN: 9780125872607. DOI: 10.1016/b978-0-12-587260-7.50018-2 (cited on page 65).

[GWG19] D. Granziol, X. Wan, and T. Garipov. *Deep Curvature Suite*. 2019. arXiv: 1912.09656 [stat.ML] (cited on pages 21, 58).

[Hig08] N. J. Higham. *Functions of Matrices*. Society for Industrial and Applied Mathematics, 2008. DOI: 10.1137/1.9780898717778 (cited on page 40).

[HL97] M. Hochbruck and C. Lubich. "On Krylov Subspace Approximations to the Matrix Exponential Operator". In: *SIAM Journal on Numerical Analysis* 34.5 (Oct. 1997), pp. 1911–1925. DOI: 10.1137/s0036142995280572 (cited on page 40).

[HS52] M. R. Hestenes and E. Stiefel. *Methods of conjugate gradients for solving linear systems*. Vol. 49. NBS Washington, DC, 1952 (cited on pages 28, 33).

[Hut89] M. Hutchinson. "A Stochastic Estimator of the Trace of the Influence Matrix for Laplacian Smoothing Splines". In: *Communications in Statistics - Simulation and Computation* 18.3 (Jan. 1989), pp. 1059–1076. DOI: 10.1080/03610918908812806 (cited on page 54).

[Iak+15] R. Iakymchuk, C. Collange, D. Defour, and S. Graillat. "ExBLAS: Reproducible and Accurate BLAS Library". In: *NRE: Numerical Reproducibility at Exascale*. Numerical Reproducibility at Exascale (NRE2015) workshop held as part of the Supercomputing Conference (SC15). Austin, TX, United States, Nov. 2015 (cited on page 50).

[ITS09] M. D. Ilic, I. W. Turner, and D. P. Simpson. "A restarted Lanczos approximation to functions of a symmetric matrix". In: *IMA Journal of Numerical Analysis* 30.4 (June 2009), pp. 1044–1061. DOI: 10.1093/imanum/drp003 (cited on page 44).

[Jac12] D. Jackson. “On Approximation by Trigonometric Sums and Polynomials”. In: *Transactions of the American Mathematical Society* 13.4 (1912), pp. 491–515. ISSN: 00029947 (cited on page 60).

[Jin+21] F. Jin, D. Willsch, M. Willsch, H. Lagemann, K. Michielsen, and H. De Raedt. “Random State Technology”. In: *Journal of the Physical Society of Japan* 90.1 (Jan. 2021), p. 012001. ISSN: 1347-4073. DOI: 10.7566/jpsj.90.012001 (cited on pages 51, 54, 58).

[Jin+23] Y. Jin, C. Musco, A. Sidford, and A. V. Singh. “Moments, Random Walks, and Limits for Spectrum Approximation”. In: *Proceedings of Thirty Sixth Conference on Learning Theory*. Ed. by G. Neu and L. Rosasco. Vol. 195. Proceedings of Machine Learning Research. PMLR, July 2023, pp. 5373–5394 (cited on page 61).

[Joz94] R. Jozsa. “Fidelity for Mixed Quantum States”. In: *Journal of Modern Optics* 41.12 (Dec. 1994), pp. 2315–2323. DOI: 10.1080/09500349414552171 (cited on page 51).

[JP94] J. Jaklič and P. Prelovšek. “Lanczos method for the calculation of finite-temperature quantities in correlated systems”. In: *Physical Review B* 49.7 (Feb. 1994), pp. 5065–5068. DOI: 10.1103/physrevb.49.5065 (cited on page 21).

[JS19] Y. Jin and A. Sidford. “Principal Component Projection and Regression in Nearly Linear Time through Asymmetric SVRG”. In: *Advances in Neural Information Processing Systems 32*. Ed. by H. Wallach, H. Larochelle, A. Beygelzimer, F. d Alché-Buc, E. Fox, and R. Garnett. Curran Associates, Inc., 2019, pp. 3868–3878. arXiv: 1910.06517 [cs.DS] (cited on page 40).

[Kal95] E. Kaltofen. “Analysis of Coppersmith’s block Wiedemann algorithm for the parallel solution of sparse linear systems”. In: *Mathematics of Computation* 64.210 (1995), pp. 777–806. ISSN: 1088-6842. DOI: 10.1090/s0025-5718-1995-1270621-1 (cited on page 72).

[Kni96] L. A. Knizhnerman. “The Simple Lanczos Procedure: Estimates of the Error of the Gauss Quadrature Formula and Their Applications”. In: *Comput. Math. Math. Phys.* 36.11 (Jan. 1996), pp. 1481–1492. ISSN: 0965-5425 (cited on pages 26, 53).

[KS72] S. Karlin and L. S. Shapley. *Geometry of moment spaces*. 3. printing. Memoirs of the American Mathematical Society. American Math. Soc, 1972. ISBN: 9780821812129 (cited on page 12).

[KV17] W. Kong and G. Valiant. “Spectrum estimation from samples”. In: *The Annals of Statistics* 45.5 (Oct. 2017), pp. 2218–2247. DOI: 10.1214/16-aos1525 (cited on pages 58, 61).

[Lan50] C. Lanczos. “An iteration method for the solution of the eigenvalue problem of linear differential and integral operators”. In: *Journal of research of the National Bureau of Standards* 45 (1950), pp. 255–282 (cited on pages 1, 7).

[Li+19] R. Li, Y. Xi, L. Erlandson, and Y. Saad. “The Eigenvalues Slicing Library (EVSL): Algorithms, Implementation, and Software”. In: *SIAM Journal on Scientific Computing* 41.4 (Jan. 2019), pp. C393–C415. DOI: 10.1137/18m1170935 (cited on page 58).

[Li22] H. Li. personal communication. 2022 (cited on page 50).

[LS13] J. Liesen and Z. Strakoš. *Krylov subspace methods: principles and analysis*. 1st ed. Numerical mathematics and scientific computation. Oxford University Press, 2013. ISBN: 9780199655410 (cited on pages 1, 33).

[Mat+26] F. Matti, M. S. Andersen, T. Chen, and D. Kressner. *Kernel-based linear system identification using augmented Krylov subspaces*. 2026. arXiv: 2605.08362 [math.NA] (cited on page 71).

[MD20] G. Meurant and J. Duintjer Tebbens. *Krylov Methods for Nonsymmetric Linear Systems: From Theory to Computations*. Springer International Publishing, 2020. ISBN: 9783030552510. DOI: 10.1007/978-3-030-55251-0 (cited on page 36).

[Meu06] G. Meurant. *The Lanczos and Conjugate Gradient Algorithms*. Society for Industrial and Applied Mathematics, 2006. DOI: 10.1137/1.9780898718140. eprint: https://epubs.siam.org/doi/pdf/10.1137/1.9780898718140 (cited on pages 1, 22, 23, 32).

[Mey+21] R. A. Meyer, C. Musco, C. Musco, and D. P. Woodruff. “Hutch++: Optimal Stochastic Trace Estimation”. In: *Symposium on Simplicity in Algorithms (SOSA)*. Society for Industrial and Applied Mathematics, Jan. 2021, pp. 142–155. DOI: 10.1137/1.9781611976496.16 (cited on pages 54–56).

[MH03] J. C. Mason and D. C. Handscomb. *Chebyshev polynomials*. Boca Raton, Fla: Chapman & Hall/CRC, 2003. ISBN: 9780849303555 (cited on page 19).

[MMM24] R. Meyer, C. Musco, and C. Musco. “On the Unreasonable Effectiveness of Single Vector Krylov Methods for Low-Rank Approximation”. In: *Proceedings of the 2024 Annual ACM-SIAM Symposium on Discrete Algorithms (SODA)*. Society for Industrial and Applied Mathematics, Jan. 2024, pp. 811–845. ISBN: 9781611977912. DOI: 10.1137/1.9781611977912.32 (cited on page 71).

[MMS18] C. Musco, C. Musco, and A. Sidford. “Stability of the Lanczos Method for Matrix Function Approximation”. In: Society for Industrial and Applied Mathematics, Jan. 2018, pp. 1605–1624. DOI: 10.1137/1.9781611975031.105 (cited on pages 40, 45).

[MS06] G. Meurant and Z. Strakoš. “The Lanczos and conjugate gradient algorithms in finite precision arithmetic”. In: *Acta Numerica* 15 (May 2006), pp. 471–542. DOI: 10.1017/s096249290626001x (cited on pages 1, 22).

[MT20] P.-G. Martinsson and J. A. Tropp. “Randomized numerical linear algebra: Foundations and algorithms”. In: *Acta Numerica* 29 (May 2020), pp. 403–572. DOI: 10.1017/s0962492920000021 (cited on page 55).

[MT24] G. A. Meurant and P. Tichý. *Error norm estimation in the conjugate gradient algorithm*. SIAM spotlights. Philadelphia: Society for Industrial and Applied Mathematics, 2024. ISBN: 9781611977868 (cited on page 31).

[Mus+25] C. Musco, C. Musco, L. Rosenblatt, and A. V. Singh. “Sharper Bounds for Chebyshev Moment Matching, with Applications”. In: *Proceedings of Thirty Eighth Conference on Learning Theory*. Ed. by N. Haghtalab and A. Moitra. Vol. 291. Proceedings of Machine Learning Research. PMLR, 30 Jun–04 Jul 2025, pp. 4309–4358. arXiv: 2408.12385 [cs.DS] (cited on page 60).

[Neu29] J. von Neumann. "Beweis des Ergodensatzes und des*H*-Theorems in der neuen Mechanik". In: *Zeitschrift für Physik* 57.1-2 (Jan. 1929). English translation https://arxiv.org/abs/1003.2133, pp. 30–70. DOI: 10.1007/bf01339852 (cited on page 54).

[NW83] A. Nauts and R. E. Wyatt. "New Approach to Many-State Quantum Dynamics: The Recursive-Residue-Generation Method". In: *Physical Review Letters* 51.25 (Dec. 1983), pp. 2238–2241. DOI: 10.1103/physrevlett.51.2238 (cited on page 40).

[OLe80] D. P. O'Leary. "The block conjugate gradient algorithm and related methods". In: *Linear Algebra and its Applications* 29 (Feb. 1980), pp. 293–322. ISSN: 0024-3795. DOI: 10.1016/0024-3795(80)90247-5 (cited on pages 66, 71).

[OST07] D. P. O'Leary, Z. Strakoš, and P. Tichý. "On sensitivity of Gauss–Christoffel quadrature". In: *Numerische Mathematik* 107.1 (Apr. 2007), pp. 147–174. DOI: 10.1007/s00211-007-0078-x (cited on page 13).

[OSV12] L. Orecchia, S. Sachdeva, and N. K. Vishnoi. "Approximating the exponential, the Lanczos method and an Õ(m)-time spectral algorithm for balanced separator". In: *Proceedings of the 44th symposium on Theory of Computing - STOC '12*. ACM Press, 2012. DOI: 10.1145/2213977.2214080 (cited on page 40).

[Pai70] C. C. Paige. "Practical use of the symmetric Lanczos process with reorthogonalization". In: *BIT* 10.2 (June 1970), pp. 183–195. DOI: 10.1007/bf01936866 (cited on pages 22, 23).

[Pai72] C. C. Paige. "Computational Variants of the Lanczos Method for the Eigenproblem". In: *IMA Journal of Applied Mathematics* 10.3 (1972), pp. 373–381. DOI: 10.1093/imamat/10.3.373 (cited on pages 22, 23).

[Pai76] C. C. Paige. "Error Analysis of the Lanczos Algorithm for Tridiagonalizing a Symmetric Matrix". In: *IMA Journal of Applied Mathematics* 18.3 (Dec. 1976), pp. 341–349. ISSN: 0272-4960. DOI: 10.1093/imamat/18.3.341 (cited on pages 22, 23).

[Pai80] C. C. Paige. "Accuracy and effectiveness of the Lanczos algorithm for the symmetric eigenproblem". In: *Linear Algebra and its Applications* 34 (1980), pp. 235–258. ISSN: 0024-3795. DOI: 10.1016/0024-3795(80)90167-6 (cited on pages 22, 23).

[Pap19] V. Papyan. *The Full Spectrum of Deepnet Hessians at Scale: Dynamics with SGD Training and Sample Size*. 2019. arXiv: 1811.07062 [cs.LG] (cited on page 58).

[Par98] B. N. Parlett. *The Symmetric Eigenvalue Problem*. Society for Industrial and Applied Mathematics, Jan. 1998. DOI: 10.1137/1.9781611971163 (cited on page 22).

[PCK22] D. Persson, A. Cortinovis, and D. Kressner. "Improved Variants of the Hutch++ Algorithm for Trace Estimation". In: *SIAM Journal on Matrix Analysis and Applications* 43.3 (July 2022), pp. 1162–1185. DOI: 10.1137/21m1447623 (cited on page 56).

[PCM26] D. Persson, T. Chen, and C. Musco. "Randomized block-Krylov subspace methods for low-rank approximation of matrix functions". In: *Linear Algebra and its Applications* 741 (2026), pp. 32–65. ISSN: 0024-3795. DOI: 10.1016/j.laa.2026.03.029 (cited on page 57).

[PK23] D. Persson and D. Kressner. "Randomized Low-Rank Approximation of Monotone Matrix Functions". In: *SIAM Journal on Matrix Analysis and Applications* 44.2 (2023), pp. 894–918. ISSN: 1095-7162. DOI: 10.1137/22m1523923 (cited on page 57).

[PL04] R. Pace and J. P. LeSage. "Chebyshev approximation of log-determinants of spatial weight matrices". In: *Computational Statistics & Data Analysis* 45.2 (Mar. 2004), pp. 179–196. DOI: 10.1016/s0167-9473(02)00321-3 (cited on page 51).

[PL86] T. J. Park and J. C. Light. "Unitary quantum time evolution by iterative Lanczos reduction". In: *The Journal of Chemical Physics* 85.10 (Nov. 1986), pp. 5870–5876. DOI: 10.1063/1.451548 (cited on page 40).

[Ple+20] G. Pleiss, M. Jankowiak, D. Eriksson, A. Damle, and J. Gardner. "Fast Matrix Square Roots with Applications to Gaussian Processes and Bayesian Optimization". In: *Advances in Neural Information Processing Systems*. Ed. by H. Larochelle, M. Ranzato, R. Hadsell, M. F. Balcan, and H. Lin. Vol. 33. 2020, pp. 22268–22281. arXiv: 2006.11267 [cs.LG] (cited on page 40).

[PMM25] D. Persson, R. A. Meyer, and C. Musco. "Algorithm-Agnostic Low-Rank Approximation of Operator Monotone Matrix Functions". In: *SIAM Journal on Matrix Analysis and Applications* 46.1 (Jan. 2025), pp. 1–21. ISSN: 1095-7162. DOI: 10.1137/23m1619435. arXiv: 2311.14023 [math.NA] (cited on page 57).

[Pol09] E. Polizzi. "Density-matrix-based algorithm for solving eigenvalue problems". In: *Physical Review B* 79.11 (Mar. 2009). DOI: 10.1103/physrevb.79.115112 (cited on page 58).

[PS75] C. C. Paige and M. A. Saunders. "Solution of Sparse Indefinite Systems of Linear Equations". In: *SIAM Journal on Numerical Analysis* 12.4 (Sept. 1975), pp. 617–629. DOI: 10.1137/0712047 (cited on page 34).

[PS79] B. N. Parlet and D. S. C. Scott. "The Lanczos algorithm with selective orthogonalization". In: *Mathematics of Computation* 33.145 (Jan. 1979), pp. 217–238. DOI: 10.1090/s0025-5718-1979-0514820-3 (cited on page 23).

[PSW06] S. Popescu, A. J. Short, and A. Winter. "Entanglement and the foundations of statistical mechanics". In: *Nature Physics* 2.11 (Oct. 2006), pp. 754–758. ISSN: 1745-2481. DOI: 10.1038/nphys444 (cited on page 54).

[RA14] F. Roosta-Khorasani and U. Ascher. "Improved Bounds on Sample Size for Implicit Matrix Trace Estimators". In: *Foundations of Computational Mathematics* 15.5 (Sept. 2014), pp. 1187–1212. DOI: 10.1007/s10208-014-9220-1 (cited on page 54).

[Riv20] T. J. Rivlin. *Chebyshev polynomials: from approximation theory to algebra & number theory*. Second edition. Mineola, New York: Dover Publications, Inc, 2020. ISBN: 9780486842332 (cited on page 19).

[Riv81] T. J. Rivlin. *An introduction to the approximation of functions*. Unabridged and corr. republication of the 1969 ed. Dover books on advanced mathematics. Dover, 1981. ISBN: 9780486640693 (cited on page 60).

[RV89] H. D. Raedt and P. de Vries. "Simulation of two and three-dimensional disordered systems: Lifshitz tails and localization properties". In: *Zeitschrift für Physik B Condensed Matter* 77.2 (June 1989), pp. 243–251. DOI: 10.1007/bf01313668 (cited on page 54).

[Saa03] Y. Saad. *Iterative Methods for Sparse Linear Systems*. Society for Industrial and Applied Mathematics, Jan. 2003. ISBN: 9780898718003. DOI: 10.1137/1.9780898718003 (cited on pages 1, 37, 38).

[Saa92] Y. Saad. "Analysis of Some Krylov Subspace Approximations to the Matrix Exponential Operator". In: *SIAM Journal on Numerical Analysis* 29.1 (1992), pp. 209–228. DOI: 10.1137/0729014. eprint: https://doi.org/10.1137/0729014 (cited on pages 40, 42).

[Sch25] M. Schweitzer. "Near Instance Optimality of the Lanczos Method for Stieltjes and Related Matrix Functions". In: *SIAM Journal on Matrix Analysis and Applications* 46.3 (Aug. 2025), pp. 1846–1865. ISSN: 1095-7162. DOI: 10.1137/25m1739650 (cited on page 44).

[Sch27] E. Schrödinger. "Energieaustausch nach der Wellenmechanik". In: *Annalen der Physik* 388.15 (1927), pp. 956–968. DOI: 10.1002/andp.19273881504 (cited on page 54).

[Sil+96] R. Silver, H. Roeder, A. Voter, and J. Kress. "Kernel Polynomial Approximations for Densities of States and Spectral Functions". In: *Journal of Computational Physics* 124.1 (Mar. 1996), pp. 115–130. DOI: 10.1006/jcph.1996.0048 (cited on pages 21, 62).

[Sim84] H. D. Simon. "The Lanczos algorithm with partial reorthogonalization". In: *Mathematics of Computation* 42.165 (1984), pp. 115–142. ISSN: 1088-6842. DOI: 10.1090/s0025-5718-1984-0725988-x (cited on page 23).

[Ski89] J. Skilling. "The Eigenvalues of Mega-dimensional Matrices". In: *Maximum Entropy and Bayesian Methods*. Springer Netherlands, 1989, pp. 455–466. DOI: 10.1007/978-94-015-7860-8_48 (cited on page 62).

[SR94] R. Silver and H. Röder. "Densities of states of mega-dimensional Hamiltonian matrices". In: *International Journal of Modern Physics C* 05.04 (Aug. 1994), pp. 735–753. DOI: 10.1142/s0129183194000842 (cited on page 62).

[SRS20] J. Schnack, J. Richter, and R. Steinigeweg. "Accuracy of the finite-temperature Lanczos method compared to simple typicality-based estimates". In: *Physical Review Research* 2.1 (Feb. 2020). DOI: 10.1103/physrevresearch.2.013186 (cited on page 51).

[SS10] R. Schnalle and J. Schnack. "Calculating the energy spectra of magnetic molecules: application of real- and spin-space symmetries". In: *International Reviews in Physical Chemistry* 29.3 (July 2010), pp. 403–452. DOI: 10.1080/0144235x.2010.485755 (cited on page 51).

[Tre19] L. N. Trefethen. *Approximation Theory and Approximation Practice, Extended Edition*. Society for Industrial and Applied Mathematics, Jan. 2019. DOI: 10.1137/1.9781611975949 (cited on pages 56, 60).

[UCS17] S. Ubaru, J. Chen, and Y. Saad. "Fast Estimation of $\mathrm{tr}(f(A))$ via Stochastic Lanczos Quadrature". In: *SIAM Journal on Matrix Analysis and Applications* 38.4 (Jan. 2017), pp. 1075–1099. DOI: 10.1137/16m1104974 (cited on pages 21, 55, 56).

[Ver18] R. Vershynin. *High-Dimensional Probability*. Cambridge University Press, Sept. 2018. DOI: 10.1017/9781108231596 (cited on pages 54, 55).

[Vil09] C. Villani. *Optimal Transport*. Springer Berlin Heidelberg, 2009. ISBN: 9783540710509. DOI: 10.1007/978-3-540-71050-9 (cited on page 60).

[Vor87] H. V. D. Vorst. "An iterative solution method for solving f(A)x = b, using Krylov subspace information obtained for the symmetric positive definite matrix A". In: *Journal of Computational and Applied Mathematics* 18.2 (May 1987), pp. 249–263. DOI: 10.1016/0377-0427(87)90020-3 (cited on page 40).

[Wan+21] S. Wang, Y. Sun, C. Musco, and Z. Bao. "Public Transport Planning". In: *Proceedings of the 2021 International Conference on Management of Data*. Association for Computing Machinery, June 2021. DOI: 10.1145/3448016.3457247 (cited on page 51).

[Wei+06] A. Weiße, G. Wellein, A. Alvermann, and H. Fehske. "The kernel polynomial method". In: *Reviews of Modern Physics* 78.1 (Mar. 2006), pp. 275–306. DOI: 10.1103/revmodphys.78.275 (cited on pages 21, 51, 58, 62–64).

[WW76] D. Weaire and A. R. Williams. "New numerical approach to the Anderson localization problem". In: *Journal of Physics C: Solid State Physics* 9.17 (Sept. 1976), pp. L461–L463. DOI: 10.1088/0022-3719/9/17/004 (cited on page 54).

[WW77] D. Weaire and A. R. Williams. "The Anderson localization problem. I. A new numerical approach". In: *Journal of Physics C: Solid State Physics* 10.8 (Apr. 1977), pp. 1239–1245. DOI: 10.1088/0022-3719/10/8/025 (cited on page 54).

[Yao+20] Z. Yao, A. Gholami, K. Keutzer, and M. W. Mahoney. "PyHessian: Neural Networks Through the Lens of the Hessian". In: *2020 IEEE International Conference on Big Data (Big Data)*. IEEE, Dec. 2020, pp. 581–590. DOI: 10.1109/bigdata50022.2020.9378171. arXiv: 1912.07145 [cs.LG] (cited on page 58).